\documentclass[11pt]{amsart}

\usepackage{bbm} 
\def\ind{\boldsymbol{\mathbbm{1}}}

\def\R{{\mathbb R}}
\def\N{{\mathbb N}}

\def\Z{{\mathbb Z}}
\def\1{{1\!\!1}}
\def\E{{\mathbb E}}
\def\P{{\mathbb P}}

\def\Q{{\mathbb Q}}
\def\cal{\mathcal}
\def\dom{\rm{dom}\,}

\def\supp{{\rm{supp}}}
\newcommand{\up}{ \tilde{\lambda}}
\newcommand{\down}{\tilde{\eta}}
\newcommand{\domain}{\Lambda_C}
\newcommand{\newg}{\tilde{g}_C}
\newcommand{\newlc}{L_{\delta}}

\newcounter{bean}
\newcommand{\benuma}{\setlength{\labelwidth}{.25in}
\begin{list}%
{(\alph{bean})}{\usecounter{bean}}}
\newcommand{\eenuma}{\end{list}}

\newcommand{\ink}{\rule{.5\baselineskip}{.55\baselineskip}}

\newtheorem{theorem}{Theorem}[section]
\newtheorem{remark}[theorem]{Remark}

\newtheorem{lemma}[theorem]{Lemma}
\newtheorem{cor}[theorem]{Corollary}
\newtheorem{defn}[theorem]{Definition}

\newtheorem{ass}[theorem]{Assumption}

\newcommand{\conv}{\overline{\mbox{\rm co}}}
\newcommand{\spaan}{\mbox{\rm span}}

\newcommand{\sm}{\setminus}

\newcommand{\noi}{\noindent }
\newcommand{\ba}{\begin{array}}
\newcommand{\ea}{\end{array}}
\newcommand{\bea}{\begin{eqnarray}}
\newcommand{\eea}{\end{eqnarray}}
\newcommand{\beas}{\begin{eqnarray*}}
\newcommand{\eeas}{\end{eqnarray*}}
\newcommand{\be}{\begin{equation}}
\newcommand{\ee}{\end{equation}}
\newcommand{\bt}{\begin{theorem}}
\newcommand{\et}{\end{theorem}}
\newcommand{\bc}{\begin{center}}
\newcommand{\ec}{\end{center}}
\newcommand{\ben}{\begin{enumerate}}
\newcommand{\een}{\end{enumerate}}
\newcommand{\lan}{\langle}
\newcommand{\ran}{\rangle}
\newcommand{\ei}{\end{itemize}}
\newcommand{\rear}{\renewcommand{\arraystretch}}

\newcommand{\ds}{\displaystyle}

\newcommand{\ve}{\varepsilon}

\def\npr{\tilde{B}}
\def\newdrift{\tilde{b}}
\def\newfil{\tilde{\cal F}}

\def\newtau{\tilde{\tau}}

\usepackage{psfrag} 
\usepackage{epsfig} 

\newcommand{\raj}{\stackrel{J_1}{\ra}}
\newcommand{\rint}{\mbox{\rm rint}}

 \newcommand{\relint}{\mbox{\rm rint}}

\newcommand{\gap}{\hspace{.15in}} 

\newcommand{\ceil}[1]{{\lceil{#1}\rceil}}

\usepackage{bbm}

\newcommand{\bI}{\mathbb{I}}
\newcommand{\bK}{\mathbb{K}} 
\newcommand{\bJ}{\mathbb{J}}

\newcommand\ddspace{{\cal D}\left[0,\infty\right)} 
\newcommand\ddgspace{{\cal D}_G\left[0,\infty\right)}
\newcommand\dcgspace{{\cal D}_{c,G} \left[0,\infty \right)} 
\newcommand\ddzspace{{\cal D}_0\left[0,\infty\right)}
\newcommand\bvspace{{\cal BV}\left[0,\infty\right)}
\newcommand\bvzspace{{\cal BV}_0\left[0,\infty\right)}
\newcommand\ccspace{{\cal C}\left[0,\infty\right)}
\newcommand\ccgspace{{\cal C}_G \left[0,\infty\right)}

\newcommand{\bGamma}{\bar{\Gamma}}

\newcommand{\cu}{ {\cal U}}

\newcommand{\ra}{\rightarrow}

\newcommand{\suli}{\sum\limits}

\newcommand{\bG}{G}
\newcommand{\cU}{{\cal V}}

\newcommand{\sek}{\sigma_{\delta K}}

\def\lcdc{L_{C,\delta_C}}
\def\kcdc{K_{C,\delta_C}}
\def\mczr{M_C(z,r)}

\def\ol{{\cal L}}

\def\oa{{\cal A}}
\def\oconv{\conv}
\def\cone{\mbox{\rm cone}}

\input{epsf.tex}

\title[Reflected diffusions defined via the extended Skorokhod map]{Reflected diffusions defined via the extended Skorokhod map}
\author{Kavita Ramanan}
\thanks{This research was supported in part by the National Science Foundation Grants 
NSF-DMS-0406191, NSF-DMI-0323668-0000000965 and NSF-DMS-0405343}
\address{Department of Mathematical Sciences\\
Carnegie Mellon University\\
Pittsburgh, PA 15213\\
USA} 
\email{kramanan@math.cmu.edu}


\begin{document}
\begin{abstract}  
This work introduces the extended Skorokhod problem (ESP) and associated extended Skorokhod map (ESM) 
that enable a pathwise construction of reflected diffusions that are not necessarily semimartingales. 
Roughly speaking, given  the closure $G$ of an open connected set in $\R^J$,   
 a non-empty convex cone $d(x) \subset \R^J$  specified at each point $x$ on the boundary $\partial G$,  
and a c\`{a}dl\`{a}g trajectory $\psi$ taking values in 
$\R^J$,  the ESM $\bGamma$ defines a constrained version $\phi$  of $\psi$ that 
takes values in $G$ and is such that the increments of  
$\phi - \psi$ on any interval $[s,t]$ lie in the 
closed convex hull of the directions  $d(\phi(u)), u \in (s,t]$. 
When the  graph of $d(\cdot)$ is closed, the following three properties are established:   
(i)  given $\psi$, if $(\phi,\eta)$ solve the ESP then $(\phi,\eta)$ solve the 
corresponding Skorokhod problem (SP) if and only if $\eta$ is of bounded variation; 
(ii) given $\psi$, any solution $(\phi,\eta)$ to the ESP  is  a solution 
to the SP on the interval $[0,\tau_0)$, but not in general on $[0,\tau_0]$, 
where $\tau_0$ is the first time that $\phi$ hits the set ${\cal V}$ of points $x \in  \partial G$ such that $d(x)$ contains a line; 
(iii)  the graph of the ESM $\bGamma$ is closed on the space of c\`{a}dl\`{a}g trajectories 
(with respect to both the uniform and the $J_1$-Skorokhod topologies). 

The paper then focuses on a class of multi-dimensional ESPs on polyhedral domains with a 
non-empty ${\cal V}$-set. Uniqueness and existence of solutions for this class of ESPs is established and 
existence and pathwise uniqueness  of strong solutions to the associated stochastic 
differential equations with reflection is derived. 
The associated reflected diffusions are also shown to satisfy the corresponding submartingale problem. 
Lastly, it is proved that these reflected diffusions are semimartingales on $[0,\tau_0]$. 
One motivation for the study of this class of reflected diffusions is that they arise as 
approximations of queueing networks in heavy traffic that use the so-called generalised processor sharing discipline. 
\end{abstract} 

\maketitle

\newpage
\bigskip
\hrule
\vspace{-1cm} 
\tableofcontents
\vspace{-1cm}
\hrule

\section{Introduction} 
\label{sec-intro}

\subsection{Background and Motivation} 
\label{subs-desc}

Let $G$ be the closure of an open, connected  domain in  $\R^J$. 
Let $d(\cdot)$ be a set-valued mapping defined on the boundary $\partial G$ of $G$ 
such that for every $x \in \partial G$, $d(x)$ is 
a non-empty, closed and convex cone in $\R^J$ with vertex at the origin $\{0\}$, and the graph 
$\{(x,d(x)): x \in \partial G\}$ of $d(\cdot)$ is closed. 
For convenience, we  extend the definition of $d(\cdot)$ to all of 
$G$ by setting $d(x) = \{0\}$ for $x$ in the interior $G^\circ$ of $G$. 
In this paper we are concerned with reflected deterministic and stochastic processes, and in particular 
reflected Brownian motion, associated with a given pair $(G,d(\cdot))$.   
Loosely speaking, reflected Brownian motion behaves like Brownian motion 
 in the interior $G^\circ$ of the domain $G$ and, whenever it  
reaches a point $x \in \partial G$,  is instantaneously 
restricted to remain in $G$ by a constraining process that pushes along one of the directions in $d(x)$. 
For historical reasons, this constraining action is referred to as instantaneous reflection, and so 
we will refer to $d(\cdot)$ as the reflection field. 
There are three main approaches to the study of 
reflected diffusions -- 
 the Skorokhod Problem (SP) approach, first introduced in \cite{sko} and subsequently 
 developed  in numerous papers such as \cite{andore,cos,dupish1,harrei,lioszn,sai,tan},  
the submartingale problem formulation, introduced in  \cite{strvar}, and Dirichlet form 
methods (used, for example, in \cite{chen1,fukbook}).

In the SP approach, the reflected 
process is represented  as the image of an unconstrained process   
under a deterministic mapping referred to as the Skorokhod Map (SM). 
A rigorous definition of the SP is given below. 
Let $\ddspace$  be the space of $\R^J$-valued, right-continuous functions on $[0,\infty)$ that 
have left limits in $(0,\infty)$.  
Unless stated otherwise, we endow $\ddspace$ with the topology of uniform convergence 
on compact intervals, and note that the resulting space is complete \cite{bil,par}. 
Let $\ddgspace$ (respectively, $\ddzspace$) be the subspace of functions $f$ in 
$\ddspace$ with $f(0) \in G$ (respectively, $f(0) = 0$) and 
 let $\bvzspace$ be the subspace of  functions in $\ddzspace$ that have 
 finite variation on every bounded interval in $[0,\infty)$.  For $\eta \in \bvzspace$ 
and $t \in [0,\infty)$, we use $|\eta|(t)$ to denote the total variation of $\eta$ on $[0,t]$. 
Also, for $x \in G$, let $d^1(x)$ denote the intersection of $d(x)$ with $S_1(0)$, the unit sphere  
in  $\R^J$ centered at the origin. A precise formulation of the SP is given as follows. 

\begin{defn} {\bf (Skorokhod Problem)} 
\label{def-sp}  
Let $(G,d(\cdot))$ and 
$\psi \in \ddgspace$  
be given. 
Then $(\phi, \eta) \in \ddgspace \times \bvzspace$ 
 solve the SP for $\psi$ 
 if $\phi(0) = \psi (0)$, 
 and if for all 
$t \in [0, \infty)$, the following properties are satisfied: 
\begin{enumerate}
\item
$\phi(t) = \psi (t) + \eta (t)$;
\item
$\phi (t) \in \bG$;
\item
$| \eta | (t) < \infty$; 
\item 
$\ds | \eta | (t) = \int_{[0,t]} 
 \ind_{\left\{ \phi (s) \in \partial G \right\} } d | \eta | (s);$ 
\item
There exists a measurable function $\gamma : [ 0, \infty) \rightarrow S_1(0)$ 
such that $\gamma (t) \in d^1( \phi (t))$
($d|\eta |$-almost everywhere)
and
\[ \eta (t) = \int_{[0,t]} \gamma (s) d | \eta | (s). \]
\end{enumerate}     
\end{defn} 
Note that properties 1 and 2 ensure that $\eta$ constrains $\phi$ to remain 
within $G$.  Property 3 requires that the constraining term $\eta$ has finite variation (on every bounded  
interval).  Property 4 allows $\eta$ to change only at times $s$ when $\phi(s)$ is on the 
boundary $\partial G$, in which case property 5 stipulates that the change be along 
one of the directions in $d(\phi(s))$. 
If $(\phi, \phi - \psi)$ solve the SP for $\psi$, 
then we write  $\phi \in \Gamma (\psi)$, 
and refer to $\Gamma$ as the Skorokhod Map (henceforth abbreviated as SM). 
Observe that in general the SM could be multi-valued. 
With some abuse of notation we write $\phi = \Gamma (\psi)$ 
when $\Gamma$ is single-valued and $(\phi, \phi - \psi)$ solve 
the SP for $\psi$.  
The set of $\psi \in \ddgspace$ for which there exists 
 a solution to the SP 
 is defined to be  the domain of the SM $\Gamma$, denoted  $\dom (\Gamma)$.

The SP was first formulated for the case  $G = \R_+$, the non-negative real line,  
and $d(0) = e_1$  by A.V. Skorokhod \cite{sko} in order 
to construct solutions to one-dimensional stochastic differential equations 
with reflection (SDERs), with a Neumann boundary condition at $0$. 
As is well-known (see, for example, \cite{andore}),  the associated one-dimensional SM, which we denote by $\Gamma_1$, 
admits the following explicit representation (here $a \vee b$ denotes the maximum of $a$ and $b$): 
\be
\label{def-1dsmap}
 \Gamma_1 (\psi) (t) \doteq \psi(t) +
\sup_{s \in [0,t]} \left[ -\psi (s) \right] \vee 0. 
\ee 
If $W$ is an adapted, standard Brownian motion defined on a filtered probability space 
$((\Omega, {\cal F}, P), \{{\cal F}_t\})$, 
then the map $\Gamma_1$ can be used to construct a reflected Brownian motion $Z$ by setting 
$Z(\omega) \doteq \Gamma_1 (W(\omega))$ for $\omega \in \Omega$. 
Since the SP is a  pathwise technique,   
it is  especially convenient for establishing  existence and pathwise uniqueness  of 
 strong solutions to SDERs.  
Another advantage of the SP is that, unlike the submartingale problem,
 it can be used to construct reflected stochastic processes 
that are not necessarily diffusions or even Markov processes. 
On the other hand, 
any reflected stochastic process defined as the image of a semimartingale 
under the SM must itself necessarily be a semimartingale (this is an immediate consequence of 
 property 3 of the SP).    Thus the SP formulation does not 
allow the construction of reflected Brownian motions that are not semimartingales.

A second, probabilistic, approach that is used to analyse reflected diffusions is the submartingale problem, which
 was first formulated in \cite{strvar} for the analysis of diffusions on smooth domains with smooth boundary conditions 
and later applied to nonsmooth domains (see, for example, \cite{debtob1,debtob2,kwonwil,taywil,varwil}).  
 The submartingale problem associated with a  class of reflected Brownian motions (RBMs) 
in the $J$-dimensional orthant that are analysed in this paper is described in Definition \ref{def-submg}. 
The submartingale formulation has the advantage that it 
can be used to construct and analyse reflected diffusions that are not necessarily 
semimartingales.  
A drawback, however, is that it only yields weak existence and uniqueness of 
solutions to the associated SDERs. 
The third, Dirichlet form, approach, has an analytic flavor and is particularly
well-suited to the study of symmetric Markov processes (e.g.\ Brownian motion with 
normal reflection) in domains with rough boundaries.  However, once again, this approach 
only yields weak existence and uniqueness of 
solutions \cite{chen1,fukbook}.

In this paper we introduce a fourth approach, which we refer to as the Extended Skorokhod Problem 
(ESP), which enables a pathwise analysis of reflected stochastic processes that are not necessarily 
semimartingales.   As noted earlier, the inapplicability of the SP for the 
construction of 
 non-semimartingale reflected diffusions is a consequence of property 3 of the SP,  which 
requires that the constraining term, $\eta$, be of bounded variation. 
This problem is further compounded by the fact that 
properties 4 and 5 of the SP are also phrased in terms of  the total variation measure $d|\eta|$.  
It is thus natural to ask if  property 3 can be relaxed, while still imposing conditions that suitably 
restrict (in the spirit of properties 4 and 5 of the SP) the times at and directions 
in which $\eta$ can constrain $\phi$. 
This motivates  the following definition.

\begin{defn} {\bf (Extended Skorokhod Problem)} 
\label{def-esp}  
Suppose $(G,d(\cdot))$ and  $\psi \in \ddgspace$ are given.  
Then $(\phi, \eta) \in \ddgspace \times \ddspace$ 
solve the ESP for $\psi$ 
  if $\phi(0) = \psi (0)$,
 and if for all 
$t \in [0, \infty)$, the following properties hold: 
\begin{enumerate}
\item
$\phi(t) = \psi (t) + \eta (t)$;
\item
$\phi (t) \in \bG$;
\item
For every $s \in [0, t]$  
\be
\label{hullprop}
 \eta (t) - \eta(s)  \in \conv \left[ \huge\cup_{u \in (s,t]} d(\phi(u)) \right], 
\ee
where $\conv [A]$ represents the closure of the convex hull
 generated by  the 
set $A$; 
\item 
$\eta (t) - \eta(t-) \in \conv \left[ d(\phi(t)) \right].$ 
\end{enumerate}  
\end{defn}  
Observe that properties 1 and 2 coincide with those of the SP. 
Property 3 is a natural generalisation of property  5 
of the SP when $\eta$ is not necessarily of bounded variation. 
 However, note that it only  guarantees that 
\[ \eta(t) - \eta(t-) \in \conv[d(\phi(t)) \cup d(\phi(t-))]  
\gap \mbox{ for } t \in [0,\infty). \]
In order to ensure uniqueness of solutions under reasonable 
conditions  for paths that 
exhibit jumps,   
it is necessary to impose property 4 as well. 
Since $d(x) = \{0\}$ for $x \in G^\circ$, 
properties 3 and  4 of the ESP together imply that 
 if $\phi(u) \in G^\circ$ for $u \in [s,t]$, 
then $\eta(t) = \eta(s-)$,  
which is a natural generalisation of property 4 of the SP.   
As in the case of  the SP, if $(\phi,\eta)$ solve the ESP for $\psi$, 
we write $\phi \in \bGamma (\psi)$, and refer to 
$\bGamma$ as the Extended Skorokhod Map (ESM), which 
could in general be multi-valued.  
The set of $\psi$ for which the ESP has a solution is 
denoted  $\dom (\bGamma)$. Once again, we will 
abuse notation and write $\phi = \bGamma (\psi)$ when 
$\psi \in \dom(\bGamma)$ and $\bGamma (\psi) = \{\phi\}$ is single-valued.

The first goal of this work is to introduce and prove some general properties of the  ESP, which 
show that the ESP is a natural generalisation of the SP. 
These (deterministic) properties are summarised in Theorem \ref{th-main1}.  
The second objective of this work is to demonstrate the usefulness of the ESP for analysing 
reflected diffusions. 
This is done by focusing on a 
class of reflected diffusions in polyhedral domains in $\R^J$  with piecewise constant reflection 
fields (whose data $(G,d(\cdot))$ satisfy Assumption \ref{as-poly}). 
As shown in \cite{dupram2,dupram4,ramrei1,ramrei2}, 
 ESPs in this class arise as models of queueing  networks  that use 
 the so-called generalised processor sharing (GPS) service discipline. 
For this class of ESPs, existence and pathwise uniqueness of strong solutions to the associated 
SDERs is derived, and the solutions are shown to also  satisfy 
the corresponding submartingale problem. 
In addition, it is shown that the  $J$-dimensional reflected diffusions are semimartingales 
on the closed interval $[0,\tau_0]$, where $\tau_0$ is the 
first time to hit the origin. 
These (stochastic) results are presented in Theorem \ref{th-main2}. 
It was shown in \cite{wil1} that when $J = 2$, RBMs in this class are not semimartingales on $[0,\infty)$. 
In subsequent work, the results derived in this paper are used to study the semimartingale property on $[0,\infty)$ of higher-dimensional 
reflected diffusions in this class.  
The applicability of the ESP to analyse reflected diffusions in curved domains will also 
be investigated in future work.   
In this context, it is worthwhile to note that the ESP coincides with the 
Skorokhod-type lemma introduced in \cite{burtob} for the particular two-dimensional 
thorn domains considered there (see Section \ref{subs-prior} for further discussion). 
The next section provides a more detailed description of the main results.

\subsection{Main Results and Outline of the Paper}
\label{subs-mainres}
 
The first main result characterises deterministic properties of the ESP on  
general domains $G$ with reflection fields $d(\cdot)$ that have a closed graph.  
As mentioned earlier, the space $\ddspace$ is endowed with the 
topology of uniform convergence on compact sets (abbreviated u.o.c.). For notational 
conciseness, throughout the symbol  $\ra$ is used to denote convergence in the u.o.c.\ topology. 
On occasion (in which case this will be explicitly mentioned), we will also consider the 
Skorokhod $J_1$ topology on $\ddspace$ (see, for example,  Section 12.9 of \cite{whibook} for 
a precise definition) and use $\raj$ to denote convergence in this topology. 
Recall $S_1(0)$ is the unit sphere in $\R^J$ centered at the origin. 
The following theorem summarises the main results of Section \ref{sec-esp}.  
Properties  1 and 2 of Theorem \ref{th-main1} correspond to Lemma \ref{lem-esp1}, property 3 is equivalent 
to Theorem \ref{th-esp4} and  property 4 follows from  Lemma \ref{lem-esp3} and Remark \ref{rem-j1}.

\begin{theorem}
\label{th-main1}
Given $(G,d(\cdot))$ that satisfy Assumption \ref{as-domgd}, 
 let $\Gamma$ and $\bGamma$ be the corresponding SM and ESM. 
Then the following properties hold. 
\begin{enumerate} 
\item 
$\dom(\Gamma) \subseteq \dom(\bGamma)$ and for $\psi \in \dom(\Gamma)$, $\phi \in \Gamma(\psi) \Rightarrow \phi \in \bGamma(\psi)$. 
\item 
Suppose 
$(\phi,\eta) \in \ddgspace \times \ddzspace$ 
solve the ESP for  $\psi \in \dom(\bGamma)$. 
Then $(\phi, \eta)$ solve the SP for $\psi$ 
if and only if $\eta \in \bvzspace$. 
\item
If $(\phi,\eta)$ solve the ESP for some $\psi \in \dom(\bGamma)$ 
and $\tau_0 \doteq \inf\{t \geq 0: \phi(t) \in {\cal V}\}$, 
where 
\[ \quad \quad \quad {\cal V} \doteq \{ x \in \partial G: \mbox{ there exists } 
d \in S_1(0) \mbox{ such that } \{d,-d\} \subseteq d(x) \},  
\]
then  $(\phi,\eta)$ also solve the SP for $\psi$ on $[0,\tau_0)$. 
In particular, if ${\cal V} = \emptyset$, then $(\phi,\eta)$ solve the 
SP for $\psi$. 
\item 
Given a sequence of functions $\{\psi_n\}$ such that  
$\psi_n \in \dom (\bGamma)$, for $n \in \N$, and $\psi_n \ra \psi$,  
 let $\{\phi_n\}$ be a corresponding sequence with $\phi_n \in \bGamma (\psi_n)$ for $n \in \N$. 
If there exists a limit point $\phi$ of the sequence $\{\phi_n\}$ with respect to the u.o.c.\ topology, 
 then  $\phi \in \bGamma (\psi)$.  The statement continues to hold 
 if $\psi_n \ra \psi$ is replaced by $\psi_n \raj \psi$ and $\phi$ is a limit point of $\{\phi_n\}$ with respect to the Skorokhod 
$J_1$ topology.  
\end{enumerate} 
\end{theorem} 
The first three results of Theorem \ref{th-main1} demonstrate in what way  the ESM $\bGamma$ is  a generalisation of the 
SM $\Gamma$.  
 In addition, Corollary \ref{cor-cteg} proves that the ESM is in fact a strict generalisation of the 
SM $\Gamma$ for a large class of ESPs with ${\cal V} \neq \emptyset$. 
Specifically, for that class of ESPs it is shown  that  there always exists a continuous function 
$\psi$ and a pair $(\phi,\eta)$ that solve the ESP for 
$\psi$ such that $|\eta|(\tau_0) = \infty$. 
The fourth property of Theorem \ref{th-main1}, stated more succinctly, says that if $d(\cdot)$ has a closed
graph on $\R^J$, then the corresponding (multi-valued) ESM $\bGamma$ also has a closed graph 
(where the closure can be taken with respect to either the u.o.c.\ or Skorokhod $J_1$ topologies).   
As shown in Lemma \ref{lem-esp4}, 
the closure property is  very  useful  for establishing existence of solutions -- the 
corresponding property does not  hold for the SM without the imposition of additional 
conditions on $(G,d(\cdot))$.  For example, 
  the completely-${\cal S}$ condition in \cite{berelk,manvan}, or  generalisations of it 
introduced in \cite{cos} and \cite{dupish1}, were imposed in various contexts 
to establish that the SM $\Gamma$ has a closed graph.  However, all these conditions 
imply that ${\cal V} = \emptyset$.  
Thus properties 3 and 4 above together imply and generalise (see Corollary \ref{cor-esp4} and Remark \ref{rem-oldres}) 
the closure property results for the SM  
established in \cite{berelk,cos,dupish1,manvan}.

While Theorem \ref{th-main1} establishes some very useful properties of the 
ESP under rather weak assumptions on $(G,d(\cdot))$,  additional conditions 
are clearly required to establish existence and uniqueness of solutions to the ESP (an obvious necessary 
condition for existence of solutions is that for each $x \in \partial G$, there exists  
a vector $d \in d(x)$ that points into the interior of $G$).  
Here we do not attempt to derive general conditions for existence and uniqueness of solutions to the 
ESP on arbitrary domains.  Indeed, despite a lot of work on the subject  (see, for example, \cite{andore,berelk,cos,dupish1,dupram3,dupram4,harrei,lioszn,tan}),   
necessary and sufficient conditions for existence and uniqueness of solutions on general 
domains are not fully understood even for the SP. 
Instead,  in Section \ref{subs-dompoly} we focus on a class of  ESPs in polyhedral domains with piecewise constant 
$d(\cdot)$.  
We  establish sufficient conditions for existence and uniqueness of solutions to ESPs in this class in Section \ref{subs-exiuni}, 
and in Theorem \ref{th-exiunigps} verify these conditions for the GPS family of ESPs  described in Section \ref{subs-gpsesp}.   
 This class of ESPs is of interest because it characterises models of networks with fully cooperative 
servers  (see, for instance, \cite{dupram2,dupram4,dupram5,ramrei1,ramrei2}). 
Applications, especially from queueing theory, 
have previously motivated the study of many polyhedral SPs with 
oblique directions of constraint (see, for example, \cite{cheman,daiwil,harrei}).

In Section \ref{sec-sder} we consider SDERs associated with the ESP. 
The next main theorem summarises the results on properties of   reflected diffusions associated with the GPS ESP, 
which has as domain $G = \R_+^J$, the non-negative $J$-dimensional orthant. 
To state these results  we need to first introduce some notation. 
For a given integer $J \geq 2$, 
let $\Omega_J$ be the set of continuous functions $\omega$ from 
$[0,\infty)$ to $\R_+^J = \{x \in \R^J:x_i \geq 0, i = 1, \ldots, J\}$. 
 For $t \geq 0$, let ${\cal M}_t$ be the $\sigma$-algebra 
of subsets of $\Omega_J$ generated by the coordinate maps $\pi_s(\omega) \doteq \omega(s)$ for 
$0 \leq s \leq t$, and let ${\cal M}$ denote the associated $\sigma$-algebra 
$\sigma\{\pi_s:0 \leq s < \infty\}$.  
 The definition of a strong solution to an SDER associated with an ESP is given in Section \ref{subs-SDER}.

\begin{theorem}
\label{th-main2}
Consider 
 drift and dispersion coefficients $b(\cdot)$ and  $\sigma(\cdot)$ that satisfy the usual Lipschitz conditions 
(stated as Assumption \ref{as-ol}(1)) and suppose that a 
$J$-dimensional, adapted Brownian motion, $\{X_t, t \geq 0\}$, defined on a filtered probability space  
$((\Omega, {\cal F}, P), \{{\cal F}_t\})$ is given. Then the following properties hold. 
\begin{enumerate} 
\item
 For every $z \in \R_+^J$, 
there exists a pathwise unique  strong solution $Z$ to the SDER associated with the 
GPS ESP with initial condition $Z(0) = z$. Moreover, $Z$ is a strong Markov process. 
\item 
Suppose, in addition, that the diffusion coefficient is uniformly elliptic (see Assumption \ref{as-ol}(2)). 
 If for each $z \in \R_+^J$, 
 $Q_z$ is the measure induced on $(\Omega_J, {\cal M})$ by the law of the pathwise unique strong solution $Z$ with initial 
condition $z$, then for $J = 2$, 
$\{Q_z, z \in \R_+^J\}$ 
 satisfies the submartingale problem associated with the GPS ESP (described in Definition \ref{def-submg}). 
\item 
Also, if the diffusion coefficient is uniformly elliptic, 
then $Z$ is a 
semimartingale on $[0,\tau_0]$, where $\tau_0$ is the first time to hit the set ${\cal V} = \{0\}$.  
\end{enumerate} 
\end{theorem}  
The first statement of Theorem \ref{th-main2} follows directly from Corollary \ref{cor-exi}, while 
the second property  corresponds to Theorem \ref{th-submg}. 
  As can be seen from the proofs of Theorem \ref{th-exi} and Corollary \ref{cor-exi},  the existence of a strong solution $Z$ 
to the SDER associated with the GPS ESP (under the standard assumptions on the drift and diffusion coefficients) and the 
fact that $Z$ is a semimartingale on $[0,\tau_0)$ are quite straightforward consequences of the 
corresponding deterministic results (specifically, Theorem \ref{th-exiunigps} and Theorem \ref{th-esp4}).  
In turn, these properties can be shown to imply the first two properties of  
the associated submartingale problem.  
The proof of the remaining third condition of the submartingale problem 
 relies on geometric properties of the GPS ESP (stated in Lemma \ref{lem-gsp}) 
that reduce the problem to the verification of a property of  one-dimensional 
reflected Brownian motion, which is carried out in Corollary \ref{lem-1dsp}.

The most challenging result to prove in Theorem \ref{th-main2} is the third property, which is stated as 
Theorem \ref{th-rbm}. 
As Theorem \ref{th-cteg} demonstrates, this result does not carry over 
from a deterministic analysis of the ESP, but instead requires a stochastic analysis.  
In Section \ref{subs-sm}, we first establish this result in a more general setting, namely 
 for strong solutions $Z$ to SDERs associated with general (not necessarily polyhedral) ESPs. 
Specifically, in Theorem \ref{th-smg2} we identify  sufficient conditions (namely  
inequalities (\ref{as-bs3}) and (\ref{eq-fin}) and Assumption \ref{as-gfn}) for the strong solution $Z$ to be 
a semimartingale on $[0,\tau_0]$. 
The first inequality (\ref{as-bs3}) requires that the drift and diffusion coefficients be uniformly bounded in a neighbourhood of 
${\cal V}$. 
This automatically holds for the 
GPS ESP with either bounded or continuous drift and diffusion coefficients since, for the GPS ESP,   
${\cal V}= \{0\}$ is bounded.  
The second  inequality (\ref{eq-fin}) is verified in Corollary \ref{th-estrbm}. 
As shown there, due to a certain relation between $Z$ and an associated 
one-dimensional reflected diffusion (see Corollary \ref{lem-1dsp} for a precise statement) the verification   of the 
relation (\ref{eq-fin}) essentially reduces to checking a property of an ordinary (unconstrained) diffusion. 
 The key condition is therefore Assumption \ref{as-gfn}, which requires the   
 existence of a test function that satisfies certain oblique derivative inequalities on the boundary of the domain. 
Section \ref{sec-appendix} is devoted to the construction of such a test function for (a slight generalisation of) the 
GPS family of ESPs.  This construction may be of independent interest (for example, for the 
construction of viscosity solutions to related partial differential equations \cite{dupish3}).

A short outline of the rest of the paper is as follows. 
In Section \ref{sec-esp}, we derive deterministic properties of the ESP on general 
domains (that satisfy the mild hypothesis stated as Assumption \ref{as-domgd}) -- the main results of 
this section were summarised above in  Theorem \ref{th-main1}.  
  In Section \ref{subs-dompoly}, we specialise to 
the class of so-called polyhedral ESPs (described in Assumption \ref{as-poly}). 
We introduce the class of GPS ESPs in Section \ref{subs-gpsesp} and prove some associated 
properties.   In Section \ref{sec-sder}, we analyse 
SDERs associated with ESPs. We discuss the existence and uniqueness of strong solutions to such 
SDERs in Section \ref{subs-SDER}, and show that the pathwise unique strong solution associated with the GPS ESP 
solves the corresponding submartingale problem in Section \ref{subs-submg}. 
In Section \ref{subs-sm}, we state general sufficient conditions for the reflected diffusion to be a semimartingale 
on $[0,\tau_0]$ and then verify them  for non-degenerate reflected diffusions associated 
with the GPS ESP. 
This entails the construction of certain test functions that satisfy Assumption \ref{as-gfn}, 
the details of which are  relegated to Section \ref{sec-appendix}.

\subsection{Relation to Some Prior Work} 
\label{subs-prior}

When $J=2$, the data $(G,d(\cdot))$ for the polyhedral ESPs studied here corresponds to 
the two-dimensional wedge model of \cite{varwil} with $\alpha = 1$ and the wedge angle less than $\pi$. 
 In \cite{varwil}, the submartingale problem approach was used to 
establish weak existence and uniqueness of the associated reflected Brownian   
motions (RBMs). 
 Corollary \ref{cor-exi} of the present paper (specialised to the case $J = 2$) establishes  strong uniqueness and 
existence of associated reflected diffusions (with drift and diffusion coefficient satisfying the 
usual Lipschitz conditions, and the diffusion coefficient possibly degenerate),  
 thus strengthening  the corresponding result (with $\alpha = 1$)
in Theorem 3.12 of \cite{varwil}. 
The associated RBM  was shown to be a semimartingale on $[0,\tau_0]$ 
in  Theorem 1 of \cite{wil1} and this result, along with additional work, 
was used to show that the RBM is not a semimartingale on $[0,\infty)$ in 
Theorem 5 of \cite{wil1}.  
An explicit semimartingale representation for RBMs in certain two-dimensional 
wedges was also given in \cite{deb90}. 
Here we employ different techniques, that are not restricted to two dimensions, to prove that the 
$J$-dimensional GPS reflected diffusions (for all $J \geq 2$)  are semimartingales 
on $[0,\tau_0]$.  
This result is used in a forthcoming paper 
to study the semimartingale property of this family of $J$-dimensional 
reflected diffusions on $[0,\infty)$.   
Investigation of the semimartingale property is important because 
semimartingales comprise the natural class of integrators for stochastic 
integrals (see, for example, \cite{bic}) and the evolution of functionals of 
semimartingales can be characterized using It\^{o}'s formula. 

Although this paper concentrates on reflected diffusions associated with  the class of GPS ESPs, or more generally 
on ESPs with polyhedral domains 
having piecewise constant reflection fields, as elaborated below, the ESP 
is potentially also useful for analysing non-semimartingale reflected diffusions 
in curved domains.  
In view of this fact, many results in the paper are stated in greater generality 
than required for the  class of polyhedral ESPs that are the focus of this paper. 
Non-semimartingale RBMs in $2$-dimensional cusps with normal 
reflection fields 
were analysed using the submartingale approach 
in \cite{debtob1,debtob2}.  
In \cite{burtob}, a pathwise approach was adopted to examine properties of reflected diffusions in 
downward-pointing $2$-dimensional 
thorns with horizontal vectors of reflection. 
Specifically, the thorns $G$ considered in \cite{burtob}  
admit the following description in terms of  
two continuous real functions $L, R$ defined on  
$[0,\infty)$, with $L(0) = R(0) = 0$ and $L(y) < R(y)$ for all $y > 0$: 
$G = \{(x,y) \in \R^2: y \geq 0, L(y) \leq x \leq R(y)\}$. 
The deterministic Skorokhod-type lemma  introduced  in Theorem 1 of \cite{burtob} 
can easily be seen to correspond  to the ESP associated with $(G,d(\cdot))$, where 
$d(\cdot)$ is specified on the boundary $\partial G$ by 
$d((x,y)) = \{ \alpha e_1, \alpha \geq 0\}$ when $x = L(y), y \neq 0$, 
$d((x,y)) = \{ -\alpha e_1, \alpha \geq 0\}$ when $x = R(y), y \neq 0$,  
$d((0,0)) = \{(x,y) \in \R^2: y \geq 0\}$ and, as usual, $d(x) = \{0\}$ for $x \in G^\circ$. 
The Skorokhod-type lemma of \cite{burtob} can thus be viewed as a particular 
 two-dimensional ESP, and existence and uniqueness for solutions to  this ESP for continuous functions 
$\psi$ (defined on $[0,\infty)$ and taking values in $\R^2$) follows from 
Theorem 1 of \cite{burtob}.  
While the Skorokhod-type lemma of \cite{burtob} was phrased in the context of the two-dimensional thorns 
considered therein, 
the ESP formulation is  applicable to more general reflection fields and  domains in higher dimensions. 
The Skorokhod-type lemma was used in \cite{burtob} to prove an interesting result  
on the boundedness of the variation of the constraining term $\eta$ during a single excursion of 
the reflected diffusions in these thorns.  
Other works  that have  studied the existence of a semimartingale decomposition for 
 symmetric, reflected diffusions associated with Dirichlet spaces on possibly  
non-smooth domains include \cite{chen1,chen2}.

\subsection{Notation}
\label{subs-notat}
 
Here we collect some notation that is commonly 
used throughout the paper. 
Given any subset $E$ of $\R^J$,   
  ${\cal D}([0,\infty):E)$ denotes 
the space of right continuous functions 
with left limits taking values in $E$, and  
  ${\cal BV}([0,\infty):E)$ and 
${\cal C}([0,\infty):E)$, respectively, denote 
 the subspace 
of functions that have bounded variation on every bounded interval and the subspace 
of continuous functions. 
Given $G \subset E \subset \R^J$, 
 ${\cal D}_G([0,\infty):E) \doteq {\cal D}([0,\infty):E) 
\cap \{f \in {\cal D}([0,\infty):E) : f(0) \in G\}$   and 
 ${\cal C}_G([0,\infty):E)$ is defined analogously.  Also, 
${\cal D}_0([0,\infty):E)$ and ${\cal BV}_0([0,\infty):E)$ are defined to be the subspace of functions $f$ that 
satisfy $f(0) = 0$ in ${\cal D}([0,\infty):E)$ and ${\cal BV}([0,\infty):E)$, respectively.  
When $E = \R^J$, for conciseness we  denote these spaces  simply by 
$\ddspace$, $\ddzspace$, $\ddgspace$, $\bvspace$, $\bvzspace$, $\ccspace$ and $\ccgspace$, respectively. 
  Unless specified otherwise, we assume that all the function spaces 
are endowed with the topology of uniform 
convergence (with respect to the Euclidean norm) 
on compact sets, and the notation 
$f_n \ra f$ implies that $f_n$ converges to $f$ in this topology, as $n \ra \infty$.    
For $f \in {\cal BV}([0,\infty):E)$ and $t \in [0,\infty)$,  
let $|f|(t)$ be the 
total variation of $f$ on $[0,t]$ with respect to the 
Euclidean norm on $E$, which is 
denoted by $|\cdot|$. 
For $f \in {\cal D}([0,\infty):E)$ and $t \in [0,\infty)$,    , 
as usual  $f(t-) \doteq  \lim_{s \uparrow t} f(s)$.  
For $U \subseteq \R^J$, we use ${\cal C}(U)$ and ${\cal C}^i(U)$, respectively, to 
denote the space of real-valued functions that are continuous 
and  $i$ times continuously differentiable on some open set 
containing $U$. 
Let $\supp[f]$ represent the support of a real-valued 
function $f$ and 
for $f \in  {\cal C}^1(E)$,   
let $\nabla f$ denote the gradient of $f$.

We use $\bK$ and $\bJ$ to denote the finite sets $\{1, \ldots, K\}$ and $\{1, \ldots, J\}$, 
respectively. 
Given real numbers $a, b$, we let 
$a \wedge b$ and $a \vee b$ denote the 
minimum and maximum of the two numbers respectively.  
For $a \in \R$, as usual $\ceil{a}$ denotes the least  
integer greater than or equal to $a$. 
Given vectors $u, v \in \R^J$, both $\lan u, v \ran$ and 
$u \cdot v$ will be used to denote inner product.  
For a finite set $S$, we use  
$\# [S]$ to denote the cardinality of the 
set $S$.  
For $x \in \R^J$, 
$d(x,A) \doteq \inf_{y \in A} |x - y|$ 
is the Euclidean distance of $x$ from the set $A$. 
Moreover, given $\delta > 0$, we let $N_\delta (A) \doteq \{y \in \R^J: d(y,A) \leq \delta \}$, 
be the closed $\delta$-neighbourhood of $A$. 
With some abuse of notation, when $A = \{x\}$ is a singleton, we
 write $N_\delta(x)$ instead of  $N_\delta(\{x\})$ and    
write $N^\circ(\delta)$ to denote the interior $(N (\delta))^\circ$ of $N_\delta$. 
$S_\delta(x) \doteq \{y \in \R^J: |y-x| = \delta\}$ 
is used to denote the sphere of radius $\delta$ centered at $x$. 
Given any set $A \subset \R^J$ we 
let $A^\circ$, $\overline{A}$ and $\partial A$  
 denote its interior, closure and  boundary  
 respectively, 
 $\ind_{A} (\cdot)$ represents the indicator function of the set $A$,    
$\conv [A]$ denotes   the (closure of the) convex hull generated by the set $A$ and 
$\overline{\cone} [A]$ represents the closure of the non-negative cone $\{ \alpha x: \alpha \geq 0, x \in A\}$  
generated by the set $A$.  
Given sets $A, M \subset \R^J$ with  $A \subset M$,  $A$ is said to be open 
relative to $M$ if $A$ is the intersection of $M$ with some open 
set in $\R^J$. Furthermore, a point $x \in A$ is said to be a relative 
interior point of $A$ with respect to $M$ if there is some 
$\ve > 0$ such that $N_\ve(x) \cap M \subset A$, and the 
collection of all relative interior points is called the 
relative interior of $A$, and denoted as $\rint (A)$.

\section{Properties of the Extended Skorokhod Problem} 
\label{sec-esp} 

   As mentioned in the introduction, throughout the paper we consider  pairs $(G,d(\cdot))$ 
that satisfy the following assumption. 

\begin{ass} {\bf (General Domains)} 
\label{as-domgd}
 $G$ is the closure of a connected,  
 open set in $\R^J$. 
For every $x \in \partial G$, 
$d(x)$ is a non-empty, non-zero, closed, convex cone with vertex at $\{0\}$, 
$d(x) \doteq \{0\}$ for $x \in G^\circ$ and 
the graph $\{(x,d(x)), x \in G\}$ of $d(\cdot)$ is closed. 
\end{ass}

\begin{remark}
\label{rem-domgd}
{\em 
Recall that by definition, the graph of $d(\cdot)$ is closed if and only if for 
every pair of convergent sequences $\{x_n\} \subset G$, $x_n \ra x$ and 
 $\{d_n\} \subset \R^J$, $d_n \ra d$ such that $d_n \in d(x_n)$ for every $n \in \N$,  
 it follows that $d \in d(x)$.  
 Now let 
\be
\label{def-d1} 
d^1(x) \doteq d(x) \cap S_1 (0) \hspace{.1in} \mbox{ for } x \in G 
\ee 
and consider the map $d^1(\cdot):\partial G \rightarrow S_1(0)$. Since 
$\partial G$ and $S_1(0)$ are closed, 
Assumption \ref{as-domgd} implies that the graph of $d^1(\cdot)$ is also closed. 
In turn, since $S_1(0)$ is compact, $d^1(x)$ is compact for every $x \in G$, and so this implies that 
$d^1(\cdot)$ is upper-semicontinuous 
(see Proposition 1.4.8 and Definition 1.4.1 of \cite{aubfra}). 
In other words,  
this means that for every $x \in \partial G$, given $\delta > 0$ there exists $\theta > 0$ such that 
\be 
\label{usc-inc1} 
\cup_{y \in N_\theta(x) \cap \partial G} d^1(y) \subseteq N_\delta (d^1(x)) \cap S_1(0). 
\ee
Since $d(x) = \{0\}$ for $x \in G^\circ$, this implies in fact that given $\delta > 0$, there exists 
$\theta > 0$ such that 
\be
\label{usc-incl2}
\conv\left[ \cup_{y \in N_\theta(x)} d(y) \right] \subseteq \overline{\cone} \left[ N_\delta \left( d^1(x) \right) \right]. 
\ee
In fact, since each $d(x)$ is a  non-empty cone,  the closure of the graph of $d(\cdot)$ 
is  in fact equivalent to the upper semicontinuity (u.s.c.) of $d^1(\cdot)$. 
The latter characterisation will sometimes turn out to be more convenient to use. 
 }
\end{remark}

In this section, we establish some useful (deterministic) properties of the 
ESP under the relatively mild condition stated in Assumption \ref{as-domgd}. 
In Section \ref{subsub-genprop}, we characterise the relationship 
between the SP and the ESP. 
 Section \ref{subsub-vset} 
 introduces the concept of the ${\cal V}$-set, which plays an important role in the 
analysis of the ESP,  and establishes its properties.

\subsection{Relation to the SP} 
\label{subsub-genprop} 

The first result is an elementary 
non-anticipatory property of solutions to the ESP, which holds when the ESM is single-valued. 
A map $\Lambda:\ddspace \ra \ddspace$ will be said to be non-anticipatory if 
for every $\psi, \psi^\prime \in \ddspace$ 
 and $T \in (0,\infty)$,
 $\psi(u) = \psi^\prime (u)$ 
for $u \in [0,T]$ implies that 
$\Lambda(\psi)(u) = \Lambda (\psi^\prime)(u)$ 
for $u \in [0,T]$.

\begin{lemma} {\bf (Non-anticipatory property)} 
\label{lem-esp0}
Suppose $(\phi, \eta)$ solve the ESP $(G, d(\cdot))$ for 
$\psi \in \ddgspace$ and suppose that for $T \in [0,\infty)$,  
\[ \phi^T(\cdot) \doteq \phi(T + \cdot), \gap  
\psi^T (\cdot) \doteq \psi(T + \cdot) - \psi(T), \gap 
\eta^T (\cdot) \doteq \eta(T + \cdot) - \eta(T).  
\]
Then  $(\phi^T, \eta^T)$ solve the ESP for 
$\phi(T) + \psi^T$. 
Moreover, if $(\phi,\eta)$ is the unique solution to the ESP 
for $\psi$ then 
 for any $[T,S] \subset [0,\infty)$, 
$\phi(S)$  depends only on $\phi(T)$ 
and the values $\{\psi(s), s \in [T,S] \}$.  In particular, in this case the ESM and 
the map $\psi \mapsto \eta$ are non-anticipatory. 
\end{lemma}
\begin{proof} 
The proof of the first statement 
follows directly from the definition of the ESP. 
Indeed, 
since $(\phi, \eta)$ solve the ESP for $\psi$, it is clear 
that for any $T < \infty$ and $t \in [0,\infty)$,  $\phi^T(t) - \eta^T(t)$ is equal to  
\[ \phi(T + t) - \eta(T + t) + \eta(T) 
   = \psi(T + t) + \phi(T) - \psi(T) =  \phi(T) + \psi^T(t), \]
which proves property 1 of the ESP.  
Property 2 holds trivially. 
Finally, for any $0 \leq s \leq t < \infty$,  $\eta^T(t) - \eta^T(s)$ is equal to  
\[ \eta(T + t) - \eta(T + s) 
\in \conv \left[ \cup_{u \in (T+s, T+t]} 
d(\phi(u)) \right] = \conv \left[ \cup_{u \in (s, t]} 
d(\phi^T(u)) \right], 
\]
which establishes property 3. Property 4 follows analogously, 
thus proving that $(\phi^T, \eta^T)$ solve the ESP for $\phi(T) + \psi^T$.

If $(\phi,\eta)$ is the unique solution to the ESP for 
$\psi$, then the first statement of the lemma implies that 
for every $T \in [0,\infty)$ and $S > T$, 
$\phi(S) = \bGamma(\phi(T) + \psi^T) (S-T)$ and $\eta(S) = \phi(S) - \psi(S)$. 
This immediately proves  
the second and third assertions of the lemma. 
\end{proof}

The next result 
describes in what sense the ESP is a generalisation of the SP. 
It is not hard to see  from  Definition \ref{def-esp} that 
any solution to the SP is also a solution to the ESP (for the same input 
$\psi$).  
Lemma \ref{lem-esp1} shows in addition that 
solutions  to the ESP for a given $\psi$ 
 are also solutions to the SP for that $\psi$ precisely 
when the corresponding 
constraining term $\eta$  is of finite variation (on bounded intervals).

\begin{lemma} {\bf (Generalisation of the SP)} 
\label{lem-esp1} 
Given data $(G,d(\cdot))$ that satisfies Assumption 
 \ref{as-domgd}, let $\Gamma$ and $\bGamma$, respectively, be the associated SM and ESM. 
 Then  the following properties hold. 
\begin{enumerate} 
 \item 
$\dom(\Gamma) \subseteq \dom(\bGamma)$ and for $\psi \in \dom(\Gamma)$, $\phi \in \Gamma (\psi)  \Rightarrow \phi \in \bGamma(\psi)$. 
 \item
 Suppose 
 $(\phi,\eta) \in \ddgspace \times \ddzspace$ 
solve the ESP for  $\psi \in \dom(\bGamma)$.  
Then  $(\phi, \eta)$ solve the SP for $\psi$ 
 if and only if $\eta \in \bvzspace$. 
 \end{enumerate} 
\end{lemma} 
\begin{proof}
 The first  assertion follows directly from the fact 
that properties 1 and 2 are common to both the SP and the ESP, 
and properties 3-5 in Definition \ref{def-sp} of the SP 
imply properties 3 and 4 in Definition \ref{def-esp} of the ESP.

For the second statement, first let $(\phi,\eta) \in \ddgspace \times \ddzspace$ solve the ESP for 
 $\psi \in \dom(\bGamma)$. 
If $\eta \not \in \bvzspace$, 
then property 3 of the SP is violated, and so 
clearly $(\phi,\eta)$ do not solve the SP for $\psi$. 
Now suppose $\eta \in \bvzspace$. 
Then $(\phi,\eta)$ automatically satisfy properties 1--3 of 
the SP. 
Also observe that $\eta$ is absolutely continuous with respect to $|\eta|$ and let $\gamma$ be the 
Radon-Nikod\`{y}m derivative $d\eta/d|\eta|$ of $d\eta$ with respect to $d|\eta|$. 
Then $\gamma$ is $d|\eta|$-measurable, 
$\gamma(s) \in S_1(0)$ for $d|\eta|$ a.e.\ $s \in [0,\infty)$  
and 
\be
\label{gamma} 
 \eta(t) = \int_{[0,t]} \gamma(s) d|\eta|(s).  
\ee
Moreover, as is well-known  (see, for example, Section X.4 of \cite{doob}),  for 
$d|\eta|$ \ a.e.\ $t \in [0,\infty)$,  
\be
\label{deriv} 
\gamma(t) 
=   \lim_{n \ra \infty} \dfrac{d\eta [t,t+\ve_n]}{d|\eta|[t,t+\ve_n]}
=  \lim_{n \ra \infty} \dfrac{\eta(t+\ve_n) - \eta(t-)}{|\eta|(t+\ve_n) - |\eta|(t-)}, \
\ee
where $\{\ve_n, n \in \N\}$ is a sequence (possibly depending on $t$) such that $|\eta|(t+\ve_n) - |\eta|(t-) > 0$ for every $n \in \N$
and $\ve_n \ra 0$ as 
$n \ra 0$ (such a sequence can always be found for $d|\eta|$ a.e.\ $t \in [0,\infty)$). 
Fix $t \in [0,\infty)$ such that (\ref{deriv}) holds.  Then properties 3 and 4 of the ESP, along 
with the right-continuity of $\phi$, show that given any $\theta > 0$, there exists $\ve_t > 0$ 
such that for every $\ve \in (0,\ve_t)$, 
\be
\label{eta-inc}
\eta(t+\ve) - \eta(t-)  \in \conv \left[ \cup_{u \in [t,t+\ve]} d (\phi(u)) \right] \subseteq \conv \left[ \cup_{y \in N_\theta(\phi(t))} d(y) \right]. 
\ee
If $\phi(t) \in G^\circ$, then since $G^\circ$ is open, there exists $\theta > 0$ such that 
$N_\theta (\phi(t)) \subset G^\circ$, and hence the fact that $d(y) = \{0\}$ for $y \in G^\circ$ 
implies that the right-hand side of (\ref{eta-inc}) is equal to $\{0\}$. 
 When combined with (\ref{deriv}) this implies that 
$\gamma(t) = 0$ for $d|\eta|$ a.e.\ $t$ such that $\phi(t) \in G^\circ$, which establishes property 4 of the SP. 
  On the other hand, if $\phi(t) \in \partial G$ then the u.s.c.\ of $d^1(\cdot)$ (in particular, relation 
(\ref{usc-incl2})) shows that given  $\delta > 0$,  there exists $\theta > 0$ such that 
\be
\label{eta-inc2}
\conv \left[ \cup_{y \in N_\theta(\phi(t))} d(y) \right] \subseteq \overline{\cone} \left[ N_{\delta} \left( d^1(\phi(t)) \right) \right]. 
\ee
Combining this inclusion with (\ref{eta-inc}), (\ref{deriv}) and the fact that $|\eta|(t+ \ve_n)  - |\eta|(t-) > 0$ 
for all $n \in \N$, we conclude that 
\[ \gamma(t) \in \overline{\cone} \left[ N_{\delta} \left( d^1 (\phi(t)) \right) \right] \cap S_1(0). \]
Since $\delta > 0$ is arbitrary, taking the intersection of the right-hand side over $\delta > 0$ 
shows that $\gamma(t) \in d^1(\phi(t))$ for $d|\eta|$ a.e.\ $t$ such that $\phi(t) \in \partial G$. 
Thus  $(\phi,\eta)$ satisfy property 5 of the SP and the proof of the lemma is complete.  
\end{proof}

Lemma \ref{lem-esp3} proves a closure property 
for solutions to the ESP: namely that the graph 
$\{(\psi,\phi): \phi \in \bGamma(\psi), \psi \in \ddgspace\}$  of the set-valued mapping 
$\bGamma$ 
is closed (with respect to both the uniform and Skorokhod $J_1$ topologies).  
As discussed after the statement of Theorem \ref{th-main1}, 
 such a closure property is valid for the SP 
only under certain additional conditions, which are in some 
instances too restrictive (since they imply ${\cal V} = \emptyset$).   
Indeed, one of the goals of this work is to define a suitable pathwise mapping $\psi \mapsto \phi$ for all 
$\psi \in \ddgspace$ even when 
${\cal V} \neq \emptyset$.

\begin{lemma} {\bf (Closure Property)} 
\label{lem-esp3} 
Given an ESP  $(G,d(\cdot))$ that satisfies 
Assumption \ref{as-domgd}, suppose  for $n \in \N$,  
$\psi_n \in \dom (\bGamma)$ and $\phi_n \in \bGamma (\psi_n)$. 
If $\psi_n \ra \psi$ and 
$\phi$ is a limit point (in the u.o.c.\ topology) of the sequence 
$\{\phi_n\}$, then  $\phi \in \bGamma (\psi)$. 
\end{lemma} 

\noi {\bf Remark.} For the class of  
polyhedral ESPs, in Section \ref{subs-exiuni} we 
establish conditions under which the sequence $\{\phi_n\}$ in Lemma \ref{lem-esp3} 
is precompact, so that a limit point $\phi$ exists.  \\

\noi 
{\bf Proof of Lemma \ref{lem-esp3}.}  
Let $\{\psi_n\}$, $\{\phi_n\}$ and $\phi$ be as in the statement of the lemma  
 and set $\eta_n \doteq \phi_n - \psi_n$ and 
$\eta \doteq \phi - \psi$.  
Since $\phi$ is a limit point of $\{\phi_n\}$, 
there must exist a subsequence $\{n_k\}$ such that   
 $\phi_{n_k} \ra \phi$ as $k \ra \infty$. 
Property 1 and (since $G$ is closed) property 2 of the ESP 
are automatically satisfied by $(\phi, \eta)$. 
Now fix $t \in [0,\infty)$.  Then given $\delta > 0$, there 
exists $k_0 < \infty$ such that for all $k \geq k_0$, 
\[ \eta_{n_k}(t) - \eta_{n_k}(t-) \in d(\phi_{n_k}(t)) \subseteq \overline{\cone} \left[ N_{\delta} \left( d^1(\phi(t)) \right) \right], 
\]
where the first relation follows from 
property 4 of the ESP and the second inclusion is a consequence of 
the u.s.c. of $d^1(\cdot)$ (see relation (\ref{usc-incl2})) and  the fact that $\phi_{n_k}(t) \ra \phi(t)$ as $k \ra \infty$.  
Sending first $k \ra \infty$ and then $\delta \ra 0$ in the last display, we conclude that 
\be
\label{eta-jump1}
 \eta(t) - \eta(t-)  \in d(\phi(t)) \quad \mbox{ for every } t \in (0,\infty), 
\ee
which shows that $(\phi, \eta)$ also satisfy property 4 of the ESP for $\psi$. 

Now, let $J_\phi  \doteq \{ t \in (0,\infty): \phi(t) \neq \phi(t-) \}$ be the set of jump points of 
 $\phi$.  Then $J_\phi$ is a closed, countable set  and so $(0,\infty)\sm J_{\phi}$ is open and can 
hence be written as the countable union of open intervals $(s_i,t_i), i \in \N$.  
Fix $i \in \N$ and let $[s,t] \subseteq [s_i,t_i]$.  
  Then for $\ve \in (0, (t - s)/2)$, property 3 of the ESP shows that  
\[ \eta_{n_k} (t - \ve) - \eta_{n_k} (s + \ve) \in \conv \left[ \cup_{u \in (s + \ve, t - \ve]} d\left(\phi_{n_k} (u)\right) \right]. 
\]
We  claim (and justify the claim below) that since $\phi_{n_k} \ra \phi$ and $\phi$ is continuous on $[s + \ve, t - \ve]$, 
given $\delta > 0$ there exists $k_* = k_*(\delta) < \infty$ such that for every $k \geq k_*$, 
\be
\label{lem3-claim}
 \cup_{u \in [s + \ve, t - \ve]} d^1\left(\phi_{n_k} (u) \right) \subseteq N_{\delta} \left( \cup_{u \in [s+\ve, t-\ve]} d^1(\phi(u)) \right) \cap S_1(0). 
\ee
If the claim holds, then the last two displays together show that 
\[ \eta_{n_k} (t - \ve) - \eta_{n_k} (s + \ve) \in \overline{\cone} \left[ \{0\} \cup N_{\delta} \left( \cup_{u \in [s + \ve, t - \ve]} 
d^1 (\phi(u)) \right) \right]. 
\]
Taking limits first as $k \ra \infty$, then $\delta \ra 0$ and lastly $\ve \ra 0$, we obtain 
\[ \eta (t-) - \eta (s) \in \conv \left( \cup_{u \in (s,t)} d(\phi(u)) \right)  \mbox{ if }  s_i \leq s < t \leq t_i \mbox{ for some } i \in \N.\]
Now, for arbitrary $(a,b) \subset (0,\infty)$,  $\eta(b) - \eta(a)$ can be decomposed into a countable sum of terms of the form 
$\eta(t) - \eta(t-)$ for some $t \in (a,b]$ and 
$\eta(t-) - \eta(s)$ for $s,t$ such that  $[s,t] \subseteq [s_i,t_i]$ for some $i \in \N$. 
Thus the last display, 
together with (\ref{eta-jump1}),
 shows that $(\phi,\eta)$ satisfy property 3 of the ESP for $\psi$.

Thus to complete the proof of the lemma, it only remains to justify the claim (\ref{lem3-claim}). 
For this we use an argument by contradiction.  
Suppose there exists some $i \in \N$, $[s,t] \subseteq [s_i,t_i]$, $\ve \in (0,(t-s)/2)$ and $\delta > 0$ such that 
the relation (\ref{lem3-claim}) does not hold. 
Then there exists  a further subsequence of $\{n_k\}$ (which we denote again by $\{n_k\}$),  
and corresponding sequences $\{u_k\}$ and $\{d_k\}$ with  $u_k \in [s+\ve,t-\ve]$, $d_k \in d^1(\phi_{n_{k}}(u_k))$ 
and  $d_k \not \in N_{\delta} \left( \cup_{u \in [s+\ve,t-\ve]} d^1(\phi(u)) \right)$ for $k \in \N$. 
Since 
$S_1(0)$ and $[s+\ve,t-\ve]$ are compact, there exist $d_* \in S_1(0)$ and $u_* \in [s+\ve,t-\ve]$ 
such that $d_k \ra d_*$, $u_k \ra u_*$ (along a common subsequence, which we denote again by $\{d_k\}$ and $\{u_k\}$). 
Moreover, it is clear that 
\be
\label{contra}
 d_* \not \in N_{\delta/2} \left( \cup_{u \in [s+\ve,t-\ve]} d^1(\phi(u)) \right). 
\ee
On the other hand, since $u_k \ra u_*$, $\phi_{n_k} \ra \phi$ (uniformly) on $[s+\ve,t-\ve]$ and $\phi$ is continuous on $(s,t)$, this implies that 
$\phi_{n_{k}} (u_k) \ra \phi(u_*)$.  By  the u.s.c.\ of $d^1(\cdot)$ at $\phi(u_*)$, this means that there exists 
$\tilde{k} < \infty$ such that for every $k \geq \tilde{k}$ the inclusion 
$d^1(\phi_{n_{k}} (u_k)) \subseteq N_{\delta/3} (d^1 (\phi(u_*)))$ is satisfied. 
Since $d_k \in d^1(\phi_{n_{k}} (u_k))$ and $d_k \ra d_*$ this implies $d_* \in N_{\delta/3} (d^1 (\phi(u_*)))$. 
However, 
 this contradicts (\ref{contra}), thus proving the claim (\ref{lem3-claim}) and hence the lemma. 
\ink \\

The three lemmas given above establish general properties of solutions 
to a broad class of ESPs (that satisfy Assumption \ref{as-domgd}), 
assuming that  solutions exist.  
Clearly, additional conditions need to be imposed on $(G,d(\cdot))$ 
 in order to guarantee existence of solutions to the ESP 
(an obvious necessary condition for existence is that for 
every $x \in \partial G$, there exists $d \in d(x)$ that points into the 
interior of $G$). 
In \cite{dupram3} conditions were established for a class of polyhedral ESPs 
(of the form described in Assumption \ref{as-poly})  
that guarantee the existence of solutions for $\psi \in \dcgspace$, 
the space of piecewise constant functions in $\ddgspace$ having a finite number of jumps. 
In the next lemma,  
 the closure property of Lemma \ref{lem-esp3} 
is invoked to show when  existence of solutions to the ESP 
on a dense subset of $\ddgspace$ implies existence on the entire space $\ddgspace$.   
This is used in Section \ref{subs-dompoly} to establish existence and uniqueness of solutions 
to the class of GPS ESPs. 

\begin{lemma} {\bf (Existence and Uniqueness)} 
\label{lem-esp4}
Suppose $(G,d(\cdot))$ is such  that the domain $\dom(\Gamma)$ of the  associated SM $\Gamma$ 
 contains a dense subset ${\cal S}$ of $\ddgspace$ (respectively, $\ccgspace$). 
Then the following properties hold.  
\begin{enumerate} 
\item
If $\Gamma$ is uniformly continuous on ${\cal S}$, 
 then there exists a solution to the ESP for all $\psi \in \ddgspace$ 
(respectively, $\psi \in \ccgspace$). 
\item 
If  $\bGamma$ is uniformly continuous on its domain $\dom(\bGamma)$, then 
$\bGamma$ is defined, single-valued and uniformly continuous on all of 
 $\ddgspace$. 
Moreover, in this case $\psi \in \ccgspace$ implies that 
$\phi = \bGamma(\psi) \in \ccgspace$. 
\end{enumerate} 
In particular, if there exists a projection $\pi: \R^J \ra G$ that satisfies 
\be 
\label{def-dipr}
 \pi(x) = x \quad \mbox{ for } x \in G \hspace{.5in} \mbox{ and } 
\hspace{.5in} 
\pi(x) - x \in d(\pi(x)) \quad \mbox{ for } x \in \partial G,  
\ee
and the ESM is uniformly continuous on its domain, then there exists 
a unique solution to the ESP for all $\psi \in \ddgspace$ and the ESM is uniformly 
continuous on $\ddgspace$.   
\end{lemma}
\begin{proof}
Fix $\psi \in \ddgspace$.  The fact that ${\cal S}$ is dense in 
$\ddgspace$ implies that there exists a sequence 
$\{\psi_n\} \subset {\cal S}$ such that $\psi_n \ra \psi$.   
Since ${\cal S} \subset \dom(\Gamma)$ and  $\Gamma$ is uniformly continuous
on ${\cal S}$, there exists a unique solution to the SP for every $\psi \in {\cal S}$. 
For $n \in \N$, let $\phi_n \doteq \Gamma (\psi_n)$.  
The uniform continuity of $\Gamma$ on ${\cal S}$ along with the
completeness of  $\ddgspace$ with respect to the u.o.c.\ metric  
implies that $\phi_n \ra \phi$ for some $\phi \in \ddgspace$. 
Since $\phi_n = \Gamma(\psi_n)$,  property 1 of Lemma \ref{lem-esp1}  shows that $\phi_n
\in \bGamma(\psi_n)$.  Lemma \ref{lem-esp3} 
then guarantees that $\phi \in \bGamma (\psi)$, from which we  conclude that $\dom(\bGamma) = 
\ddgspace$. This establishes the first statement of the lemma.

Now suppose $\bGamma$ is uniformly continuous on $\dom(\bGamma)$.  
 Then it is automatically single-valued on its domain and so   Lemma \ref{lem-esp1} implies   that  
$\Gamma (\psi) = \bGamma (\psi)$ for $\psi \in \dom(\Gamma)$. 
Thus, by the first statement just proved, we must have $\dom(\bGamma) = \ddgspace$.  
In fact, in this case the proof of the first statement shows that 
$\bGamma$ is equal to the unique uniformly continuous extension of $\Gamma$ from ${\cal S}$ to 
$\ddgspace$ (which exists by p.\ 149 of \cite{roy}). 
In order to prove the second assertion of statement 2 of the lemma,  
fix $\psi \in \ccgspace$ and let $\phi \doteq \bGamma (\psi)$. 
 Since $\phi$ is right-continuous,  it suffices to show that 
for every $\ve > 0$ and $T < \infty$, 
\be
\label{want}
 \lim_{\delta \downarrow 0} \sup_{t \in [\ve,T]} |\phi(t) - \phi(t- \delta)| = 0. 
\ee
Fix $T < \infty$ and $\ve > 0$, and choose $\tilde{\ve} \in (0,\ve)$ and $\delta \in (0,\tilde{\ve})$. 
Define $\psi_1 \doteq \psi, \phi_1 \doteq \phi$, 
\[ \psi_2^{\delta} (s) \doteq  \phi(\tilde{\ve} - \delta) + \psi(s -\delta) - \psi(\tilde{\ve} - \delta) \quad \mbox{ and } \quad  
 \phi_2^{\delta} (s) \doteq \phi (s - \delta)  \quad \mbox{ for } s \in [\ve,T]. \] 
By Lemma \ref{lem-esp0}, it follows that $\phi_2^\delta = \bGamma (\psi_2^\delta)$. Moreover, by the uniform 
continuity of $\bGamma$,  for some function $h:\R_+ \ra \R_+$ such that $h(\eta) \downarrow 0$ 
as $\eta \downarrow 0$, we have for each $\delta \in (0,\tilde{\ve})$, 
\[
\ba{l}
 \sup_{t \in [\ve,T]} |\phi(t) - \phi(t - \delta)|  \\
\quad \quad \quad \quad  \quad \quad =  \sup_{t \in [\ve,T]} |\phi_1(t) - \phi_2^\delta (t)| \\
 \quad \quad \quad \quad \quad \quad\leq  h\left( \sup_{t \in [\ve,T]} |\psi_1(t) - \psi_2^\delta (t)| \right) \\
\quad \quad  \quad \quad \quad \quad = h \left( \sup_{t \in [\ve,T]} |\psi(t) - \psi(t - \delta) + \psi(\tilde{\ve} - \delta) - \phi(\tilde{\ve} - \delta) | \right).  
\ea
 \]
Sending first  $\delta \ra 0$ and then $\tilde{\ve} \ra 0$, and  using the continuity of $\psi$, the right-continuity 
of $\phi$  and the fact that $\phi(0) = \psi(0)$,  we obtain (\ref{want}).

Finally, we use the fact that the existence of a projection is  equivalent to the existence of solutions 
to the SP for $\psi \in \dcgspace$
(see, for example,  \cite{cos,dupish1,dupram3}).  
The last statement of the lemma is then a direct consequence of the first assertion in statement 2 and the fact that 
$\dcgspace$ is dense in $\ddgspace$.  
\end{proof}

\subsection{The ${\cal V}$-set of an ESP} 
\label{subsub-vset} 

In this section we introduce a special set ${\cal V}$ associated with 
an ESP, which plays an important role in characterising the semimartingale 
property of reflected diffusions defined via the ESP (see Section 
\ref{subs-sm}). 
We first establish properties of the set ${\cal V}$ in 
Lemma \ref{lem-flat}.  Then,  
in Theorem \ref{th-esp4},  
we show that solutions to the ESP satisfy the SP until 
the time to hit an arbitrary small neighbourhood of 
${\cal V}$. When ${\cal V} = \emptyset$, this implies 
that any solution to the ESP is in fact a solution to the SP.
A stochastic analogue of this result is presented in 
Theorem \ref{th-exi}.

\begin{defn}
\label{dfn-vset} {\bf (The ${\cal V}$-set of the ESP)}
Given an ESP $(G,d(\cdot))$, 
we define 
\be
\label{def-vset} 
{\cal V} \doteq \{ x \in \partial G: \mbox{ there exists } 
d \in S_1(0) \mbox{ such that } \{d,-d\} \subseteq d^1(x) \}. 
\ee
\end{defn}

Thus the ${\cal V}$-set of the ESP is the set of 
points $x \in \partial G$ such that the set of 
directions of constraint $d(x)$ contains a line. 
Note that $x \in G\sm {\cal V}$ if and only if $d^1(x)$ is  
contained in an open half space of $\R^J$, which is 
equivalent to saying that
\be 
\label{vsetprop} 
 \max_{u \in S_1(0)} \min_{d \in d^1(x)} \lan d,u \ran > 0 \quad \mbox{  for } \quad x  \in G \sm {\cal V}. 
\ee
Following the convention that the minimum over an empty set is infinity, 
the above inequality holds trivially for $x \in G^\circ$. 
 We now prove some useful  properties of the set ${\cal V}$. 
Below $N_\delta^\circ(A)$ denotes the open $\delta$-neighbourhood of the set $A$.

\begin{lemma} 
{\bf (Properties of the ${\cal V}$-set)} 
\label{lem-flat} 
Given an ESP $(G,d(\cdot))$ that satisfies 
Assumption \ref{as-domgd}, let the associated ${\cal V}$-set
 be  defined by (\ref{def-vset}).  
 Then ${\cal V}$ is closed.  
Moreover, given any $\delta > 0$ and 
$L < \infty$ there 
 exist  
 $\rho > 0$ and a finite set $\bK \doteq \{1, \ldots, K\}$ and 
collection  $\{{\cal O}_k, k \in \bK \}$  
of open sets and associated vectors 
 $\{v_k \in S_1(0), k \in \bK \}$   
that satisfy the following two properties. 
\begin{enumerate}
\item
$\left[\{x \in G: |x| \leq L\} \sm N^\circ_\delta ({\cal V}) \right] 
\subseteq 
\left[ \huge\cup_{k \in \bK} {\cal O}_k \right].$ 
\item 
If $y \in \{x \in \partial G: |x| \leq L\}  \cap 
N_\rho ({\cal O}_k)$ for some $k \in \bK $  
then 
\[ \lan d, v_k \ran > \rho 
\hspace{.1in} \mbox{ for every  } \hspace{.1in} 
d \in d^1(y).
\] 
\end{enumerate}   
\end{lemma} 
\begin{proof}
The fact that ${\cal V}$ is closed follows directly from Assumption \ref{as-domgd}, 
which ensures that the graph of $d^1(\cdot)$ is closed (as observed in Remark \ref{rem-domgd}).  
Fix $\delta > 0$ and 
for brevity of notation, define $G_L \doteq \{x \in G: |x| \leq L\}$.  
To establish the second assertion of the theorem, 
we first observe that by (\ref{vsetprop}), for any 
$x \in G\sm {\cal V}$ there exist $\rho_x, \kappa_x > 0$ and 
$v_x \in S_1(0)$ such that 
\be
\label{eq-lbound}
\min_{d \in N_{\kappa_x}(d^1(x))} 
\lan d, v_x \ran > \rho_x.  
\ee
Due to  
the upper semicontinuity of $d^1(\cdot)$ (in particular, property (\ref{usc-incl2})),  given 
any $x \in G_L \sm N_\delta({\cal V})$  
there exists $\ve_x > 0$ such that 
\be
 \label{dyindx} 
y \in N_{3 \ve_x} (x) \hspace{.1in} 
 \Rightarrow \hspace{.1in}  d^1(y) \subseteq N^\circ_{\kappa_x}(d^1(x)). 
\ee
Since 
$G_L \sm N^\circ_\delta({\cal V})$  is compact, 
there exists a finite set $\bK = \{1, \ldots, K\}$ and a finite collection of points 
$\{x_k,k \in \bK\} \subset G_L \sm N^\circ_\delta({\cal V})$ 
such that the corresponding open neighbourhoods ${\cal O}_k \doteq N^\circ_{\ve_{x_k}}(x_k)$, $k \in \bK$, 
form a covering, i.e.,  
\[ G_L \sm N^\circ_\delta ({\cal V}) \subset 
\left[ \cup_{k \in \bK} {\cal O}_k \right]. 
\]
This establishes property 1 of the lemma.

Next, note that  if 
$\tilde{\rho} \doteq \min_{k \in \bK} 
\rho_{x_k}$ and  
$v_k \doteq v_{x_k}$ for $k \in \bK$, by (\ref{eq-lbound}) we obtain the inequality 
\[ 
\min_{d \in N_{\kappa_{x_k}}(d^1(x_k))} 
\lan d, v_k \ran > \tilde{\rho} > 0 
\hspace{.1in} 
\mbox { for 
 } \hspace{.1in} k \in \bK.  
\] 
Let $\rho \doteq \tilde{\rho} \wedge \min_{k \in \bK} \ve_{x_k} > 0$.  
Then combining (\ref{dyindx}) with the last inequality 
we infer that for $k \in \bK$,  
\[ y \in \partial G \cap G_L  
\cap N_\rho ({\cal O}_k)  \Rightarrow 
y \in N_{2 \ve_{x_k}}(x_k) 
\Rightarrow   d^1(y) \subseteq N_{\kappa_{x_k}} (d^1(x_k)), \]
which implies that   $\min_{d \in d^1(y)} \lan d, v_k \ran > \rho$. 
This establishes property 
2 and completes the proof of the lemma.  
\end{proof}

\begin{theorem} 
\label{th-esp4} 
Suppose the ESP $(G, d(\cdot))$ satisfies 
Assumption \ref{as-domgd}.  Let 
 $(\phi, \eta)$ solve the ESP for  $\psi 
\in \ddgspace$ and 
let the associated ${\cal V}$-set be 
as defined in (\ref{def-vset}). 
Then $(\phi, \eta)$ solve the SP on 
$[0,\tau_0)$, where 
\be
\label{def-tau}
 \tau_0 \doteq \inf \{t \geq 0 : \phi(t) \in 
{\cal V} \}. 
\ee
\end{theorem} 
\noi
\begin{proof} 
For $\delta > 0$,  define 
\be 
\label{def-taudel} 
 \tau_\delta \doteq \inf \{t \geq 0 : \phi(t) \in 
 N_\delta({\cal V}) \}. 
\ee
Since $[0,\tau_0) \subseteq \cup_{\delta > 0} [0, \tau_\delta)$  (in fact 
equality holds if $\phi$ is continuous) and 
$(\phi, \eta)$ solve the ESP for $\psi$, 
in order to show that $(\phi, \eta)$ 
solve the SP for $\psi$ on $[0,\tau_0)$, 
by  Lemma \ref{lem-esp1}(2)  it suffices to show that 
\be
\label{finbv}
 |\eta|(T \wedge \tau_\delta-) < \infty  \hspace{.3in} 
\mbox{ for every } \hspace{.1in}  \delta > 0 \gap 
\mbox{ and } \gap T < \infty. 
\ee
Fix $\delta > 0$ and $T < \infty$ and let 
$L \doteq \sup_{t \in [0, T]} |\phi(t)| \vee |\psi(t)|$. 
Note that $L < \infty$ since 
$\phi, \psi \in \ddgspace$, 
and define  $G_L \doteq \{x \in G: |x| \leq L\}$. 
Let $\bK = \{1, \ldots, K\}$, $\rho > 0$, $\{{\cal O}_k, k \in \bK\}$  
and $\{v_k, k \in \bK \}$ 
 satisfy properties 1 and 2 of Lemma \ref{lem-flat}. 
If $\phi(0) \in N_\delta ({\cal V})$, 
$\tau_\delta = 0$ and (\ref{finbv}) follows trivially. 
So assume that $\phi(0) \not \in N_\delta ({\cal V})$, 
which in fact implies that  $\phi(0) \in G_L \sm N_\delta ({\cal V})$. 
Then by Lemma \ref{lem-flat}(1), there exists 
 $k_0 \in \bK$ 
  such that $\phi(0) \in {\cal O}_{k_0}$. 
Let $T_0 \doteq 0$ and consider the sequence $\{T_m,k_m\}$ generated 
recursively as follows. 
For $m = 0, 1, \ldots$, whenever $T_m <  \tau_\delta$,  
 define 
\[ T_{m+1} \doteq \inf \{ t > T_{m}: \phi (t) \not \in 
 N^\circ_\rho ({\cal O}_{k_m}) \mbox{ or } \phi (t) \in 
 N_\delta ({\cal V}) \}.
\]   
If $T_{m+1} < T \wedge \tau_{\delta}$, it follows that 
$\phi(T_{m+1}) \in G_L \sm N_\delta ({\cal V})$ 
and so by Lemma \ref{lem-flat}(1), there exists $k_{m+1}$ such that 
$\phi (T_{m+1}) \in {\cal O}_{k_{m+1}}$. 
Since $\phi \in \ddgspace$ and $\rho > 0$, there 
exists a smallest integer $M < \infty$ such that $T_M \geq T \wedge \tau_\delta$.  
We redefine $T_M \doteq T \wedge \tau_\delta$. 
For $m = 1, \ldots, M$, 
let ${\cal J}_m$ be the jump points of $\eta$ in 
$[T_{m-1},T_m)$ and define 
${\cal J} \doteq \cup_{m=1}^M {\cal J}_m$. 
Given any finite partition $\pi_m = \{T_{m-1}  = t_0^{m} < t_1^m < \ldots t_{j_m}^m 
= T_m\}$ of $[T_{m-1}, T_m]$, we claim (and justify below) that  
\[ 
\ba{rcl} 
4L & \geq &  
 \lan \eta (T_m -) -  
 \eta(T_{m-1} -), v_{k_{m-1}} \ran \\[0.7em]  
&  = & \ds
\left \lan \sum\limits_{i=1}^{j_m} 
\left[\eta(t_i^m  -) - \eta(t_{i-1}^m) \right] 
+  \suli_{t \in {\cal J}_m \cap \pi_m} \left[ \eta(t) - \eta(t-)\right],  
 v_{k_{m-1}}  \right\ran \\[1.0em]
& \geq & 
\ds \rho \left[ \sum\limits_{i=1}^{j_{m}} | \eta(t_{i}^m -) - 
 \eta(t_{i-1}^m) | + \sum_{t \in {\cal J}_m \cap \pi_m} 
|\eta(t) - \eta(t-)| \right].  
\ea
\]
The first inequality above follows from the relation $\eta(t) = \phi(t) - \psi(t)$ and the definition of $L$, while  
 the last inequality uses properties 3 and 4 of the ESP, the fact that 
 $\phi(t) \in N_\rho({\cal O}_{k_{m-1}})$ 
for $t \in [T_{m-1}, T_m)$ and 
  Lemma \ref{lem-flat}(2). 
In turn, this bound implies that 
\[ 
\ba{rcl}
|\eta| (T \wedge \tau_\delta-) & =  & 
\ds \sup_{\pi}
\left[ \suli_{i=1}^{j_\pi} |\eta (t_{i}-) - 
\eta (t_{i-1})| + \sum_{t \in {\cal J} \cap \pi} 
|\eta(t) - \eta(t-)| \right] \\ [0.8em]
& = & \ds \suli_{m=1}^M \ds \sup_{\pi_m}
\left[ 
 \suli_{i=1}^{j_m} 
|\eta (t_i^m-) - \eta (t_{i-1}^m)| + \sum_{t \in {\cal J}_m \cap \pi_m} |\eta(t) - \eta(t-)| \right]\\ [1.2em]
& \leq & \dfrac{4 LM}{\rho},
\ea
\]
where the supremum in the first line is over all finite partitions $\pi = \{0 = t_0 < t_1 < \ldots  < t_{j_\pi} =  T_M\}$  of $[0,T_M]$ and 
the supremum in the second line is over all finite partitions $\pi_m = \{T_{m-1} = t_0^m < \ldots t_1^m < \ldots < t_{j_m}^m = T_m\}$ 
 of $[T_{m-1},T_m]$. 
This establishes (\ref{finbv}) and thus proves the theorem. 
\end{proof}

\begin{cor} 
\label{cor-esp4} 
Suppose the ESP  $(G,d(\cdot))$  satisfies Assumption 
\ref{as-domgd} and has an empty ${\cal V}$-set,  
where ${\cal V}$ is  defined by (\ref{def-vset}).   
If $(\phi,\eta)$ solve the ESP for $\psi \in \ddgspace$, 
then  $(\phi,\eta)$ solve the SP for $\psi$. 
Moreover, if there exists a sequence  
$\{\psi_n\}$ with $\psi_n \ra \psi$ 
such that for every $n \in \N$, 
$(\phi_n, \eta_n)$ solve the SP for 
$\psi_n$, then 
any limit point $(\phi, \eta)$ of 
$\{(\phi_n,\eta_n)\}$ solves the SP for $\psi$. 
\end{cor} 
\begin{proof} 
The first statement follows from Theorem \ref{th-esp4} and 
the fact that $\tau_0 = \infty$ when ${\cal V} = \emptyset$. 
By property 1 of Lemma \ref{lem-esp1}, for every $n \in \N$, 
$(\phi_n, \eta_n)$ also solve the ESP for $\psi_n$. 
Since   $(\phi,\eta)$ is a limit point of  
$\{(\phi_n, \eta_n)\}$,  the closure property of   
Lemma \ref{lem-esp3} then shows that 
 $(\phi,\eta)$ solve the ESP for $\psi$. 
The first statement of the corollary then shows that 
$(\phi, \eta)$ must in fact also solve the SP for $\psi$. 
\end{proof}

\begin{remark}
\label{rem-j1}
{\em  From the definitions of the SP and ESP it is easy to verify that given any 
time change $\lambda$ on $\R_+$ (i.e.,  any continuous, strictly increasing function 
$\lambda:\R_+ \ra \R_+$ with $\lambda(0) = 0$ and $\lim_{t \ra \infty} \lambda(t) = \infty$), 
the pair $(\phi,\eta)$ solve the SP (respectively, ESP) for $\psi \in \ddgspace$ if and only 
if $(\phi \circ \lambda, \eta \circ \lambda)$ solve the SP (respectively, ESP) for 
$\psi \circ \lambda$.  From the definition of the $J_1$-Skorokhod  topology (see, for example, 
Section 12.9 of \cite{whibook}), it then automatically follows that the statements in Lemmas \ref{lem-esp3} and 
 \ref{lem-esp4} and Corollary \ref{cor-esp4} also hold when $\ddgspace$ is endowed with the $J_1$-Skorokhod topology and 
an associated  metric that makes it a complete space, in place of the u.o.c.\ topology and associated metric. }
\end{remark}

\begin{remark} 
\label{rem-oldres}
{\em 
The second statement of Corollary \ref{cor-esp4} 
shows that solutions to the SP are closed under 
limits when ${\cal V} = \emptyset$, 
and thus is a slight generalisation of  related results in 
\cite{berelk}, \cite{cos} and \cite{dupish1}. 
The closure property for polyhedral SPs of the 
form $\{(d_i, e_i, 0), i = 1, \ldots, J\}$ was established in 
\cite{berelk} under what is known as the 
completely-${\cal S}$ condition, 
 which implies the condition ${\cal V} = \emptyset$.  
In Theorem 3.4 of \cite{dupish1},  
this result was generalised to polyhedral SPs 
under a condition (Assumption 3.2 of \cite{dupish1}) that also implies that ${\cal V} = \emptyset$.  
In Theorem 3.1 of \cite{cos}, the closure property was established for 
 general SPs that  satisfy Assumption \ref{as-domgd}, 
have ${\cal V}= \emptyset$ and satisfy the additional conditions   (2.15) and (2.16) of  \cite{cos}.  
The proof given in this paper uses the ESP and is thus different from those given in 
\cite{berelk}, \cite{cos} and \cite{dupish1}. } 
\end{remark}

\section{Polyhedral Extended Skorokhod Problems} 
\label{subs-dompoly} 

In this section, we focus on the class of ESPs $(G,d(\cdot))$ with polyhedral domains 
and piecewise constant reflection fields described in 
 Assumption \ref{as-poly} below.  Henceforth, we refer to this class as 
polyhedral ESPs. 

\begin{ass}
\label{as-poly}
There exists a finite set $\bI = \{1, \ldots, I\}$ and $\{(d_i,n_i,c_i) \in \R^J \times S_1(0) \times \R: 
\lan d_i,n_i\ran  = 1, i \in \bI \}$ such that 
\[ G \doteq \cap_{i\in \bI} \{x \in \R^J: \lan x, n_i \ran \geq c_i \}, \]
$d(x) = \{0\}$ and for $x \in \partial G$, 
\be
\label{def-dx}
 d(x) \doteq \left\{ \sum\limits_{i \in I(x)} \alpha_i d_i: \alpha_i \geq 0 \mbox{ for } i \in I(x)\right\},
\ee
where $I(x) \doteq \{i \in \bI:\lan x, n_i \ran = c_i\}$. 
\end{ass} 
Note that a polyhedral ESP has a polyhedral domain $G$, which consequently admits a 
representation  as the intersection of a finite number of half-spaces. Moreover, it  has a 
reflection field $d(\cdot)$ that is constant and  equal to the ray along the positive 
$d_i$ direction for points in the relative interior of the $(J-1)$-dimensional face 
$\partial G^i \doteq  \{x \in G:\lan x, n_i\ran = c_i\}$ and,   
 at the   intersections of multiple $\partial G^i$, 
 is equal to the convex cone generated by the corresponding $d_i$. 
Thus polyhedral ESPs are completely characterized 
by a finite set of triplets $\{(d_i,n_i,c_i),i \in \bI\}$ (though this representation need not be 
unique).  
The condition $\lan d_i,n_i \ran > 0$ is clearly necessary for existence of
solutions:  the conditions that $n_i \in S_1(0)$ and 
$\lan n_i, d_i \ran = 1$ are just  convenient normalizations.

In Section \ref{subs-exiuni}, we present sufficient conditions for existence and uniqueness 
of solutions to polyhedral ESPs and Lispchitz continuity of the associated ESMs.  
In Section \ref{subs-gpsesp}, we introduce the family of multi-dimensional GPS ESPs, and show in 
Theorem \ref{th-exiunigps} that they satisfy the conditions of Section \ref{subs-exiuni}. 
The ${\cal V}$-set associated with an ESP 
was shown in Section \ref{sec-esp} to play an important role in determining its properties --  particularly in determining its relation 
to the SP. 
In Section \ref{subsub-vsets}, we study the implications of the structure of ${\cal V}$-sets for the properties 
of solutions to polyhedral ESPs.

\subsection{Existence and Uniqueness of Solutions} 
\label{subs-exiuni}

Lemma \ref{lem-esp4} showed that 
the existence of unique solutions for the ESP on all of 
$\ddgspace$ follows if there exists  a projection 
operator for the corresponding SP and the associated ESM is Lispchitz continuous on its domain.  
In this section, we describe sufficient conditions for existence of the projection operator and 
Lipschitz continuity of the ESM  associated with a polyhedral ESP.  
These conditions are verified for   the GPS ESP in Theorem \ref{th-exiunigps}.

Given data $\{(d_i,n_i,c_i),i \in \bI \}$, 
the condition stated below as Assumption \ref{as-setb}  was 
first introduced as Assumption 2.1 of \cite{dupish1}, where it 
was  shown (in Theorem 2.2 of 
\cite{dupish1})  to be sufficient for 
Lipschitz continuity of the SM corresponding to the associated (polyhedral) SP. 
In Theorem \ref{th-lc} below we show that this is also a sufficient 
condition for Lipschitz continuity of the corresponding ESM.  
 Assumption \ref{as-setb} is expressed in terms of the existence of a convex 
set $B$ whose inward normals satisfy certain geometric properties 
dictated by the  data. 
We will see in Section \ref{subs-gfn} that the existence of such a 
set $B$ also plays a crucial role in 
establishing a semimartingale property for 
reflected diffusions in polyhedral domains defined 
via the ESP 
(see the proof of Theorem \ref{lem-locgfn}  given in Section \ref{subs-ap1}). 
A more  easily verifiable  ``dual'' condition that implies 
the existence of a set $B$ satisfying Assumption \ref{as-setb} was introduced in  
 \cite{dupram3, dupram4}. 
As demonstrated in \cite{dupram4,dupram5}, this dual condition is often more convenient 
to use in practice.  

We now state the condition. 
Given a convex set $B \subset \R^J$ and $z \in \partial B$, 
 we let $n(z)$ denote the set of unit  inward normals to the 
set at the point $z$. 
In other words,  
\[ \nu(z) \doteq \{\nu \in S_1(0): 
\lan  \nu, x - z \ran \geq 0 \mbox{ for all } x \in B \}. 
\] 

\begin{ass} {\bf (Set \boldmath B\unboldmath)} 
\label{as-setb} 
There exists a compact, convex, symmetric set
$B$ with $0 \in B^{\circ}$,
and  $\delta > 0$ such that  
for $i \in \bI$, 
\be 
\label{eq-propA}
\left\{ \begin{array}{l} z \in \partial B \\
 | \langle z,n_{i} \rangle | < \delta 
\end{array} \right\} \Rightarrow
\langle \nu, d_{i} \rangle = 0
\mbox{ for all } \nu \in \nu (z).
\ee
\end{ass}

\begin{theorem}
\label{th-lc}
If Assumption \ref{as-setb} is satisfied for the   ESP $\{(d_i,n_i,c_i),i \in \bI\}$, then 
the associated ESM is Lipschitz continuous on its domain of definition. 
\end{theorem}
\begin{proof}  
This proof involves a  straightforward modification of 
the proof of Theorem 2.2 in \cite{dupish1} -- the only difference 
being that here we need to allow for constraining terms  of unbounded variation. 
Thus we only provide arguments that 
 differ from those used in \cite{dupish1} and refer the reader to \cite{dupish1} for the 
remaining details.   
Let $B$ be the  convex set associated with the ESP that satisfies Assumption \ref{as-setb}. 
Suppose for $i = 1, 2$, $\psi_i \in \dom (\bar{\Gamma})$ and 
$(\phi_i, \eta_i)$ solve the ESP for $\psi_i$.  Moreover, let 
$c \doteq \sup_{t \in [0,T]} |\psi_1(t) - \psi_2(t)|$. 
Then we argue by contradiction to show that $\eta_1(t) - \eta_2(t) \in cB$ for all $t \in [0,T]$. 
Suppose there exists $a \in (c,\infty)$ such that $\eta_1(t) - \eta_2(t) \not \in (a B)^\circ$ for 
some $t \in [0,T]$ and let $\tau \doteq \inf\{t \in [0,T]: \eta_1(t) - \eta_2(t) \not \in (aB)^\circ\}$. 
As in \cite{dupish1} we consider two cases. \\
{\em Case 1.}  Suppose $\eta_1 (\tau-) - \eta_2(\tau-) \in \partial (a B)$. In this case, 
let $z \doteq \eta_1(\tau-) - \eta_2(\tau-)$ and  $\nu \in \nu(z/a)$. 
Then for every $t \in (0,\tau)$, the fact that $\eta_1(t) - \eta_2(t) \in (aB)^\circ$, $\nu$ is an inward normal 
to $aB$ at $z$ and $B$ is convex implies that 
\[ \lan z - \eta_1(t) + \eta_2(t), \nu \ran < 0  \quad \mbox{ for every }  t \in (0,\tau). \]
Since $z - \eta_1(t) + \eta_2 (t) = (\eta_1(\tau-) - \eta_1(t)) - (\eta_2(\tau-) - \eta_2(t))$, 
this implies that there must exist a sequence $\{t_n\}$ with $t_n \uparrow \tau$ along which 
at least one of the following relations must hold: 
\[(i) \quad \quad  \lan \eta_1(\tau-) - \eta_1(t_n), \nu \ran < 0 \quad \mbox{ for every } n \in \N; \]
\[ (ii) \quad \quad \lan \eta_2(\tau-) - \eta_2(t_n), \nu \ran > 0 \quad \mbox{ for every } n \in \N. \]
By property 3 of the ESP,  for any $n \in \N$,  
\[ \eta_1(\tau-) - \eta_1(t_n)  \in \oconv \left[\cup_{u \in (t_n, \tau)} d(\phi_1(u)) \right]. 
\]
Also, note that $d(\phi_1(u)) = \{0\}$ if $\phi_1(u) \in G^\circ$ and by (\ref{def-dx}), we have 
\[ d(\phi_1(u)) = \left\{ \sum\limits_{i \in \bI:\lan \phi_1(u),n_i \ran = c_i} \alpha_i d_i: \alpha_i \geq 0 \right\} \quad \mbox{ if } 
\phi_1(u) \in \partial G.  
\]
Therefore,  if (i) holds, there must exist $u_n \in (t_n, \tau)$ and $i_n \in \bI$, 
for $n \in \N$, such that 
$\lan \phi_1(u_n), n_{i_n} \ran = c_{i_n}$ and $\lan d_{i_n}, \nu \ran < 0$. 
This implies that there exists $i \in \bI$ such that 
(possibly along a subsequence, which we relabel again by $n$)   $u_n \in (t_n, \tau)$, 
$\lan \phi_1(u_n), n_i \ran = c_i$ and $\lan d_i, \nu \ran < 0$. 
Taking limits as $n \ra \infty$ and using the fact that $t_n \uparrow \tau$, it follows that 
$\lim_{n \ra \infty} \phi_1(u_n) = \phi_1(\tau-)$ and  $\lan \phi_1(\tau-), n_i \ran = c_i$. 
This establishes relations (19) and (20) in \cite{dupish1}.  The remaining argument given there can then be used to  
arrive at a contradiction for Case 1(i). By symmetry, Case 1(ii) is proved analogously.   \\
{\em Case 2.} Now suppose $\eta_1(\tau-) - \eta_2(\tau-) \in (aB)^\circ$. 
In this case, set $z \doteq \eta_1(\tau) - \eta_2(\tau)$ and let $r \in [a,\infty)$ 
be such that  $\eta_1(\tau) - \eta_2(\tau) \in \partial (r B)$.  Let 
$\nu \in \nu(z/r)$ and note that by the convexity of $B$ we have 
$\lan z - \eta_1(\tau-) + \eta_2(\tau-), \nu \ran < 0.$ 
Noticing that $z - \eta_1(\tau-) + \eta_2(\tau-) = (\eta_1(\tau) - \eta_1(\tau-)) - (\eta_2(\tau) - \eta_2(\tau-))$ 
observe that at least one of the 
inequalities below must hold: 
\[(i) \quad \quad  \lan \eta_1(\tau) - \eta_1(\tau-), \nu \ran < 0 \quad \mbox{ for every } n \in N \] 
\[ (ii) \quad \quad \lan \eta_2(\tau) - \eta_2(\tau-), \nu \ran > 0 \quad \mbox{ for every } n \in N. \]
By property 4 of the ESP, this correspondingly implies that  at least one of the following two relations must be satisfied: 
\[ (i) \mbox{ there exists } i \in \bI \mbox{ such that }  \lan \phi_1(\tau), n_i \ran = c_i  \mbox{ and } 
\lan d_i, \nu \ran < 0 \]
\[ (ii)  \mbox{ there exists } j \in \bI \mbox{ such that }  \lan \phi_2(\tau), n_j \ran = c_j   \mbox{ and } 
\lan d_j, \nu \ran > 0. \]
The rest of the argument leading to a contradiction in Case (ii) now follows as in \cite{dupish1} (from relations (24) and (25) onwards). 
\end{proof}

Sufficient conditions for the existence of projections for SPs 
were derived in \cite{berelk,cos,manvan}.  However, these are not  applicable 
in the present context since they all assume conditions that imply ${\cal V} = \emptyset$. 
So, instead, we refer to the general results on existence of a projection for polyhedral SPs with ${\cal V} \neq \emptyset$ that were 
obtained in Section 4 of \cite{dupram3} (also see \cite{dupram5} and Section 3 of \cite{dupram4} for 
application of these methods to concrete SPs).

\subsection{The GPS Family of ESPs} 
\label{subs-gpsesp}

Generalized processor sharing (GPS) is a service discipline used in high-speed 
networks that allows for the efficient sharing of a single resource amongst 
traffic of different classes.  
The GPS SP was introduced in \cite{dupram2,dupram4} to analyse  
the behaviour of the GPS discipline.  
The GPS  SP admits a representation of the form $\{(d_i,n_i,0),i = 1, \ldots, J+1\}$,  
where $n_i = e_i$ for $i \in \bJ  \doteq \{ 1, \ldots, J\}$ (here $\{e_i, i \in \bJ\}$ is the standard 
orthonormal basis in $\R^J$), $n_{J+1} = \sum_{i=1}^J e_i/\sqrt{J}$, $d_{J+1} = \sum_{i=1}^J e_i/\sqrt{J}$ 
 and the reflection directions $\{d_i, i \in \bJ\}$ 
are defined as follows in terms of a ``weight'' vector $\bar{\alpha} \in (\R^J_+)^\circ$ 
that satisfies $\sum_{i=1}^J \bar{\alpha}_i = 1$: for $i,j \in \bJ$, 
\[ 
(d_{i})_j = 
\left\{
\ba{rl}
-\dfrac{\bar{\alpha}_j}{1 - \bar{\alpha}_i} & \mbox{ for } j \neq i, \\
1 & \mbox{ for } j = i.  
\ea
\right.
\]

We now recall a property of the GPS ESP proved in \cite{dupram4}, and establish 
a useful corollary that will be used in Section \ref{sec-sder} (see Theorem \ref{th-submg} and 
Section \ref{subs-ver1}) to analyse properties of RBMs associated with the GPS ESP. 
 A $J \times J$ matrix $A$ is said to be completely-${\cal S}$ if for every 
principal submatrix $\tilde{A}$ of $A$, there exists a vector $\tilde{y} \geq 0$ such that 
$\tilde A \tilde y > 0$ (here the inequalities hold componentwise).  Completely-${\cal S}$ matrices 
were studied in the context of the Skorokhod Problem and semimartingale reflecting diffusions in 
\cite{berelk,manvan,reiwil}, 

\begin{lemma}
\label{lem-gsp}
The GPS ESP has ${\cal V} = \{0\}$, 
and  the vector $d_{J+1}$ is perpendicular to $\spaan [d(x), x \in \partial G \sm \{0\}]$. 
Moreover, for every  $j \in \bJ$, the $J \times J$ matrix 
$A_j$, whose columns are given by vectors in the set 
$\{d_i, i = 1, \ldots, J+1\} \sm \{d_j\}$, is 
completely-${\cal S}$. 
\end{lemma}
\begin{proof} 
The first statement was proved in Lemma 3.1 of \cite{dupram4} and the 
second statement follows from the proof of Theorem 3.8 of \cite{dupram4}. 
\end{proof}

\begin{cor}
\label{lem-1dsp}
Suppose $(\phi,\eta)$ solve the GPS ESP 
for $\psi \in \ddgspace$. 
Then $\lan \phi, d_{J+1} \ran = \Gamma_1 (\lan \psi, d_{J+1} \ran)$, 
where $\Gamma_1$ is the one-dimensional SM 
 defined in (\ref{def-1dsmap}). 
\end{cor} 
\begin{proof}  
Since $(\phi,\eta)$ solve the GPS ESP for $\psi$, by properties 3  and 4 of the 
ESP we know that 
\be
\label{loc-prop3}
 \eta(t) - \eta(s) \in \oconv \left[ \cup_{u \in (s,t]} d(\phi(u)) \right] \quad 
\mbox{ and } \quad \eta(t) - \eta(t-) \in \oconv[\phi (t)].   
\ee  
By Lemma \ref{lem-gsp}, it follows that 
$d_{J+1}$ is perpendicular to $\spaan [d(x), x \in \partial \R_+^J \sm \{0\}]$. 
Also,   
for every $u \in [0,\infty)$,  
since  $\phi(u) \in \R^J_+$ we have 
$\lan \phi(u), d_{J+1} \ran \geq 0$ and $\lan \phi(u), d_{J+1} \ran = 0$ if and only if $\phi(u) = 0$. 
Define the (set-valued) mapping $\hat{d}$ that takes $\R_+$ to subsets of $\R$ by 
$\hat{d}(x) = \{0\}$ if $x \neq 0$ and $\hat{d}(0) = \R_+$. 
Then taking the inner product of both sides of both terms in (\ref{loc-prop3}) with $d_{J+1}$ we obtain  
\[ \lan \eta(t) - \eta(s), d_{J+1} \ran  \in \oconv \left[ \cup_{u \in (s,t]} \lan d(\phi(u)), d_{J+1} \ran \right] 
= \oconv \left[ \cup_{u \in (s,t]}  \hat{d} \left(\lan \phi(u)), d_{J+1} \ran \right) \right]
\] 
and, similarly, 
 \[  \lan \eta (t) - \eta(t-), d_{J+1}  \ran  \in \lan \conv (d(\phi(t)), d_{J+1} \ran = \hat{d}\left(\lan \phi(t), d_{J+1} \ran\right) 
\]
(where for a set $A \subset \R^J$, we let $\lan A, d_{J+1} \ran $ denote the set 
$\{ \lan x, d_{J+1} \ran : x \in A\}$). 
The above properties imply that $(\lan \phi, d_{J+1} \ran, \lan \eta, d_{J+1} \ran)$ solve the
ESP $(\R_+,\hat{d}(\cdot))$ for $\lan \psi, d_{J+1} \ran$. 
 Moreover, this also shows that 
$\lan \eta, d_{J+1} \ran$ is non-decreasing and so, in particular, lies in 
$\bvzspace$.  It then follows from Lemma \ref{lem-esp1}(2) that 
$(\lan \phi, d_{J+1} \ran, \lan \eta, d_{J+1} \ran)$ solve the one-dimensional SP for 
$\lan \psi, d_{J+1} \ran$.    The lemma then follows from the well-known fact that 
solutions to the one-dimensional SP are unique \cite{sko}. 
\end{proof}

Theorems 3.7 and 3.8 of \cite{dupram4} show that Assumption \ref{as-setb} is satisfied and a projection 
exists for the GPS ESP. 
Combining these  results with Lemma \ref{lem-esp4} and Theorem \ref{th-lc} we obtain 
the following result. 

\begin{theorem} 
\label{th-exiunigps}
For every integer $J \geq 2$ and $\psi \in \ddgspace$, 
 there exists a unique solution $(\phi,\eta)$ 
 associated with the $J$-dimensional GPS ESP.  
Moreover, the associated GPS ESM is Lipschitz continuous on $\ddgspace$. 
\end{theorem}

\begin{remark}
\label{rem-uv}
{\em  Polyhedral data $(G,d(\cdot))$  with more complicated ${\cal V}$-sets  arise 
in the analysis of queueing networks with cooperative servers. 
In particular, in \cite{dupram5} a family of SPs was introduced to analyse 
the fluid model of a two-station queueing network with each station 
serving two classes and in which the GPS discipline is used at each station. 
This SP has domain $\R_+^4$ and an unbounded 
${\cal V}$-set equal to $\{x \in \R^4_+: x_1 = x_4 = 0 \} 
\cup \{x \in \R^4_+: x_2 = x_3 = 0\}$. 
For a certain parameter regime it was shown in Theorems 5 and 10 of \cite{dupram5} that 
 a set $B$ satisfying Assumption \ref{as-setb} and a projection exists for these 
ESPs. As in the GPS case, an application of 
Theorem \ref{th-lc} and Lemma \ref{lem-esp4} then yields 
existence and uniqueness of solutions on $\ddgspace$ for the associated ESPs.} 
 \end{remark}

\subsection{Structure of ${\cal V}$-sets for Polyhedral ESPs} 
\label{subsub-vsets}

In Theorem \ref{th-esp4}, we showed that any 
 solution to an ESP  is also a  solution to the 
associated SP on $[0,\tau_0)$.  
Here we show that this result  
does not in general hold if $[0,\tau_0)$ is replaced 
by $[0,\tau_0]$. 
Specifically, Theorem \ref{th-cteg} 
and Corollary \ref{cor-cteg} below show that 
for a large class of polyhedral ESPs
 with ${\cal V} \neq \emptyset$, 
there exist $(\phi, \eta)$ that solve the ESP 
for  some  $\psi \in \ccgspace$, but do not solve the SP 
on $[0,\tau_0]$ for $\psi$. 
In Section \ref{subs-sm}, we consider the stochastic analogue of this question.

\begin{figure}
\hspace*{-5.6in}
\centerline{\epsfxsize1in}\epsfbox{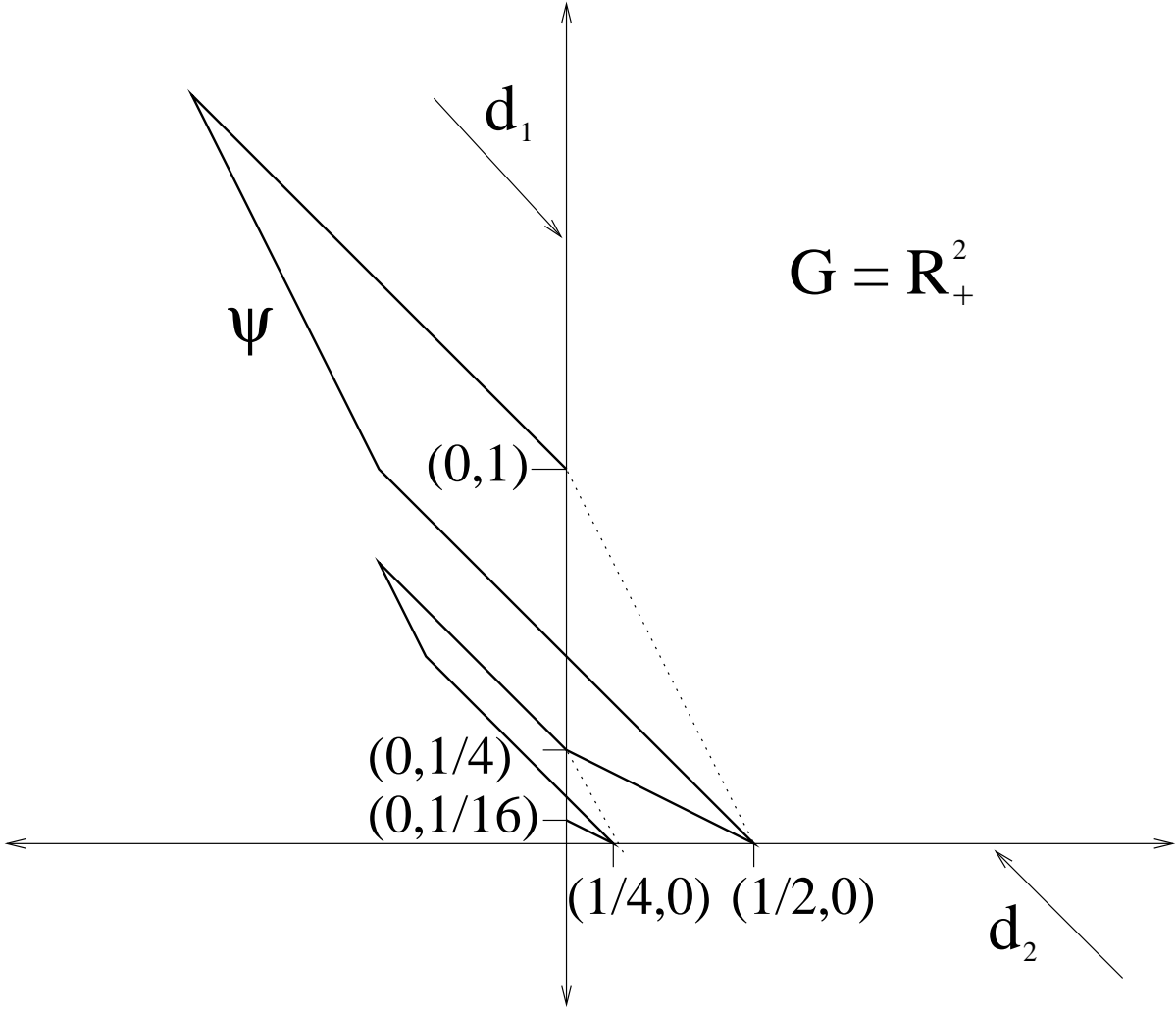} 
\caption{An input trajectory $\psi$ for 
the 2-dimensional GPS ESP}
\label{fig-psi} 
\end{figure}

\begin{theorem}
\label{th-cteg} 
Let $J = 2$ and let $\{(d_i,n_i,0),i = 1,2,3\}$ be the $2$-dimensional 
polyhedral GPS ESP. 
Then there exists $\psi \in \ccgspace$ 
such that $(\phi,\eta)$ solve the ESP for $\psi$ and 
$|\eta|(\tau_0) = \infty$, where $\tau_0 = \inf \{t \geq 0: \phi(t) = 0\}$.  
\end{theorem} 
\begin{proof} Fix $J = 2$. 
From the definition  of the GPS data given in Section \ref{subs-gpsesp}, it is easy 
to see that regardless of the value of the weight vector 
$\bar{\alpha}$,  
the $2$-dimensional GPS ESP has $n_1 = e_1, d_1 = (1,-1)$,  
$n_2 = e_2, d_2 = (-1,1)$ and $n_3 = (1,1)/\sqrt{2}$, $d_3 =  (1,1)/\sqrt{2}$. 
Also, we clearly have ${\cal V} = \{0\}$. 
For notational conciseness,  
let $t_0 \doteq 0$, $\beta \doteq 1/2$, 
and define $t_n \doteq \sum_{i=1}^n \beta^i$. 
Note that then $\lim_{n \ra \infty} t_n = 1$. 
Now, consider the piecewise linear function $\psi \in \ccgspace$ 
defined as follows: $\psi(t_0) = \psi(0) \doteq (0,1)$ and for 
$n \in \N$, recursively define 
\[ 
\ba{rcl}
\psi(t_{n-1}  + \beta^n/4) & = & \psi(t_{n-1}) + (-1,1)/n \\
\psi(t_{n-1}  +  \beta^n/2) & = & \psi(t_{n-1}  + \beta^n/4) + \beta^{2n-2} (\beta,-1) \\
\psi(t_{n-1}  + 3 \beta^n/4) & = & \psi(t_{n-1} + \beta^n/2) + (1,-1)/n \\
\psi(t_{n}) & = & \psi(t_{n-1}  + 3\beta^n/4) + \beta^{2n-1} (-1,\beta),  
\ea
\]
define $\psi(t)$ by linear interpolation for  $t \in [0,1)$ and 
set  $\psi(t) \doteq (0,0)$ for $t \geq 1$ (see Figure \ref{fig-psi} for an illustration of $\psi$).  
It is immediate from the construction that $\psi$ is continuous on $(0,1)$.  
Using the fact that $1 - \beta^2 = 3 \beta^2$ (since $\beta = 1/2$) 
we see that for $i \in \N$, 
\[ \psi(t_i) - \psi(t_{i-1}) = \beta^{2i-2} (\beta, -1) + \beta^{2i-1} (-1, \beta) 
= -3 \beta^{2i} (0,1), 
\] 
and so for $n \in \N$, 
\[ \psi (t_n) = \psi(t_0) + \sum_{i=1}^n \left( \psi(t_{i}) - \psi (t_{i-1}) \right) = 
\left[1 - 3 \sum_{i=1}^n \beta^{2i} \right] (0,1). 
\]
Thus $|\psi(t_n)| = 1/4^n$ 
and $\lim_{n \ra \infty} \psi (t_n) = (0,0) = \psi(1)$ and 
so $\psi \in \ccgspace$. 
Also elementary calculations show that  for $n \in \N$,  
\[ \sup_{t \in [t_{n-1},t_n]} 
|\psi(t) - \psi(t_{n-1})| 
\leq \dfrac{\sqrt{5}}{n}. \]
Hence if we define $k(t)$ to be the unique $k \in \N$ such that 
$t \in [t_{k-1},t_k)$, we deduce that 
\[ \sup_{t \in [t_n,1)} 
|\psi(t) - \psi(t_n)| 
\leq \sup_{t \in [t_n,1)} \left[ |\psi(t) - \psi(t_{k(t)})| 
+ |\psi(t_{k(t)}) - \psi (t_n)| \right] \\ 
\leq  \dfrac{\sqrt{5}}{n} + \dfrac{1}{4^n}.    
\]

\begin{figure}
\hspace*{-4.3in}
\centerline{\epsfxsize1in}\epsfbox{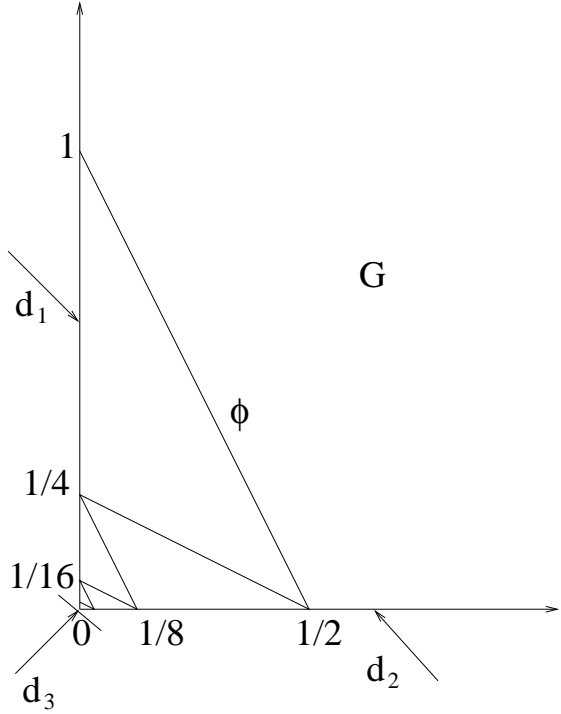} 
\caption{The  solution $(\phi,\eta)$ to 
the 2-dimensional GPS ESP with $|\eta|(\tau_0) = \infty$}
\label{fig-esmeg} 
\end{figure}

Now, let $\phi$ be the piecewise linear trajectory defined as follows: 
$\phi(t_0) = \phi(0) \doteq (0,1)$ and, for $n \in \N$, let  
\[ 
\ba{rcl}
\phi(t_{n-1}  + \beta^n/4) & = & \beta^{2n-2} (0,1),  \\
\phi(t_{n-1}  +  \beta^n/2) & = & \beta^{2n-2} (0,1) + \beta^{2n-2} (\beta,-1) = \beta^{2n-1}(1,0),  \\
\phi(t_{n-1}  + 3 \beta^n/4) & = & \beta^{2n-1} (1,0),  \\
\phi(t_{n}) & = & \beta^{2n-1} (1,0) + \beta^{2n-1} (-1,\beta) = \beta^{2n} (0,1), 
\ea
\]
with $\phi(t)$ defined by linear interpolation for all other 
$t \in (0,1)$ and $\phi(t) \doteq (0,0)$ for $t \geq 1$ (see Figure \ref{fig-esmeg} 
for an illustration of $\phi$).  
It is easy to verify that $\lim_{n \ra \infty} \phi(t_n) = \phi(1) 
= (0,0)$  and so $\phi$ is continuous and $\tau_0 = 1$. 
In addition, it is also straightforward to check that 
 $\eta \doteq \phi - \psi$ is a piecewise linear trajectory that satisfies 
$\eta(0) \doteq (0,0)$ and for $n \in \N$, 
\[ 
\ba{rcl}
\eta(t_{n-1}  + \beta^n/4) & = & \eta(t_{n-1}) + (1,-1)/n,  \\
\eta(t_{n-1}  +  \beta^n/2) & = & \eta(t_{n-1}  + \beta^n/4),  \\
\eta(t_{n-1}  + 3 \beta^n/4) & = & \eta(t_{n-1} + \beta^n/2) + (-1,1)/n, \\
\eta(t_{n}) & = & \eta(t_{n-1}  + 3\beta^n/4),  
\ea
\]
with $\eta(t)$  defined by linear interpolation for all other $t \in (0,1)$ and 
$\eta(t) = (0,0)$ for $t \geq 1$.

Next, define $\psi_n(\cdot) \doteq \psi(\cdot \wedge t_n)$ and, likewise,  
let $\phi_n(\cdot) \doteq \phi(\cdot \wedge t_n)$ 
and $\eta_n(\cdot) \doteq \eta(\cdot \wedge t_n)$. 
From the properties of $\psi_n$ and $\phi_n$ 
stated above, it is clear that $\phi_n = \Gamma (\psi_n)$,  
$\psi_n \ra \psi$, $\phi_n \ra \phi$ and hence 
$\eta_n \ra \eta$.  
Thus by the closure property established in 
Lemma \ref{lem-esp3} 
it follows that $(\phi,\eta)$ satisfy the ESP for 
$\psi$. 
Note that $\tau_0 = 1$ and  
\[ 
\ba{l}
|\eta|(1)\\  
\quad \geq \ds \sum_{n=1}^\infty 
\left[ \left|\eta \left(t_{n-1} + \dfrac{\beta^n}{4}\right) - 
\eta(t_{n-1}) \right| 
+  \left|\eta \left(t_{n-1} +  \dfrac{3\beta^n}{4}\right) - 
\eta\left(t_{n-1}+ \dfrac{\beta^n}{2} \right) \right| \right] \\ 
 \quad \geq  \ds \sum_{n=1}^\infty \dfrac{2 \sqrt 2}{n} \\
 \quad =  \infty, 
\ea
\]
which completes the proof.  
\end{proof}

We now show that the argument in the last proof can 
be generalised to a large class of polyhedral ESPs 
with ${\cal V} \neq \emptyset$.  
Given any ESP, define $\tilde{\cal V}$ to be the set of points $x \in \partial G$ 
such that  there exist $x \in \partial G$, $\rho_0 > 0$ and 
$d \in S_1(0)$ such that for all 
$\rho \in (0,\rho_0]$,  
\be
\label{eq-dminusd} 
\{ d, -d \} \subseteq \overline{\cone} \left[ \cup_{z \in N_\rho(x)\sm \{x\}} d^1(z) \right]. 
\ee 
We then have the following result. 

\begin{cor}
\label{cor-cteg}
Given any polyhedral ESP $\{(d_i,n_i,c_i), i \in \bI \}$, 
if $\tilde{\cal V} \neq \emptyset$,  then there exists 
$\psi \in \ccgspace$ such that $(\phi, \eta)$ solve the ESP for 
$\psi$ and $|\eta|(\tau_0) = \infty$, where $\tau_0 = \inf \{t > 0: \phi(t) \in {\cal V}\}$.  
\end{cor}
\begin{proof}  Let $x \in \partial G$ be such that (\ref{eq-dminusd}) holds. 
Since we are dealing with polyhedral SPs, this implies that there exist 
vectors   $z^{(l)}$ and $v^{(l)}, l = 1, \ldots, L$, such that 
 $\rho z^{(l)}/\rho_0  \in N_\rho (x) \sm \{x\} \cap \partial G$, $v^{(l)} \in d(z^{(l)}) \sm \{0\}$ and 
\[ \sum_{i=1}^L v^{(l)} = 0. \]
Define $z^{(L+1)} \doteq z^{(1)}$. 
Define $\kappa \doteq 1/2L$ and, as in the proof of Theorem 
\ref{th-cteg}, let $\beta \doteq 1/2$, $t_0 \doteq 0$,  
$t_n \doteq \sum_{i=1}^n \beta^i$ and note 
that $\lim_{n \ra \infty} t_n = 1$. 
Now define $\psi \in \ccgspace$ recursively as follows: let 
  $\psi(t_0) = \psi(0) \doteq z^{(1)}$ and 
for $n \in \N$ and $i = 0, \ldots, L-1$,  let
\[ 
\ba{rcl} 
\psi\left(t_{n-1} + (2i+1) \kappa \beta^n\right)  & = & \psi\left(t_{n-1} + 2i \kappa \beta^n\right) - v^{(i)}/n \\
\psi \left(t_{n-1} + (2i+2) \kappa \beta^n \right) & = & \psi \left(t_{n-1} + (2i+1) \kappa \beta^n \right)  \\ 
\quad & &  + \beta^{L(n-1) + i} \left( \beta z^{(i+2)} - z^{(i+1)}\right)
\ea
\]
with $\psi(t)$ defined by linear interpolation for $t \in (0,1)$ and $\psi(t) \doteq 0$ for 
$t \geq 1$. 
It is easy to see that $\lim_{n \ra \infty} \psi(t_n) = \psi(t)$ 
and hence $\psi \in \ccgspace$.   
Moreover, let $\phi(0) \doteq z^{(1)}$, 
for $n \in \N$ and $i = 0, \ldots, L-1$, define 
\[ 
\ba{rcl}
\phi(t_{n-1} + (2i+1) \kappa \beta^n) & = & \phi(t_{n-1} + 2i\kappa \beta^n)   = \beta^{L(n-1)+i}z^{(i+1)} \\
\phi(t_{n-1} + (2i+2) \kappa \beta^n) & = & \phi(t_{n-1} + (2i+ 1)\kappa \beta^n) \\
& & \quad + 
\beta^{L(n-1)+i} \left( \beta z^{(i+2)} - z^{(i+1)} \right), 
\ea
\]
define $\phi(t)$ by linear interpolation for $t \in (0,1)$ and let $\phi(t) \doteq 0$ for 
$t \geq 1$. 
Finally, let $\eta \doteq \phi - \psi$. 
Then arguments analogous to those used in Theorem \ref{th-cteg} 
show that $(\phi, \eta)$ solve the ESP 
for $\psi$, $\tau_0 = 1$ and $|\eta|(\tau_0) = \infty$. 
\end{proof}

\begin{remark} {\em  
It is easy to see that for polyhedral ESPs,  $\tilde{\cal V} \subseteq {\cal V}$. 
Indeed,  if $x \in \partial G$ satisfies (\ref{eq-dminusd}), then 
  by the definition (\ref{def-dx}) of $d(\cdot)$ for polyhedral ESPs, 
there exists $\rho > 0$ such that 
\[  \left[ \cup_{y \in N_\rho(x)\sm \{x\}} d^1(y) \right] 
\subseteq d^1(x). \] 
Along with (\ref{eq-dminusd}), this implies 
that $\{d, -d \} \subseteq d^1(x)$ (and hence that 
$x \in {\cal V}$).

However, there exist polyhedral ESPs for which 
${\cal V} \neq \emptyset$ but  $\tilde{\cal V} = \emptyset$, i.e., (\ref{eq-dminusd}) does 
not hold for any $x \in {\cal V}$. 
For example, consider $\{(d_i,n_i,0), i = 1,2,3\}$, 
where $n_i = e_i$ for $i = 1, 2$, $n_3 = (1,2)/\sqrt{5}$, 
$d_1 = (1,-1), d_2 = (0,1)$ and $d_3 = \sqrt{5}(-1,1)$. 
Then $G = \R_+^2$, $\{(1,-1), (-1,1) \} \subset d(0)$ and 
so  ${\cal V} = \{0\}$.   
On the other hand, for any $\rho > 0$, 
\[ \overline{\cone} [\cup_{y \in N_{\rho}(0)\sm\{0\}} d(y)] = \{\alpha_1 (1,-1)/\sqrt{2} + \alpha_2 (0,1): \alpha_1, \alpha_2 \geq 0 \} 
\]
 and  so $\lan (1,1), d\ran > 0$ for all $d \in \overline{\cone} [\cup_{y \in N_{\rho}(0)\sm\{0\}} d(y)]$. 
Thus (\ref{eq-dminusd}) is not satisfied when $x = 0$.   
However, the ESM is still a strict generalisation of the SM for such ESPs.   
Indeed, given an input trajectory $\psi$ that satisfies $\psi(0) = 0$, is of 
unbounded variation and has an image that lies exclusively on the line $\{ \alpha (1,-1), \alpha \in \R\}$, 
it is clear that $(0,-\psi)$ satisfies the ESP for $\psi$. 
However, $\eta = - \psi$ is of unbounded variation and 
so does not solve the SP for $\psi$.   
Nevertheless, arguments similar to that used in the 
proof of Theorem \ref{th-esp4}  can be used to show that 
 $|\eta|(\tau_0) < \infty$. 
}
\end{remark}

\section{Stochastic Differential Equations with Reflection} 
\label{sec-sder}

In this section, we 
construct and analyse properties of a general  
class of reflected diffusions.      
Throughout, we let 
 $(\Omega, {\cal F}, \P)$ be a complete probability space with 
the right continuous filtration $\{{\cal F}_t, t \geq 0\}$ and 
assume ${\cal F}_0$ contains all  $\P$-negligible sets.   
We will refer to $((\Omega, {\cal F}, \P), \{{\cal F}_t\})$ 
as a filtered probability space. 
In Section \ref{subs-SDER}, we use the ESP   
to construct solutions to 
 a class of stochastic differential 
equations with reflection (SDER).  
In Section \ref{subs-submg}, we show that reflected diffusions 
associated with the GPS family of ESPs 
satisfy the associated submartingale problem. 
When $J=2$, this was proved in \cite{varwil} for the case of 
reflected Brownian motion.

\subsection{Existence and Uniqueness of Solutions to SDERs} 
\label{subs-SDER}

Let $W = \{W_t, t \geq 0\}$ be a $K$-dimensional, $\{{\cal F}_t\}$-adapted, 
standard Brownian motion, and  
consider functions $b(\cdot)$ and $\sigma(\cdot)$   on 
$G$ that take values in 
$\R^J$ and  $\R^J \otimes \R^K$, respectively,  
and satisfy the following standard 
Lipschitz continuity and uniform ellipticity conditions. 
\begin{ass}
\label{as-ol} The functions $b(\cdot)$ and $\sigma(\cdot)$ satisfy the following conditions. 
\begin{enumerate} 
\item 
There exists a constant $\tilde{L} < \infty$ such that for all $x,y \in G$, 
\rear{2.5}
\be
\label{as-bs}
 | \sigma(x) - \sigma(y)| +  |b(x) - b(y)| \leq \tilde{L} |x - y|.  
\ee
\item 
The covariance function $a:G \ra \R^J \otimes \R^J$ defined by 
$a(\cdot) \doteq \sigma (\cdot) \sigma^T (\cdot)$ is uniformly 
elliptic, i.e.,  there exists $\lambda > 0$ such that 
\[ \sum_{i,j=1}^J a_{ij} (x) \xi_i \xi_j \geq \lambda |\xi|^2 \quad \mbox{ for all } \xi \in \R^J \mbox{ and } x \in G. 
\]
\end{enumerate} 
\end{ass}

 We now introduce the 
definition of a strong solution $Z$ to an  SDER associated 
with an ESP $(G,d(\cdot))$, drift coefficient $b(\cdot)$ and dispersion coefficient $\sigma(\cdot)$.  This is  a  
straightforward generalisation of the standard definitions used for 
SDERs defined via  the SP (as, for example,  in Section 5 of \cite{cos}).

\begin{defn} 
\label{def-SDER} 
Given $(G,d(\cdot))$, $b(\cdot)$, $\sigma(\cdot)$ 
and an $\{ {\cal F}_t\}$-adapted $K$-dimensional Brownian motion $W$ on a filtered probability space 
$((\Omega, {\cal F}, \P), \{{\cal F}_t\})$, 
a continuous  $\{{\cal F}_t\}$-adapted process 
$Z(\cdot)$ is a strong solution  to the associated SDER 
if $\P$ a.s.\  for all $t \in [0,\infty)$, 
$Z(t) \in G$ and 
\be
\label{SDER1}
Z(t) = Z(0) + \int_0^t b(Z(s)) \, ds + \int_0^t \sigma(Z(s))\cdot dW(s)   
 + Y(t), 
\ee 
where 
\[ Y(t) - Y(s) \in \conv\left[ \cup_{u \in (s,t]} d(Z(u))\right]. \]
\end{defn}  

\noindent 
In other words, 
$(Z(\cdot),Y(\cdot))$ should solve 
(on a $\P$ a.s.\ pathwise basis) the ESP 
for  
\be
\label{SDER} 
 X(\cdot) \doteq Z(0) + \int_0^\cdot b(Z(s)) \, ds + \int_0^\cdot \sigma(Z(s))
 \cdot dW(s). 
\ee
We will use 
$\P_z$ and $\E_z$ to denote 
probability and expectation, respectively, 
conditioned on $Z(0) = z$.

In the next  theorem,  
we  state sufficient conditions for the existence 
of   a pathwise unique, strong solution to the SDER. 
Since the proof employs     
 standard arguments that are used to construct reflected diffusions  
 defined via the SP (see, for instance,  
\cite{dupkus, tan}) or, more generally, to construct solutions to SDEs 
(as in Section 5.2 of \cite{karshr}),  
we provide only a  rough sketch of the proof. 
We also show that the solution satisfies a semimartingale 
decomposition property -- this is essentially a direct consequence of the corresponding property
for the deterministic ESP that was proved in Theorem \ref{th-esp4}. 
Recall  the definition (\ref{def-vset}) of the ${\cal V}$-set associated with an ESP. 
The theorem below  shows, in particular, that given  a strong solution to the SDER, it 
is a semimartingale if ${\cal V} = \emptyset$.

\begin{theorem}
\label{th-exi}    
Suppose the ESP $(G,d(\cdot))$ satisfies 
Assumption \ref{as-domgd} and 
the associated ESM is well-defined and Lipschitz 
continuous on $\ccgspace$.  
If the coefficients $b(\cdot), \sigma(\cdot)$ satisfy Assumption \ref{as-ol}(1), 
then there exists a pathwise unique, strong  
solution to the associated SDER. 
In addition, the process is strong Markov. 
Furthermore, if we define 
\be
\label{tau}
 \tau_0 \doteq \inf\{t \geq 0: Z(t) \in  {\cal V} \}, 
\ee 
then $Z$ is a semimartingale on  $[0, \tau_0)$, which 
$\P$ a.s.\ admits the decomposition 
\be 
\label{decomp} 
Z (\cdot) = Z(0) + M (\cdot) + A (\cdot),  
\ee
 where for $t \in [0,\tau_0)$, 
\be
\label{MA} 
 M(t) \doteq \int_0^t \sigma(Z(s)) \cdot dW(s)  \quad \mbox{ and } \quad  
A(t) \doteq  \int_0^t b(Z(s)) \, ds  + Y(t), 
\ee
where $Y$ has bounded variation on $[0,t]$ and satisfies  
\be
\label{prop-Y}
 Y(t) = \int_0^t \gamma(s) \, d|Y|(s), 
\ee
and $\gamma(s) \in d^1(Z(s))$ $d|Y|$ a.e.\ $s \in [0,t]$.   
\end{theorem} 
\begin{proof} 
Given that $b(\cdot)$ and 
$\sigma(\cdot)$ satisfy  Assumption \ref{as-ol}(1) 
and the ESM is Lipschitz continuous, 
 we   use the same standard approximations as those  
used to prove 
existence and uniqueness of solutions to SDERs defined via   
a Lipschitz continuous SM (see \cite{andore,tan}). 
For $i = 1, 2$, 
given a continuous $\{{\cal F}_t\}$-adapted process $Z^i$,    
define  $\tilde{X}^{i}$ by 
the right hand side of (\ref{SDER}) with 
$Z$ replaced by $Z^i$, 
and let  
 $X^{i} \doteq \bGamma (\tilde{X}^{i})$. 
Then  the martingale property of the stochastic integral, and 
the Lipschitz continuity of the ESM,  
$b(\cdot)$ and $\sigma(\cdot)$ guarantee  
the existence of $C < \infty$ such that 
\be
\label{be-eq}
  \E \left[\sup_{s \in [0,t]}|X^1 (s) - X^2 (s)|^2  \right]
\leq C \int_{0}^t \E\left[\sup_{r \in [0,s]}|X^1(r) - X^2 (r)|^2\right] ds. 
\ee 
(See Section 5.2 of \cite{karshr} for more details of the 
derivation of this inequality.) 
Applying the usual Picard iteration technique 
along with this bound, 
 and using Gronwall's inequality, \v{C}ebysev's inequality and 
the Borel-Cantelli lemma,  
 one can show the existence of a pathwise unique 
solution to the SDER (following the same arguments as, for example, 
in page 17 of \cite{dupkus}).

Given any $z^1, z^2 \in G$, let $Z^1$ and $Z^2$ be the associated unique strong solutions and 
 for $i = 1,2$, let $\tilde{X}^i$ be equal to the right hand side of 
(\ref{SDER}) with $Z(0)$ replaced by $z^i$. If we also 
 define $X^i \doteq \bGamma (\tilde{X}^i)$, for $i = 1,2$, 
then, just as in  (\ref{be-eq}),
 one can obtain the estimate 
\[ \E[|Z^1 (t) - Z^2 (t)|^2] \leq \tilde{C} |z^1 - z^2|^2. \]
This establishes that any strong solution is a Feller process   
and is therefore strong Markov.

The decomposition in (\ref{decomp})-(\ref{MA}) for the strong solution $Z$ is an immediate consequence of 
Definition \ref{def-SDER}, which 
requires that $\P$ a.s.\ $(Z,Y)$ satisfy the 
ESP for  $X$, 
where $X$ is given by (\ref{SDER}). 
From Theorem \ref{th-esp4} we then conclude 
that $\P$ a.s.\   $(Z,Y)$ satisfy the SP for 
$X$ on $[0,\tau_0)$  and therefore  
$Y$ is an $\{{\cal F}_t\}$-adapted process  of bounded variation on $[0,\tau_0)$ that  
satisfies (\ref{prop-Y}).    
Lastly, since $b(\cdot)$ and $\sigma(\cdot)$  satisfy (\ref{as-bs}), 
 $M$ is a continuous local $\{{\cal F}_t\}$-martingale and 
$\int_{0}^t b(Z_s) ds$ is $\{{\cal F}_t\}$-adapted and of bounded variation. 
This shows that $A$ is $\{{\cal F}_t\}$-adapted and $\P$ a.s.\ of bounded variation,  which 
completes the proof. 
\end{proof}

Combining the above result with Theorem \ref{th-exiunigps} yields the following result for the GPS ESP. 

\begin{cor}
\label{cor-exi} Fix $J \geq 2$ and suppose that 
the coefficients $b(\cdot)$ and  $\sigma(\cdot)$ 
 satisfy Assumption \ref{as-ol}(1). Then  given 
$z \in \R_+^J$, the SDER associated with the GPS ESP has a pathwise unique, strong solution $Z$ with 
initial condition $z$. 
Moreover,  $Z$ is a 
strong Markov process and is a semimartingale on $[0,\tau_0)$ with decomposition (\ref{decomp}). 
\end{cor}

\subsection{The Submartingale Problem and the GPS ESP} 
\label{subs-submg}

Fix an integer $J \geq 2$ and let $\{(d_i,n_i,0),i=1,\ldots,J+1\}$ be the representation for the GPS ESP 
given in Section \ref{subs-gpsesp}. 
Recall that the domain of the GPS ESP is $G = \R_+^J$ and 
for $i = 1, \ldots, J+1$, define $\partial G^i \doteq \{x \in \R^J_+: x_i = 0 \}$ and 
let $S \doteq \cup_{i=1}^J \rint[\partial G^i]$ be the smooth part of the boundary 
$\partial G$. 
As shown in Corollary \ref{cor-exi}, under Assumption \ref{as-ol}(1), 
the SDER associated with the GPS ESP has a unique strong solution. 
In this section we show that this solution also solves the  corresponding submartingale problem.

Recall the definitions of  ${\cal C}^2(G)$ and ${\cal C}^2_b(G)$ given in Section \ref{subs-notat} and, 
given drift and dispersion coefficients $b(\cdot)$ and  $\sigma(\cdot)$, recall  the definition of $a(\cdot)$ stated 
in Assumption \ref{as-ol}. 
 Consider the operator $\ol$ defined by 
\be
\label{def-ol}
\ol f(x) = \sum_{i=1}^J b_i(x) \dfrac{\partial f(x)}{\partial x_i}  + \sum_{i,j=1}^J a_{ij}(x) \dfrac{\partial^2 f(x)}{\partial x_i \partial x_j} 
\quad \mbox{ for } f \in {\cal C}^2(G).  
\ee 
We now define the submartingale problem corresponding to the GPS polyhedral ESP 
 $\{(d_i,n_i,0),i=1,\ldots,J+1\}$ and operator $\ol$ defined above. 
For $J = 2$, $b(\cdot) = 0$ and $a(\cdot) = 1$,  
 this corresponds to the problem analysed in \cite{varwil} with parameter $\alpha = 1$.  
The definition below refers to the canonical filtered probability space $(\Omega_J, {\cal M}, \{{\cal M}_t\})$ that was 
introduced before the statement of Theorem \ref{th-main2}.

\begin{defn} (Submartingale Problem) 
\label{def-submg} 
A family $\{\Q_z, z \in G\}$ of probability measures on $(\Omega_J, {\cal M})$ is a solution 
to the submartingale problem associated with the GPS ESP  $\{(d_i,n_i,c_i), i = 1, \ldots, J+1\}$, 
drift $b(\cdot)$ and dispersion $\sigma(\cdot)$ 
if and only if for each $z \in \R^J_+$,  $\Q_z$ satisfies the following three properties. 
\begin{enumerate} 
\item 
$\Q_z(\omega(0) = z) = 1$;  
\item 
For every $t \in [0,\infty)$ and $f \in {\cal C}_b^2(\R^J_+)$ such that 
$f$ is constant in a neighbourhood of $\partial G\sm S$ and $\lan d_i, \nabla f(x) \ran \geq 0$ for $x \in  
\rint[\partial G^i]$, $i = 1, \ldots, J$, 
\be
\label{eq-submg}
f(\omega(t))  - \int_0^t \ol f(\omega(u)) \, du 
\ee
is a $\Q_z$-submartingale on $(\Omega_J, {\cal M}, \{{\cal M}_t\})$; 
\item
\be
\label{zero}
 \E_{\Q_z} \left[ \int_0^\infty \ind_{\partial G \sm S} (\omega(s)) \, ds \right] = 0. 
\ee
\end{enumerate} 
In this case, $\Q_z$ is said to be the solution of the submartingale problem starting at $z$. 
\end{defn}

The following theorem shows that the family of laws induced by the unique strong solutions of the 
GPS ESP satisfy the associated submartingale problem.

\begin{theorem}
\label{th-submg} Fix $J = 2$. 
Given drift and dispersion coefficients $b(\cdot)$ and $\sigma(\cdot)$
 that satisfy Assumption \ref{as-ol}
and an adapted $K$-dimensional Brownian motion $W$ defined on a filtered probability space 
$(\Omega, {\cal F}, \{{\cal F}_t\}, \P)$, for $z \in \R^J_+$ let 
$\Q_z$ be the measure  on $(\Omega_J, {\cal M}, \{{\cal M}_t\})$ induced  by the unique strong solution 
$Z$  to the SDER associated with the GPS ESP that has initial condition $z$. 
Then $\{\Q_z,z \in \R_+^J\}$ satisfies the GPS submartingale problem. 
\end{theorem}
\begin{proof}
Fix $z \in \R_+^J$ and let $Z$ be the pathwise unique, strong solution associated with the GPS ESP that 
has initial condition $z$ (which exists by Corollary \ref{cor-exi}).  
By definition,  $\Q_z(\omega(0) = z) = \P_z (Z(0) = z) = 1$ and so the first property of 
Definition \ref{def-submg} is trivially satisfied.  

  Now, let $f \in {\cal C}_b^2(\R^J_+)$ be as in the statement of property 2 of the submartingale 
problem stated in Definition \ref{def-submg}, and 
fix  $\ve > 0$ such that $f(x) = 0$ for $x \in N_\ve(0) \cap G$. 
Define $\theta_0 \doteq 0$ and for $n \in \N$, let 
\[ \sigma_n \doteq \inf \{t > \theta_{n-1}: |Z(t)| \leq \ve/2 \}  \quad \mbox{ and } \quad \theta_n \doteq 
\inf \{t > \sigma_n: |Z(t)| \geq \ve  \} 
\] 
where, by the usual convention, the infimum over an empty set is taken to be $\infty$. 
Since $[0,\ve/2]$ and $[\ve,\infty)$ are closed sets, $\sigma_n$ and $\theta_n$ are $\{{\cal F}_t\}$-stopping times. 
Consider the case when $\sigma_n, \theta_n \rightarrow \infty$ as $n \ra \infty$ (the other case can be dealt with in a similar 
manner and is thus left to the reader). 
Then for $t \in [0,\infty)$, we can write 
\[ 
 \ba{rcl}
f(Z(t)) - f(Z(0)) & = & 
 \sum\limits_{n=1}^\infty \left[ \ind_{\{\sigma_n \leq t\}} [f(Z(t \wedge \theta_n)) - f(Z(\sigma_n))] \right. \\
& & \quad \quad + \left. 
\ind_{\{\theta_{n-1} \leq t\}} [f(Z(t \wedge \sigma_n)) - f(Z(\theta_{n-1}))] \right] \\
&   = &  \sum\limits_{n=1}^\infty \ind_{\{\theta_{n-1} \leq t\}} [f(Z(t \wedge \sigma_n)) - f(Z(\theta_{n-1}))], 
\ea
\]
where the last equality is a result of the fact that  $Z(t) \in N_\ve(0)$ for $t \in [\sigma_n,\theta_n]$ on the set $\{\theta_n < \infty\}$ 
(and for $t \in [\sigma_n, \infty)$ on the set $\{\sigma_n < \infty, \theta_n = \infty\}$) and $f$ is constant on $N_\ve(0)$. 
Fix $n \in \N$.  Then the uniqueness of the GPS ESP proved in Theorem \ref{th-exiunigps} along with Lemma \ref{lem-esp0} and Theorem \ref{th-esp4} show that 
on the set $\{\theta_{n-1} \leq t \}$, the process $Z(\cdot \wedge \sigma_{n}) - Z(\theta_n)$
admits the  decomposition (\ref{decomp}) with $0$ replaced by $\theta_{n-1}$ and $t$ replaced by $t \wedge \sigma_n$. 
Therefore, applying It\^{o}'s formula, we obtain the following equality $\P_z$ a.s. on the set $\{\theta_{n-1} \leq t\}$:  
\[ 
\ba{rcl}
f(Z(t \wedge \sigma_n)) - f(Z(\theta_{n-1})) & = & \ds \int_{\theta_{n-1}}^{t\wedge \sigma_n} \ol f(Z(s))\, ds \\
& & + 
 \ds \int_{\theta_{n-1}}^{t \wedge \sigma_n} \nabla f(Z(s)) \cdot dM(s) \\
& & + \ds 
\int_{\theta_{n-1}}^{t\wedge \sigma_n}  \nabla f(Z(s)) \cdot \gamma(s) \, d|Y|(s)   
\ea
\]
with $\gamma(s) \in d(Z(s))$   $d|Y|$ a.e.\ $s \in [\theta_{n-1}, \sigma_n \wedge t]$. 
Multiplying both sides of the last display by $\ind_{\{\theta_{n-1}\leq t\}}$, summing over $n \in \N$ and observing that 
due to the fact that $\nabla f$ and $\ol f$ are identically zero in an $\ve$-neighbourhood of $0$ 
(since $f$ is constant there), we have the equalities 
\[  \suli_{n=1}^\infty \ind_{\{\theta_{n-1} \leq t\}} \int_{\theta_{n-1}}^{t \wedge \sigma_n} \nabla f(Z(s))\cdot dM(s) = 
\int_0^t \nabla f(Z(s)) \cdot dM(s) \]
and 
\[ \ds \suli_{n=1}^\infty  \ind_{\{\theta_{n-1} \leq t\}} \int_{\theta_{n-1}}^{t\wedge \sigma_n} \ol f(Z(s)) \, ds = \int_0^t \ol f(Z(s)) \,ds. \]
Thus we conclude that 
\[ 
\ba{rcl} 
f(Z(t)) - f(Z(0)) & = &  \ds \int_{0}^{t} \ol f(Z(s))\, ds  + \ds \int_{0}^{t} \nabla f(Z(s)) \cdot dM(s) \\
& & + \ds \suli_{n=1}^\infty \ind_{\{\theta_{n-1} \leq t\}} \int_{\theta_{n-1}}^{t\wedge \sigma_n}  \nabla f(Z(s)) \cdot \gamma(s) \, d|Y|(s). 
\ea
\]
Since $\nabla f$ is bounded, the second term on the right-hand side  is an  $\{{\cal F}_t\}$-martingale.  
On the other hand, the last term is non-negative (for every $t \geq 0$) due to the assumed derivative condition on $f$ and the fact that $\P_z$ a.s. 
$\gamma(s) \in d(Z(s))$  $d|Y|$  a.e.\ $s \in [0,t]$.  
Rearranging the terms above,  we see that the process $\{f(Z(t)) - f(Z(0)) - \int_0^t \ol f(Z(s)) ds, t \in [0,\infty)\}$ is a 
$\P_z$-submartingale. By the definition of $\Q_z$, this in turn immediately implies the second property (\ref{eq-submg}). 

For the third property,  we first show that the time the process spends at the origin $\{0\}$ 
has zero Lebesgue measure.  First, observe that  $Z(s) \in \R_+^J$ implies $Z(s) = 0$ if and only if 
$ Z(s) \cdot d_{J+1}  = 0$. 
 Thus by the definition of $\Q_z$, in order to prove property (\ref{zero})  with $\partial \sm S$ 
replaced by $\{0\}$, 
 it suffices to show that for every $T < \infty$, 
\be
\label{surrogate}
 \E_{z} \left[\int_0^T \ind_{\{0\}} (Z(s) \cdot d_{J+1}) \, ds\right] = 0.   
\ee
Since  $Z$ is a strong solution, it 
satisfies the ESP for $X$  
and so by  Corollary \ref{lem-1dsp} we know  that for every $s \in [0,\infty)$, 
$Z(s) \cdot d_{J+1}  = \Gamma_1 (X \cdot d_{J+1} )(s)$. 
Thus $Z \cdot d_{J+1}$ is a one-dimensional reflected diffusion with diffusion coefficient 
$\sum_{i,j=1}^J a_{ij}/J > 0$, where the  strict inequality is due to Assumption \ref{as-ol}(2). 
Hence (\ref{surrogate}) is a consequence of  the well-known result that  
 a non-degenerate one-dimensional reflected diffusion spends 
a.s.\ zero Lebesgue time at the origin $0$ (see, for example, page 90 of \cite{fre}).  
This completes the proof of the theorem when $J = 2$. 

  \end{proof}

\begin{remark}
{\em 
The proof of Theorem \ref{th-submg} in fact shows 
that the first two properties of Definition \ref{def-submg} 
are satisfied for any $J$-dimensional GPS ESP for $J \geq 2$ and, 
in addition, that  the corresponding process spends 
zero Lebesgue time at the origin. 
In order to complete the verification of Definition \ref{def-submg} for 
arbitrary $J$-dimensional GPS ESP, 
it thus only remains to show that 
 the process spends zero Lebesgue time on the $k$-dimensional faces, $1 \leq k \leq J-2$,  
 of the boundary.   
This can be done using the fact that the local restriction of the GPS reflection matrix to such a face 
satisfies the  completely-${\cal S}$ condition.  
The details are omitted as a more general 
study of the boundary property of these diffusions is carried out in a forthcoming paper.  }
\end{remark}

\section{The Semimartingale Property on $[0,\tau_0]$} 
\label{subs-sm}

In this section, we show that the reflected diffusion $Z$ obtained as a strong solution 
of the SDER associated with the GPS ESP  is a semimartingale on $[0,\tau_0]$. 
It was shown in Theorem \ref{th-exi} that it is a semimartingale  on 
the interval $[0,\tau_0)$.  However, as demonstrated below,   
the transition from establishing the property on the half-open interval $[0,\tau_0)$ to the 
closed interval $[0,\tau_0]$ is more subtle. 
First, in Section \ref{subs-smggen}, we 
formulate general sufficient conditions under which the strong 
solution $Z$ of an SDER associated with a general ESP has the required semimartingale property. 
In Section \ref{subs-gfn}, these sufficient conditions are verified for the GPS ESP. 
This verification involves establishing the existence of certain test functions that 
satisfy Assumption \ref{as-gfn}.  
 The details of this proof are deferred to Section \ref{sec-appendix}.

\subsection{Sufficient Conditions for General ESPs} 
\label{subs-smggen}

The first condition, Assumption \ref{as-gfn}, is the existence of 
a sufficiently smooth function  that satisfies certain oblique derivative 
conditions and whose second derivatives satisfy a certain growth condition. 
Recall that ${\cal V}$ is defined by 
(\ref{def-vset}) and let ${\cal U} \doteq \partial G \sm {\cal V}$. 
Also recall that $\supp [g]$  denotes the support of the function $g$.

\begin{ass} 
\label{as-gfn} 
There exist constants $L, R < \infty$ and $\beta > 0$,  
and a function  $g \in {\cal C}^2 (G^\circ \cup {\cal U})$ that 
satisfy the following properties. 
\begin{enumerate} 
\item 
$\supp [g] \subset N_R({\cal V})$. 
\item 
 There exist $\theta > 0$ and $r \in (0,R)$  
 such that 
\[
\ba{rl} 
a) \hspace{.1in} \lan \nabla g (x), d \ran \geq 0 & 
\mbox{ for } d \in d^1(x)  \mbox{ and } x \in {\cal U},  \\
b) \hspace{.1in}  \lan \nabla g (x), d \ran \geq \theta &  
\mbox{ for } d \in d^1(x)  \mbox{ and } 
x \in {\cal U} \cap N_r({\cal V}). 
\ea
\] 
\item 
$\sup_{x\in G} |g(x)| \vee |\nabla g(x)| \leq L$ and 
for $x \in G^\circ \cup {\cal U}$ 
\[ \dfrac{1}{2} \suli_{i,j=1}^J 
\left| \dfrac{ \partial^2 g(x)}{\partial x_i \partial x_j}  \right| 
\leq \dfrac{L}{[d(x,\cU)]^\beta}. 
\]
\end{enumerate} 
\end{ass} 

We now state the main theorem of this section.  

\begin{theorem}
\label{th-smg2} 
Suppose the 
drift and dispersion coefficients 
$b(\cdot)$ and  $\sigma(\cdot)$ satisfy 
Assumption \ref{as-ol}.  
If there exist  constants $L, R < \infty$ and $\beta > 0$ such that 
\begin{enumerate} 
\item
the ESP   $(G,d(\cdot))$ satisfies Assumption \ref{as-gfn} with 
those constants; 
\item 
 the drift and dispersion coefficients satisfy the 
following bound:   
\be
\label{as-bs3}
 \sup_{x \in N_R({\cal V})} 
\left[ \sum_{i,j=1}^J |a_{i,j}(x)| \vee |b(x)| \right] < \infty;    
\ee 
\item 
given an $\{{\cal F}_t\}$-adapted $K$-dimensional Brownian motion $W$,  for every  $z \in G$, there exists a 
strong Markov, strong solution $Z$ to the associated SDER, which has initial condition $z$ 
and satisfies  
\be
\label{eq-fin} 
 \E_z \left[  \int_0^{t \wedge \tau_0} 
\dfrac{\ind_{(0,R]}\left(d(Z(s), {\cal V})\right)}{[d(Z(s), {\cal V})]^\beta} \, 
ds \right] 
< \infty, 
\ee
for $t \in [0,\infty)$, 
where $\tau_0 \doteq \inf \{t \geq 0: Z(t) \in {\cal V}\}$.
\end{enumerate}  
Then  $Z(\cdot \wedge \tau_0)$ is an $\{{\cal F}_t\}$-semimartingale under $\P_z$. 
\end{theorem}

\begin{remark}
\label{rem-cond5.39}
{\em 
Note that since $b(\cdot)$ and $\sigma(\cdot)$ satisfy the Lipschitz condition (\ref{as-bs}), 
the inequality (\ref{as-bs3}) is automatically 
satisfied if ${\cal V}$ is bounded. 
Although our main application of this result 
to the GPS ESP in Section \ref{subs-gfn} has bounded ${\cal V}$,  
 in anticipation of applications for which 
 ${\cal V}$ is unbounded (see  Remark \ref{rem-uv} in Section 
\ref{subsub-vsets}), we consider the general unbounded case here. }
\end{remark} 

\noi 
{\bf Proof of Theorem \ref{th-smg2}.} 
Due to (\ref{as-bs3}), by taking the constant $L$ in the statement larger, if necessary, 
we can assume that 
\be
\label{lbound}
  \sup_{x \in N_R({\cal V})} 
\left[ \sum_{i,j=1}^J |a_{i,j}(x)| \vee |b(x)| \right] < L.
 \ee
Since $Z$ solves the SDER associated with the ESP $(G, d(\cdot))$, by  (\ref{SDER1}) 
 it follows that  $\P_z$ a.s., 
\[ 
 Z(t \wedge \tau_0) = z + M(t \wedge \tau_0) + A(t \wedge \tau_0), \]
where $M$ and $A$ are defined by 
\[ 
 M(t) \doteq \int_0^t \sigma(Z(s)) \cdot dW(s) \quad \quad  \quad  
A(t) \doteq  \int_0^t b(Z(s)) \, ds  + Y(t), 
\]
and for $0 \leq s \leq t < \infty$, 
 $Y$ satisfies 
\[ Y(t) - Y(s) \in \conv \left[ \cup_{u \in (s,t]} 
d(Z(u)) \right].   
\]
If $z \in \cU$,  
 the theorem follows trivially from the 
fact that $\tau_0 = 0$ and $Y(0) = |Y|(0) = 0$ $\P_z$ \,  a.s.

Hence suppose that 
 $z \in G \sm {\cal V} =  G^\circ \cup {\cal U}$. 
 Due to Assumption \ref{as-ol}(1), 
 $M$ is a martingale and $\int_0^{t \wedge \tau_0} b(Z(s)) \, ds$ 
has bounded variation for every $t \in [0,\infty)$. 
Thus to prove the theorem it suffices to show 
that 
for every $z \in G^\circ \cup {\cal U}$, 
\be
\label{subthm}
  |Y|(t \wedge \tau_0)  < \infty \quad \quad 
\mbox{ for every } t \in [0,\infty) \quad \quad \P_z \mbox{ a.s.\ } 
\ee  
For notational conciseness, we introduce the operator $\oa$ defined by 
\[ \oa g(x) =  \dfrac{1}{2} \suli_{i,j=1}^J  a_{ij} (x) 
\dfrac{\partial^2 g(x)}{\partial x_i\partial x_j} \quad \mbox{ for } g \in {\cal C}_b^2 (G). 
\] 
Also, for 
$\delta \in (0, d(z, \cU))$ and 
$K \in (|z|, \infty)$, define 
\[ 
\sek \doteq \inf \{ t \geq 0: Z(t) \in N_\delta ({\cal V}) 
\mbox{ or } 
 |Z(t)| \geq K \}. 
\]   
Let  the function $g$  and constants $r$, $R$, $\theta$, $\beta$ and  $L$ 
satisfy  Assumption \ref{as-gfn}.  
Since $Z$ is continuous,  \ $t \wedge \sek  < \tau_0$      $\P_z$ a.s. 
Due to the semimartingale decomposition for $Z$  on $[0,\tau_0)$ 
 established in Theorem \ref{th-exi} and the fact that  $Z \in G^\circ \cup {\cal U}$ on $[0,\tau_0)$ $\P_z$ a.s., and  
 $g \in {\cal C}^2(G^\circ \cup {\cal U})$,  It\^{o}'s formula yields 
\be
\label{eqn1}
\ba{l} 
 g (Z(t \wedge \sek)) - g(z) \\
\quad \quad \quad  =  \ds  \int_{0}^{t \wedge \sek} 
\nabla g (Z(s)) \cdot b(Z(s)) \, ds 
+ \int_{0}^{t \wedge \sek} \nabla g (Z(s)) \cdot \gamma(s) \, d|Y|(s) \\ 
 \quad \quad \quad \quad + \ds \int_{0}^{t \wedge \sek}  \nabla g(Z(s)) \cdot \sigma(Z(s)) \, dW(s) 
 + \ds \dfrac{1}{2} \int_{0}^{t \wedge \sek} \oa g (Z(s)) \,  ds,  
\ea
\ee
where, $\P_z$ a.s., $\gamma(s) \in d^1(Z(s))$  $d|Y|$ a.e.\ on 
$[0, t \wedge \sigma_{\delta K}]$. 
Using Assumption \ref{as-gfn}(2)  note that 
\[ 
\ba{l}
\ds \int_0^{t \wedge \sek} 
 \nabla g(Z(s)) \cdot \gamma(s)\,  d|Y|(s) \\
 \quad \quad \quad \quad \quad \geq  
\ds \int_{0}^{t \wedge \sek} 
\ind_{[0,r]} \left(d(Z(s),{\cal V})\right)  \nabla g(Z(s)) \cdot \gamma(s)\,  d|Y|(s) \\
\quad \quad \quad \quad \quad  \geq   \ds \theta \int_{0}^{t \wedge \sek} 
\ind_{[0,r]} (d(Z(s),{\cal V}))\, d|Y|(s). 
\ea
\]
In addition, (\ref{lbound}) and 
properties 1 and 3 of Assumption \ref{as-gfn} imply that 
\[
\int_{0}^{t \wedge \sek}  
\oa g(Z(s)) \, ds
\leq L^2 \int_{0}^{t \wedge \sek}  
\dfrac{\ind_{(0,R]} \left(d(Z(s),{\cal V})\right)}{[d(Z(s),{\cal V})]^\beta} \, ds, 
\] 
and also that 
\[
\ba{rcl}  
\ds \int_{0}^{t \wedge \sigma_{\delta K}} 
\left| \nabla g (Z(s)) \cdot b(Z(s)) \right| \, ds 
& \leq & \ds L^2  \int_{0}^{t \wedge \sigma_{\delta K}} 
 \ind_{(0,R]} \left(d(Z(s),{\cal V}) \right)  \, ds \\
& \leq &  
 \ds L^2 R^\beta \int_{0}^{t \wedge \sigma_{\delta K}} 
\dfrac{\ind_{(0,R]} \left(d(Z(s),{\cal V})\right)}{[d(Z(s),{\cal V})]^\beta} \, ds.  
\ea 
\]   
Taking expectations (conditioned on $Z_0 = z$) of both sides of (\ref{eqn1}), 
the stochastic integral on the right hand side vanishes since 
 $\nabla g$ and $\sigma$ are bounded on $\{x \in G: |x| \leq K\}$. 
Rearranging terms, using the last three displays and 
the bound on $g$ and $\nabla g$ in Assumption \ref{as-gfn}(3), we obtain 
\[
\ba{l}  
\ds \theta \E_z \left[ \int_0^{t \wedge \sek} 
\ind_{[0,r]}\left(d(Z(s),{\cal V})\right) \, d|Y| (s) \right]   \\
\hspace{.3in} \leq  
\ds  \E_z\left[ g(Z(t \wedge \sek)) - g(z)\right]
+ \E_z \left[  
\ds \int_0^{t \wedge \sek} \left|\nabla g(Z(s)) \cdot b(Z(s)) \right| \, ds \right] \\
\hspace{1.0in}
 + \dfrac{1}{2} \ds   \E_z \left[ 
\int_{0}^{t \wedge \sek} \left| \oa g (Z(s))  \right| 
ds \right] \\
   \hspace{.3in} \leq  2L + 
 \ds L^2 (R^\beta + 1) \E_z \left[\int_0^{t \wedge \sek} 
\dfrac{\ind_{(0,R]}\left(d(Z(s), \cU)\right)}{[d(Z(s), \cU)]^\beta} \, ds \right].  
\ea
\]
Let $\tilde{L} \doteq L^2(R^\beta+1)$. 
Using the fact that $P_z$ a.s.\  
$t \wedge \sigma_{\delta K}  \uparrow t \wedge \tau_0$ as $K \uparrow \infty$ and $\delta \downarrow 0$,   
we first let $K \uparrow \infty$ and then  $\delta \downarrow 0$, and 
 use the monotone convergence theorem 
to obtain 
\[
  \E_z \left[  
\int_{0}^{t \wedge \tau_0} \ind_{[0,r]} \left(d(Z(s),{\cal V})\right) \, d|Y|(s)
      \right] 
\leq 
  \dfrac{2L}{\theta} + 
\dfrac{\tilde{L}}{\theta} \E_z \left[   \int_0^{t \wedge \tau_0} \dfrac{\ind_{(0,R]}
\left(d(Z(s),{\cal V})\right)}{[d(Z(s), \cU)]^\beta} \, ds \right]. 
\]
When combined with  (\ref{eq-fin}), this shows that 
\be
\label{eq-penul}
 \int_{0}^{t \wedge \tau_0} \ind_{[0,r]} (d(Z(s),{\cal V})  \, d|Y|(s)  < \infty  \quad \quad \P_z \, \mbox{ a.s. }  
\ee
Now define the random time 
\[ \kappa \doteq \sup \{ t \leq \tau_0: d(Z(t), {\cal V}) > r \}. \]
Since $Z$ has continuous paths, 
$\kappa < \tau_0$ $\P_z$ a.s. on $\{\tau_0 < \infty\}$. 
Also, trivially, we have $t < \tau_0$ $\P_z$ a.s.\ on $\{\tau_0 = \infty\}$. 
Together, this implies $t \wedge \kappa  < \tau_0$ $\P_z$ a.s. 
Therefore from Theorem \ref{th-exi}, it follows that 
$|Y|(t \wedge \kappa) < \infty$ $\P_z$ a.s. 
When combined with (\ref{eq-fin}) 
and (\ref{eq-penul}),  this yields for every $t \in [0,\infty)$, 
\rear{2.1} 
\[ 
\ba{rcl} 
 |Y|(t \wedge \tau_0) 
& = &  \ds |Y|(t \wedge \kappa)  
+   \int_{t \wedge \kappa}^{t \wedge \tau_0} \, d|Y|(s)  \\
& = & 
\ds   |Y|(t \wedge \kappa ) + 
\int_{t \wedge \kappa}^{t \wedge \tau_0 } 
\ind_{[0,r]} (d(Z(s),{\cal V})  \, d|Y|(s)  \\ 
& < & \infty \quad \quad \quad \P_z  \mbox{ a.s., } 
\ea
\]
which establishes (\ref{subthm}) and thus 
proves the theorem. 
 \hfill \ink \\

\subsection{Verification of Sufficient Conditions for  GPS RBMs} 
\label{subs-gfn} 

In this section, we verify condition (\ref{eq-fin}) and Assumption \ref{as-gfn} 
 for reflected diffusions  associated with  (a slight generalization of) the GSP  
ESP. When $b(\cdot)$ and $\sigma(\cdot)$ are either continuous or locally bounded, 
 note that condition (\ref{as-bs3}) holds trivially since the set ${\cal V} = \{0\}$ is bounded 
for the GPS ESP.

\subsubsection{Verification of condition (\ref{eq-fin}) for GPS reflected diffusions}
\label{subs-ver1}

Given drift and dispersion coefficients $b(\cdot)$ and $\sigma(\cdot)$ that satisfy 
Assumption \ref{as-ol} and a $J$-dimensional GPS ESP with $J \geq 2$, 
let $Z$ be the unique strong solution to the corresponding SDER (which exists by 
Corollary \ref{cor-exi}).  Since  $G = \R_+^J$ and ${\cal V} = \{0\}$, 
 $d(Z, {\cal V})$ is proportional to $\lan Z, d_{J+1}\ran$.  On the other hand, by  
Corollary \ref{lem-1dsp},  $\lan Z, d_{J+1} \ran$ is a non-degenerate, one-dimensional 
reflected diffusion. 
Thus verifying the condition (\ref{eq-fin})  reduces to checking a property (see (\ref{prop5}) below) of 
one-dimensional reflected diffusions. 
We first prove an estimate for one-dimensional RBMs in Lemma \ref{le-1dbm} 
and then extend the result  in Lemma \ref{le-1ddiff} to one-dimensional reflected diffusions using Girsanov transormations  
and a time-change argument.  
Below, $\Gamma_1$ is the one-dimensional SM defined in (\ref{def-1dsmap}). 

\begin{lemma} 
\label{le-1dbm} 
Given a   standard one-dimensional BM $W_1$ defined on a filtered probability space 
$((\Omega, {\cal F}, \P), \{{\cal F}_t\})$,  
 let   
\be
\label{rbmtau} 
\tau_0 \doteq \inf \{t \geq 0: W_1(t) = 0 \}.
\ee
Then for any $z \in \R_+$, $R \in (0,\infty)$ and  $\ve \in [0,1)$,  
\be
\label{prop4}
  \E_z \left[ \int_0^{t \wedge \tau_0} \dfrac{\ind_{(0,R]} 
(W_1(s))}{W_1^{1+\ve}(s)} ds \right] < \infty  \quad  \mbox{ for } t \in [0,\infty) 
\ee
where $\E_z$ denotes expectation with respect to $\P$, conditioned on $W_1(0) = z$. 
\end{lemma} 
\begin{proof} When $z = 0$, $\tau_0 = 0$ and the lemma holds trivially. Fix $z, R \in (0,\infty)$ and $\ve \in [0,1)$.  
In order to prove (\ref{prop4}) 
we will use the well-known fact (see page 379 of \cite{kal}) that 
for any Borel measurable function $f$, 
\be
\label{green}
 \E_z \left[ \int_{0}^{\tau_c \wedge \tau_a} f(W_1(s)) ds\right]
= \int_{a}^c g_{a,c} (z,y) f(y) m(dy), 
\ee
where $\tau_a$ and $\tau_c$ are defined by (\ref{rbmtau}) with 
$0$ replaced by $a$ and $c$, respectively, $m(dy) =  2 dy$ is 
the speed measure for BM (see Section II.4 in Appendix I.13 of \cite{borsal}) and  
$g_{a,c}(\cdot, \cdot)$ is the Green's function for standard  
one-dimensional BM on $(a,c)$, which is given by 
\[ g_{a,c}(z,y) \doteq  \dfrac{2 (z \wedge y - a) (c - z \vee y)}{c-a}  \quad \mbox{ for } 0 < a \leq c. \] 
Let $a \doteq 0$  and choose $c \geq R \vee z$.  
Substituting the measurable function $f(y) = \ind_{(0,R]}(y)/y^{1+\ve}$ into 
(\ref{green}), we then obtain 
\[ \ds \E_z \left[ \int_{0}^{\tau_c \wedge \tau_0 } 
\dfrac{\ind_{(0,R]}(W_1(s))}{W_1^{1+\ve}(s)} \, ds \right] 
= 2 \int_0^R  \dfrac{g_{0,c}(z,y)}{y^{1+\ve}} \, dy.   
\]
On the other hand, by the definition of $g_{0,c}(\cdot,\cdot)$, 
we see that 
\[ \dfrac{1}{2} \int_0^R  \dfrac{g_{0,c}(z,y)}{y^{1+\ve}}  \, dy
= \left\{ 
\ba{rl} 
\ds \int_z^R \dfrac{z(c-y)}{cy^{1+\ve}} \, dy + \int_{0}^z \dfrac{y(c-z)}{cy^{1+\ve}} \, dy & \mbox{ if } z \leq R \\
 \ds \int_{0}^R \dfrac{y(c-z)}{cy^{1+\ve}} \, dy & \mbox{ if } R \leq z,  
\ea
\right.
\]
and elementary calculations show that 
\[ \dfrac{1}{2} \int_0^R  \dfrac{g_{0,c}(z,y)}{y^{1+\ve}} \, dy  = 
\left\{ 
\ba{rl} 
\dfrac{z^{1-\ve}}{\ve(1-\ve)} - \dfrac{zR^{-\ve}}{\ve} - \dfrac{zR^{1-\ve}}{c(1-\ve)}  & 
\mbox{ if } \ve \in (0,1), z \leq R \\
\dfrac{R^{1-\ve}}{1-\ve} - \dfrac{zR^{1-\ve}}{c(1-\ve)}  & \mbox{ if } \ve \in (0,1), z \geq R \\
z \ln \left( \dfrac{R}{z} \right) + z - \dfrac{zR}{c}  & \mbox{ if } \ve = 0, z \leq R \\
R \left( 1 - \dfrac{z}{c} \right) & \mbox{ if } \ve = 0, z \geq R.  
\ea
\right. 
\]
Invoking the monotone convergence theorem, using the non-negativity of the integrand 
and referring to the last display, we conclude that 
\[
\ba{rcl} \ds
\E_z \left[ 
\int_{0}^{t \wedge \tau_0}   \dfrac{\ind_{(0,R]}(W_1(s))}{W_1^{1+\ve}(s)} \, ds \right] 
& = & \ds \lim_{c \uparrow \infty} \E_z \left[  \int_{0}^{t \wedge \tau_c \wedge \tau_0}  
 \dfrac{\ind_{(0,R]} (W_1(s))}{W_1^{1+\ve}(s)} \, ds \right] \\
& \leq & \ds 2\lim_{c \uparrow \infty} \int_0^R  \dfrac{g_{0,c}(z,y)}{y^{1+\ve}} \, dy  \\
& <  &  \infty, 
\ea
\]
which establishes (\ref{prop4}). 
\end{proof}

\begin{lemma} 
\label{le-1ddiff} 
Given a filtered probability space $((\Omega, {\cal F}, \P), \{{\cal F}_t\})$, suppose 
$b^*$ is an $\{{\cal F}_t\}$-adapted, real-valued process 
and $M^*$ is a continuous $\{{\cal F}_t\}$-martingale whose  
 quadratic variation process $V$  is $\P$ a.s.\ continuously differentiable, and 
there exist constants $\lambda > 0$ and $\Lambda < \infty$ such that  $\P$ a.s., 
\be
\label{deriv-conds} 
\sup_{s \in [0,\infty)} b^*(s) \leq \Lambda \quad \mbox{ and } \quad 
\lambda \leq \inf_{s \in [0,\infty)} V^\prime(s) < \sup_{s \in [0,\infty)} V^\prime(s)  \leq \Lambda,  
\ee
where $V^\prime$ denotes the process obtained as the pathwise derivative of $V$. If 
\be
\label{old-sde}
 B(t) \doteq B(0) + \int_0^t b^*(s) \, ds + M^* (t) \quad \mbox{ for } t \in [0,\infty) 
\ee 
and 
\be
\label{rdifftau} 
\tau_0 \doteq \inf \{t \geq 0: B(t) = 0 \}, 
\ee
then for any $z \in \R_+$ and $R \in (0,\infty)$, 
\be
\label{prop5}
  \E_z \left[ \int_0^{t \wedge \tau_0} \dfrac{\ind_{(0,R]} 
(B(s))}{B(s)} \, ds \right] < \infty  \quad \mbox{ for } t \in [0,\infty),  
\ee
where $\E_z$ denotes expectation with respect to $\P$, conditioned on $B(0) = z$. 
\end{lemma} 
\begin{proof} 
We first use a time-change argument to show that we can restrict ourselves, without loss of generality, 
to the case when $M^*$ is a one-dimensional standard BM. 
Define the ``inverse'' $T$ of $V$ by 
\[ T(t) \doteq \inf \{s \geq 0: V(s) > t \} \quad \mbox{ for } t \in [0,\infty). \]
The assumed properties of $V$ ensure that $\P$ a.s.,  both $T$ and $V$ are strictly increasing, continuously
 differentiable functions on $[0,\infty)$ and  $V(T(t)) = T(V(t)) = t$ for every $t \in [0,\infty)$. 
Now let $W_1 \doteq M^*(T(t))$, $\newfil = {\cal F}$, $\newfil_t \doteq {\cal F}_{T(t)}$, $\npr(t) \doteq B (T(t))$ 
for $t \in [0,\infty)$, and define $\newtau_0 \doteq \inf \{ t \geq 0: \npr (t) = 0\}$.   
Then $W_1$ is an $\{\newfil_t\}$-adapted,  standard one-dimensional BM (see Theorem 4.6 of \cite{karshr}) and 
 $\npr$ is given by 
\be
\label{newsde}
   \npr (t) =  B(0) + \int_0^t \newdrift(s) \, ds +  W_1 (t), 
\ee
where  $\newdrift$ is an $\{\newfil_t\}$-adapted process that satisfies 
\[
 \newdrift (s) = \dfrac{b^*(T(s))}{V^\prime(T(s))} \leq \dfrac{\Lambda}{\lambda} \quad \mbox{ for } s \in [0,\infty). 
\]
Moreover, $\tau_0 = T(\newtau_0)$ and 
\rear{2.4}
\[
\ba{rcl}
 \E_z \left[ \ds\int_0^{t \wedge \tau_0} \dfrac{\ind_{(0,R]} (B(s))}{B(s)} \, ds \right] 
 & = &   \ds \E_z \left[ \int_0^{t \wedge T(\newtau_0)} 
\dfrac{\ind_{(0,R]} (\npr(V(s)))}{\npr(V(s))} \, ds \right] \\
& = & \E_z \left[ \ds \int_0^{V(t) \wedge \newtau_0} 
\dfrac{\ind_{(0,R]}(\npr(r))}{\npr(r)\ V^\prime(T(r))} \, dr \right] \\
& \leq & \dfrac{1}{\lambda} \E_z \left[ \ds \int_0^{V(t) \wedge \newtau_0} 
\dfrac{\ind_{(0,R]}(\npr(r))}{\npr(r)} \, dr \right]. 
\ea
\]
This shows that in order to prove the theorem, it suffices to 
establish (\ref{prop5}) for processes 
$\tilde{B}$ of the form (\ref{newsde}), with $\tilde{b}$ uniformly bounded.

We shall now simplify the problem further by applying a Girsanov transformation to remove the drift $\newdrift$ from the 
process $\npr$. 
Fix $t \in [0,\infty)$ and define 
\[ H(s) \doteq \exp \left( - \int_0^s \newdrift(r)\, dW_1(r) - \dfrac{1}{2} \int_0^s \newdrift^2(r)\, dr \right) \quad \mbox{ for } s \in [0,t]. \]
Since $\newdrift$ is bounded, the process $H = \{H(s), s \in [0,t]\}$ is an $\{\newfil_s\}$-martingale with expectation $1$.  
Then by  Girsanov's theorem (see, for example, Theorem 5.1 of \cite{karshr}),  
 the process $\npr = \{\npr_s, s \in [0,t]\}$ is a standard, one-dimensional BM on $((\Omega, \newfil, \Q), \{\newfil_s, s \in [0,t]\})$, 
where $\Q$ is the probability measure on $(\Omega, \newfil)$ defined by 
\[ \Q(A) = \P ( H(s) A)  \quad \mbox{ for every }   A \in \newfil_s,  \quad s \in [0,t].  \]
Also, consider the process $N = \{N(s), s \in [0,t]\}$ defined by  
\be
\label{def-ns}
 N(s)   \doteq \exp \left( \int_0^s \newdrift (r) \, d\npr(r) - \dfrac{1}{2} \int_0^t \newdrift^2(r)\, dr \right) 
\quad \mbox{ for } s \in [0,t] 
\ee
and note that  $N(\cdot) = 1/H(\cdot)$. 
 So, $N$ is  an $\{\newfil_s\}$-martingale under $\Q$ 
and  $\P(A) = \Q( N(s)  A)$ for $A \in \newfil_s$, $s \in [0,t]$.  Let $\E^{\Q}_z$ denote expectation 
with respect to $\Q$, conditioned on $\npr(0) = z$ and, for greater clarity, we denote the corresponding expectation $\E_z$ 
with respect to $\P$ by $\E^{\P}_z$.   Then, using Fubini's theorem, the properties of  
$\Q$ and $\P$ stated above and 
H\"{o}lder's inequality, we obtain  
\rear{2.4}
\[ \ba{l}
 \ds \E_z^{\P} \left[ \ds\int_0^{t \wedge \tau_0} \dfrac{\ind_{(0,R]} (\npr(s))}{\npr(s)}\, ds \right] \\
\quad \quad \quad  =  \ds \int_0^t \E_z^\P \left[ \ind_{[0,\tau_0]}(s) \dfrac{\ind_{(0,R]} (\npr(s))}{\npr(s)}  \right] \, ds \\
\quad \quad \quad =  \ds \int_0^t \E_z^\Q \left[ \ind_{[0,\tau_0]}(s) \dfrac{\ind_{(0,R]} (\npr(s)) }{\npr(s)} N(s) \right]  
\, ds  \\
\quad \quad \quad \leq \ds \left(\int_0^t \E_z^\Q \left[ \ind_{[0,\tau_0]}(s) \dfrac{\ind_{(0,R]} (\npr(s))}{\npr^{3/2}(s)}\right] \, ds \right)^{2/3} 
\left(\int_0^t \E_z^\Q \left[ \left(N(s) \right)^3\right] \, ds \right)^{1/3}   \\
 \quad \quad \quad =  \ds \left(\int_0^t \E_z^\P \left[ \ind_{[0,\tau_0]}(s) \dfrac{\ind_{(0,R]} (W_1(s))}{W_1^{3/2}(s)} \right] \, ds \right)^{2/3} 
\left(\int_0^t \E_z^\Q \left[ \left(N(s) \right)^3\right] \, ds \right)^{1/3}. 
\ea
\]
It only remains to show that each of the two terms in the 
last line are finite. 
An    application of Lemma \ref{le-1dbm} with $\ve = 1/2$ immediately shows that the first term  is finite. 
For the second term, let $N^{3 \tilde{b}(\cdot)}$  be defined like $N$ in (\ref{def-ns}),  but with $\tilde{b}(\cdot)$ replaced 
everywhere by 
$3 \tilde{b}(\cdot)$.  Then it is easy to see that $N^{3 \tilde{b}(\cdot)}$ is an 
$\{\newfil_s\}$-martingale under $\Q$ with expectation $1$. 
When combined with the fact that for $s \in [0,t]$, 
\[  \left(N(s) \right)^3 = N^{3 \tilde{b}(\cdot)}(s) \exp \left(2 \int_0^s \tilde{b}^2(r) \, dr \right)
 \leq N^{3 \tilde{b}(\cdot)} (s) \exp \left(2 t \left[\sup_{r \in [0,t]} \tilde{b}^2(r)\right]\right), 
\]
this shows that the second term is also finite and thus completes the proof.  
\end{proof}

\begin{cor} 
\label{th-estrbm}  
Given  an $\{{\cal F}_t\}$-adapted,  
$K$-dimensional Brownian motion defined on a filtered probability space 
$((\Omega, {\cal F}, \P), {\cal F}_t)$, the $J$-dimensional GPS ESP and  drift and dispersion 
coefficients $b(\cdot)$ and  $\sigma(\cdot)$ satisfying 
Assumption \ref{as-ol}, let  $Z$ be the unique, strong solution to the associated 
SDER and let $\tau_0$ be given by (\ref{tau}). 
Then for  every $z \in \R_+^J$ and  $R \in (0,\infty)$, it follows  that 
\be
\label{estrbm} 
\E_z \left[
\int_0^{t \wedge \tau_0} \dfrac{\ind_{(0,R]}(d(Z(s), {\cal V}))}{d(Z(s),{\cal V})} ds
\right] < \infty \quad \mbox{ for } t \in [0,\infty),
\ee
where $\E_z$ is expectation with respect to $\P$, conditioned on $Z(0) = z$. 
\end{cor}
\begin{proof}
By Lemma \ref{lem-gsp}, we know that  
${\cal V} = \{0\}$ for the family of GPS ESPs and therefore $d(z,{\cal V}) = |z|$.  
Now let $f(z) \doteq \lan z, d_{J+1}\ran$. 
Then $f(z) = 0$ for $z \in \R_+^J$ if and only if  $|z| = 0$, and 
there exist $0 < k_1 < k_2 < \infty$ such that 
\[ k_1 f(z) \leq d(z, {\cal V}) = |z| \leq k_2 f(z) \quad \quad \mbox{ for every } z \in \R_+^J. \]
Hence if $B \doteq f(Z)$ then $\tau_0 = \inf \{ t > 0: B(t) = 0 \}$ and 
(\ref{estrbm})  holds if 
\be
\label{estrbm2}
 \tilde{\E}_{z^*} \left[ \int_0^{t \wedge \tau_0 } \dfrac{\ind_{(0,\tilde{R}]}(B(s))}{B(s)} ds
\right] < \infty  \quad \quad \mbox{ for every } t \in [0, \infty), 
\ee  
where  $\tilde{R} \doteq R/k_1$, $z^* \doteq \lan z, d_{J+1}\ran$ and for $z^* \in \R_+$, 
 $\tilde{\E}_{z^*}$ is expectation with respect to $\P$, 
conditioned on $B(0) = z^*$.   
According to Definition \ref{def-SDER}, $Z = \bGamma (X)$ $\P$ a.s., where $X$ is given by 
(\ref{SDER}); by Corollary \ref{lem-1dsp}, 
$B \doteq f(Z) = \Gamma_1 \left( \lan X, d_{J+1} \ran \right)$, where 
$\Gamma_1$ is the 1-dimensional SM. 
Thus for $t \in [0,\tau_0)$,   $B (t) = \lan X (t), d_{J+1} \ran$, which can be written explicitly as  
\[ B (t) = \int_0^t b^*(s) ds  + M^*(t), \]
where $b^*$ is the $\{{\cal F}_t\}$-adapted process defined by $b^*(s) = \sum_{i=1}^J b_i(X(s))/\sqrt{J}$ and 
$M^*$ is an $\{{\cal F}_t\}$-adapted continuous martingale with quadratic variation process $V$ given by 
\[ V(t) = \sum_{i,j=1}^J \int_0^t \dfrac{a_{ij} (X(s))}{J}, \]
where $a(\cdot) = \sigma (\cdot) \sigma^T(\cdot)$ is the diffusion coefficient defined in Assumption \ref{as-ol}(2). 
The fact that $b(\cdot)$ and $\sigma(\cdot)$ satisfy Assumption \ref{as-ol} ensures that the processes 
$b^*$ and $V$ satisfy the conditions of Lemma \ref{le-1ddiff}. So 
the estimate (\ref{estrbm2})  follows  from Lemma \ref{le-1ddiff}, 
which in turn establishes (\ref{estrbm}).  
\end{proof}

\subsubsection{Verification of Assumption \ref{as-gfn}} 
\label{subs-pfmain2}

In this section, we consider ESPs on conical polyhedral domains with vertex at 
the origin that satisfy Assumption  \ref{as-setb} (the ``set B'' condition) 
and have ${\cal V} = \{0\}$.  
The main result of this section is Theorem \ref{th-gfn},  
which proves the existence of a function that 
satisfies Assumption \ref{as-gfn} for this class of ESPs. 
The theorem relies on two results.  
The first result, stated as Theorem \ref{lem-locgfn}, 
establishes the existence of a family of ``local'' functions 
$\{g_{z,r}, z \in {\cal U}\}$ where, roughly speaking, each $g_{z,r}$ satisfies  
the properties of Assumption \ref{as-gfn} in a  bounded, convex neighbourhood 
of $z$.  
This result holds for any polyhedral ESP 
 that satisfies Assumption \ref{as-setb}. 
The second result is Lemma \ref{lem-cover}, which   
establishes a covering property that shows that the 
local functions constructed in Theorem \ref{lem-locgfn} can be patched 
together to yield a function $g$ that satisfies Assumption \ref{as-gfn}. 
We first introduce some notation, then 
state Theorem \ref{lem-locgfn} 
and Lemma \ref{lem-cover},  and then present the 
proof of the main result, Theorem \ref{th-gfn}.  
The proofs of Theorem \ref{lem-locgfn} and Lemma \ref{lem-cover} 
 are relegated to Section \ref{sec-appendix}.

Given $\bI = \{1, \ldots, I\}$ and a polyhedral ESP $\{(d_i,n_i,c_i), i \in \bI\}$ 
with domain $G$, 
recall that $I(x) \doteq \{i \in \bI: \lan x, n_i \ran = c_i\}$ for 
 $x \in \partial G$.
For $C \subseteq \bI$,  let $F_C$ be defined by 
 \[ F_C \doteq \left\{ x \in \partial G: 
 I(x) = C \right\},   \]
and note that $\partial G$ is the disjoint union of 
$F_C, C \subseteq \bI$, $C \neq \emptyset$.  
Given $C \subseteq \bI$, 
we shall refer to $F_C$ as a facet, to distinguish it from 
its closure 
\[ \{x \in G:\lan x, n_i \ran = c_i \mbox{ for  every } i \in C \}, \] 
which we will refer to as a face (the two definitions coincide if and only if $F_C$ is a point).  
For $z \in F_C \subset \partial G$,    
define 
\be
\label{def-rz}
 r_z \doteq d\left(z, \partial G\sm \left[\cup_{C^\prime \subseteq C} 
 F_{C^\prime}\right]\right), 
\ee
with the convention that the distance of $z$ to the empty set $\emptyset$ is equal to 
zero.  
If $z \in F_C$,  $r_z$ is the minimum distance from $z$ to any face on which 
it does not lie.  
It is not hard to see that $F_C$ is relatively open and $r_z > 0$ as long as $z$ is not a vertex.  
On the other hand, $r_z = 0$ when $z$ is the sole vertex because  
in that case $C = \bI$ and so $\partial G\sm \left[ \cup_{C^\prime \subseteq C} F_{C^\prime} \right] = \emptyset$.  
From the definition of polyhedral ESPs it follows that 
\be
\label{inclusion}
 x \in N^\circ_{r_z}(z) \cap G 
\quad \mbox{ implies that } \quad  d(x) \subseteq d(z).  
\ee
Given any subset $U \subset \partial G$, 
we define 
\be 
\label{def-pu}
 {\cal P} (U) \doteq \{ C \subseteq \bI: x \in F_C 
\mbox{ for some } x \in U \}
\ee
to be the collection of sets $C$ such that $U$ has a non-empty intersection with 
the corresponding facet $F_C$. 
Observe that then the set ${\cal U} \doteq \partial G \sm {\cal V}$  
can be written in the form ${\cal U} = \cup_{C \in {\cal P}({\cal U})} F_C$. 
For example, in the case of the GPS ESP with ${\cal V} = \{0\}$, we have 
${\cal P}({\cal U}) = \{C \subset \bI: C \neq \bI, C \neq \emptyset\}$.  
We can now state the main results of the section. 
Recall that $N_1(0)$ is the open unit ball centered at $0$. 

\begin{theorem} 
\label{lem-locgfn} 
Suppose the polyhedral ESP $\{(d_i,n_i,c_i), i \in \bI \}$ 
satisfies Assumption \ref{as-setb}, 
and let \, ${\cal U} = \partial G \sm {\cal V}$, 
where the ${\cal V}$-set is defined by (\ref{def-vset}).  
There exists a function  
$A:{\cal U} \ra [0,\infty)$, constants $A^\prime < \infty$ and $\theta > 0$,  bounded convex 
sets $Q^+_C, Q^-_C,  C \in  {\cal P}({\cal U})$, such that 
$0 \in (Q_C^-)^\circ$ and $Q_C^- \subset (Q_C^+)^\circ \subset N_1(0)$,  
and a family of functions 
$\{g_{z,r}: z \in {\cal U}, r \in (0, r_z) \}$ 
satisfying the following properties.  
\begin{enumerate} 
\item 
$g_{z,r} \in {\cal C}^\infty (G) \mbox{ and }  g_{z,r} (x) = \nabla g_{z,r} 
(x) = 0$ for $x \in {\cal V}$; 
\item 
 $\supp [g_{z,r}] \cap G \subset  \left( z + r Q^+_{I(z)}\right)^\circ$; 
\item 
 $\sup_{x \in G} |g_{z,r}(x)| \leq A^\prime r$ for every $z \in {\cal U}, r \in (0,r_z)$; 
\item 
$ \lan  \nabla g_{z,r} (x), d \ran  \geq   0 
\mbox{ for } d \in d^1(x) \mbox{ and }  x \in {\cal U}$;  
\item   
$ \lan \nabla g_{z,r}(x), d \ran  \geq   \theta   
\mbox{ for } d \in d^1(x)   \mbox{ and } 
x \in  \left[z + r Q^-_{I(z)}\right] \cap \partial G$;  
\item 
 $\sup_{x \in G} \suli_{i,j=1}^J \left| \dfrac{\partial^2 g_{z,r}(x)}{\partial x_i \partial x_j} 
\right|  \leq  \dfrac{A (z)}{r}$;  
\item 
For every $R < \infty$, $\tilde{A}_R \doteq \sup_{z \in N_R({\cal V}) \cap {\cal U}} A(z) < \infty$. 
\end{enumerate}
\end{theorem}

\begin{lemma}
\label{lem-cover}
Suppose the polyhedral  ESP $\{(d_i,n_i,0),i \in \bI \}$ has  
a conical domain  with angle less than $\pi$ and has 
${\cal V} = \{0\}$, where ${\cal V}$ is defined by (\ref{def-vset}).  
Let $Q^+_C$, $Q^-_C$, $C \in {\cal P}({\cal U})$,  
 be  convex sets such that $0 \in (Q_C^-)^\circ$ 
and $Q_C^- \subset (Q_C^+)^\circ \subset N_1(0)$. 
Then, given $0 < r <  R < \infty$,  
there exists a countable set of vectors 
$S \subset {\cal U}$ 
and a corresponding set of scalars $\{\rho_z: z \in S, \rho_z < r_z\}$ 
such that the sets $Q_z^- \doteq z + \rho_z Q^-_{I(z)}$ and  
$Q_z^+ \doteq z + \rho_z Q^+_{I(z)}$, $z \in S$,  satisfy the following five properties.   
\begin{enumerate} 
\item 
There exists $N < \infty$ such that 
for every $x \in G$ 
\be
\label{card}
  \# \left[\{z \in S: x \in 
Q_z^+\}\right] \leq N \, ;
\ee
\item 
There exists $\nu \in  (0, \infty)$ such that 
if $x \in  Q_z^+$ then 
$\nu \rho_z \geq d (x, {\cal V})$;  
\item  
If $x \in \partial G \cap Q_z^+$ then 
$d^1(x) \subseteq d^1(z)$; 
\item 
${\cal U} \cap N_r({\cal V}) \subset [\cup_{z \in S} Q^-_z]$; 
\item  
$\left[ \cup_{z \in S} 
Q_z^+ \right] \subset N_R ({\cal V})$. 
\end{enumerate} 
\end{lemma}

Observe that Theorem \ref{lem-locgfn} (in particular, the existence and shape of the convex sets 
$Q_C^-$ and $Q_C^+$) is  heavily dependent on the geometry of the directions of constraint 
and is not much concerned with the structure of the set ${\cal U}$. 
On the other hand, the  covering result in Lemma \ref{lem-cover}  depends more on the 
geometry of ${\cal U}$.  
The proofs of Theorem \ref{lem-locgfn} and Lemma \ref{lem-cover} are given in 
Sections \ref{subs-ap1} and \ref{subs-ap2} respectively. 
Here, we show how these results can be combined to 
construct a function $g$ that satisfies Assumption \ref{as-gfn}.

\begin{theorem} 
\label{th-gfn} 
Suppose the polyhedral ESP  $\{(d_i,n_i,0), i \in \bI \}$ has a conical 
domain with angle less than $\pi$, 
satisfies Assumption \ref{as-setb} and has 
${\cal V} = \{0\}$.  
Then for any $R < \infty$, there exists $L < \infty$ such that the ESP  satisfies Assumption \ref{as-gfn} with $\beta = 1$. 
\end{theorem}  
\begin{proof}
Fix $R < \infty$. Let $\tilde{A}_R < \infty$, $\theta > 0$, 
$\{g_{z,r}: z \in {\cal U}, r \in (0, r_z) \}$ 
and $Q_C^+$, $Q_C^-$, $C \in {\cal P}({\cal U})$,   
 satisfy the properties stated in Theorem \ref{lem-locgfn}. 
Fix $r \in (0,R)$ and 
corresponding to $Q_C^+$, $Q_C^-$, $C \in {\cal P}({\cal U})$, 
choose a countable set of points $S \subset {\cal U}$  and 
$\rho_z \in (0,r_z)$,  $z \in S$, 
such that  the properties stated in 
Lemma \ref{lem-cover} are satisfied, and 
 let the corresponding sets $Q_z^-$ and $Q_z^+$, $z \in S$, be as defined in Lemma \ref{lem-cover}.  
Also define 
\be
\label{eq-gfn} 
g(x) \doteq 
 \suli_{z\in S} g_{z,\rho_z} (x) \quad  \mbox{ for } x \in G.  
\ee 
For $x \in G^\circ \cup {\cal U}$, let  $J(x) \doteq \{z \in S: x \in (Q_z^+)^\circ\}$. 
Since, by property 1 of Lemma \ref{lem-cover}, the cardinality of $J(x)$ is finite (it is in fact uniformly bounded by $N$), 
 there exists $\ve > 0$ such that 
$N_\ve(x) \subset \cap_{z \in J(x)} (Q_z^+)^\circ$, and so 
for every $y \in N_\ve(x)$, $J(y) = J(x)$.  
Thus for every $x \in G^\circ \cup {\cal U}$, there exists $\ve > 0$ such that 
\be
\label{rel-gfn}
 g(y) =  \suli_{z \in J(x)} g_{z,\rho_z} (y)   \quad \quad \mbox{ for } y \in N_\ve(x). 
\ee
Along with properties 1--3 of Theorem \ref{lem-locgfn}, this guarantees that 
$g$  lies in 
${\cal C}^\infty(G^\circ \cup {\cal U})$, satisfies $g(x) = 0$ for $x \in {\cal V}$ and is 
a continuous function on $G$.

Now, property 5 of Lemma \ref{lem-cover} ensures 
 that $\supp [g] \subset N_R({\cal V})$,  
which establishes Assumption \ref{as-gfn}(1). 
On the other hand, (\ref{rel-gfn}), with $y = x$, combined with 
 properties 4 and 5 of Theorem \ref{lem-locgfn}, 
 property 4 of Lemma \ref{lem-cover}  and the fact that 
$\rho_z < r_z$, implies 
Assumption \ref{as-gfn}(2). 
 Furthermore, (\ref{rel-gfn}), along with properties  6  and  7 of 
   Theorem \ref{lem-locgfn} and properties 1 and 2 of Lemma 
\ref{lem-cover}, implies that 
for $x \in G^\circ \cup {\cal U} \cap N_R({\cal V})$, 
\[
\suli_{j,k=1}^J  \left|
\dfrac{\partial^2 g (x)}{\partial x_j \partial x_k}  \right|
 \leq 
\suli_{z\in J(x)} \suli_{j,k=1}^J \left| \dfrac{\partial^2  
g_{z,\rho_z} (x)}{\partial x_j \partial x_k}  \right| 
 \leq  \dfrac{\tilde{A}_R N }{\min_{z \in J(x)} \rho_z} 
 \leq  \dfrac{\tilde{A}_R N \nu}{d(x, {\cal V})}. 
\]
Since $g$ is identically zero outside $N_R({\cal V})$, this
shows that Assumption \ref{as-gfn}(3)
is satisfied with $L = \tilde{A}_R N \nu$ and $\beta = 1$, thus 
 completing the proof of the theorem. 
\end{proof}

\begin{theorem} 
\label{th-rbm}
Let $Z$ be the pathwise unique, strong solution to the SDER associated with the 
GPS ESP and drift and dispersion coefficients $b(\cdot)$ and 
$\sigma (\cdot)$ that satisfy Assumption \ref{as-ol}. 
Then $Z(\cdot \wedge \tau_0)$ is a semimartingale,  
where $\tau_0 = \inf\{t \geq 0: Z(t) = 0\}$. 
\end{theorem}
\begin{proof}
The bound (\ref{as-bs3}) holds by Remark \ref{rem-cond5.39}. 
Moreover,   given any $R < \infty$, Theorem \ref{th-gfn} shows there exist $L < \infty$ and $\beta = 1$ such 
that  Assumption \ref{as-gfn} is satisfied 
for the GPS ESP. 
The existence of a pathwise unique, strong solution follows 
from Corollary \ref{cor-exi}, and 
Corollary \ref{th-estrbm} shows that condition (\ref{eq-fin}) is satisfied 
for any $R < \infty$ and  $\beta = 1$. 
Thus the result follows from Theorem \ref{th-smg2}.  
\end{proof}

\section{Construction of Test Functions for the GPS Family}
\label{sec-appendix}

In this section, we consider a slight generalization of the family of GPS ESPs, 
namely polyhedral ESPs in convex conical polyhedral domains  with vertex at the 
origin that satisfy Assumption \ref{as-setb}  and have ${\cal V} = \{0\}$.  
In Section \ref{subs-ap1}, we prove Theorem \ref{lem-locgfn} 
and in Section \ref{subs-ap2}, we establish Lemma \ref{lem-cover}. 
Together with Theorem \ref{th-gfn}, this demonstrates the existence of a 
function $g$ that satisfies Assumption \ref{as-gfn} for this class of ESPs. 

\subsection{Proof of Theorem \ref{lem-locgfn}} 
\label{subs-ap1}

We first prove some preliminary results in Lemmas \ref{lem-dupset}-\ref{lem-qsets}.  
The following notation is used throughout this section. 
For any set $A \subset \R^J$, 
$\relint (A)$ is used to denote the relative interior of the set $A$ 
(see \cite{roc} for a precise definition).  
For conciseness,  in this section 
 we will often use $A^{\ve}$ to denote 
${N}_\ve (A) = \{ x: d(x, A) \leq \ve \}$ for $\ve > 0$. 
Recall the definitions given in Section \ref{subs-pfmain2} 
of $F_C$, $C \subseteq \bI$, and ${\cal P}({\cal U})$, where 
${\cal U} = \partial G\sm {\cal V}$. 
For $C \subset \bI$, $C \neq \emptyset$, let 
the cone $L_C$ be defined by 
\be 
\label{def-LC}
L_C \doteq \left\{ -\suli_{i \in C} a_i d_i :a_i \geq 0 \right\},  
\ee
and the set $K_C$ by 
\be 
\label{KC} 
  K_C \doteq \left\{ -\suli_{i \in C} a_i d_i/|d_i| :a_i \geq 0, 
\suli_{i \in C} a_i = 1 \right\}.  
\ee 
Note that  $K_C$ is a convex, compact subset of $\R^J$, and 
for $x \in F_C$, 
$L_C = -d(x)$ and $-d^1(x) \subseteq \{t K_C: t \geq 1\}$. 
We first state a useful consequence of the existence of a set $B$ that satisfies 
Assumption \ref{as-setb} for the ESP. 
This result was proved in \cite{dupish3}. 

\begin{lemma} 
\label{lem-dupset} 
Suppose the  polyhedral ESP $\{(d_i,n_i,c_i), i \in \bI\}$ 
 satisfies Assumption \ref{as-setb}.  
For $C \subseteq \bI$, $C \neq \emptyset$,   if 
 $L_C$ is defined by (\ref{def-LC})   
then   
\[ \min_{i \in C} \lan n_i , d \ran < 0 \hspace{.1in} \mbox{ for every } 
d \in L_C\sm \{0\}. 
\] 
\end{lemma} 
\begin{proof}
This lemma corresponds to Lemma A.3 in 
\cite{dupish3} (specialised to the case when 
the vector field $\gamma_i$ in \cite{dupish3} 
is constant and equal to $-d_i$),  
 with the caveat that  $n_i$ used 
in \cite{dupish3} 
represents an outward normal, 
while in this paper $n_i$ denotes an  
inward normal to the domain $G$. 
Note that the condition (A.1) specified in 
\cite{dupish3} follows from Assumption \ref{as-setb} 
due to Lemma 2.1 of \cite{dupish1}. 
\end{proof}

Since Assumption \ref{as-setb} holds for the 
polyhedral ESPs under consideration,  Lemma \ref{lem-dupset} is applicable. 
If $x \in {\cal U}$, 
then the vectors $d_i, i \in I(x),$ are linearly 
independent. 
Hence for $C \in {\cal P}({\cal U})$, 
$0 \not \in K_C$. 
By  
 Lemma \ref{lem-dupset} 
this implies that for $C \in {\cal P}({\cal U})$, 
\[ \min_{i \in C} \lan n_i, d \ran  < 0  \hspace{.1in} 
\mbox{ for all } d \in K_C. \]
Since $K_C$ is compact, 
there exists $\delta_C > 0$  
  such that 
\be
\label{ineq-min}
 \min_{i \in C} \lan n_i, d \ran < 0 \hspace{.1in} 
\mbox{ for all } d \in K_C^{\delta_C}. 
\ee  
Let $\kcdc$  be a closed, convex set  that has a ${\cal C}^\infty$ boundary and 
satisfies 
\be
\label{def-kcdc}
K_C^{\delta_C/2} \subset  (\kcdc )^\circ  \subset  \kcdc  \subset K_C^{\delta_C}.  
\ee
(Here a convex set $F \subset \R^J$ is said to have a 
$C^\infty$ boundary if for every point $y \in \partial F$, there exists 
a (relative) neighbourhood of $y$ in $\partial F$ that  is  a ${\cal C}^\infty$ submanifold 
of $\R^J$, appropriately modelled on some hyperplane of $\R^J$;   
for closed convex sets $F$, this has been shown in \cite{hol} to be equivalent 
to the gauge function of $F$ being $C^\infty$ in a neighbourhood of the boundary of $F$.) 
Also, define  
\be
\label{LC}
 L_{C,\delta_C} \doteq \cup_{t \geq 0} t \kcdc.  
\ee
Then the inequality (\ref{ineq-min}) clearly holds with 
$K_C^{\delta_C}$ replaced by $\kcdc$. 
This  in turn implies that 
 $0 \not \in \kcdc$, 
  that $L_{C,\delta_C}$ is a (half) cone with  vertex at the origin,  
and that there exists 
 $\beta_C > 0$  
such that 
\be 
\label{eq-theta} 
 \min_{i \in C} \lan n_i, d \ran \leq -2 \beta_C |d| 
\hspace{.1in} \mbox{ for all } d \in L_{C,\delta_C}. 
\ee 
It is also clear that, since $\kcdc$ has a ${\cal C}^\infty$ boundary, 
 $\lcdc$ is a cone whose boundary is ${\cal C}^\infty$ 
everywhere except at the vertex $\{0\}$.

In Lemma \ref{lem-gc} below, we construct a family $\{g_C, C \in {\cal P}({\cal U})\}$ of functions, where     
$g_C$ is  associated with the facet $F_C$ in ${\cal U}$. 
These functions serve as the basic building blocks for the construction of 
the family of functions 
$\{g_{z,r}, z \in {\cal U}, r \in (0,r_z)\}$ of Theorem \ref{lem-locgfn}; indeed the latter 
will essentially be obtained as suitably scaled translates of the 
functions $g_C$, $C \in {\cal P}({\cal U})$.  
Each function $g_C$ is constructed as (a suitable approximation of) 
the distance function to the cone $\lcdc$. 
As shown below in Lemma \ref{lem-gc}, 
the geometry of the directions of constraint (imposed by 
Assumption \ref{as-setb}) ensures that this distance function 
 locally satisfies the necessary gradient conditions 
(see property 2 of the lemma). 
This observation was first made in \cite{dupish3}  when considering 
an SP with ${\cal V} =\emptyset$, and was used there  
 to construct a ${\cal C}^1$ function that satisfies gradient conditions 
similar to those in (\ref{gradgc}) and (\ref{newgradgc}). 
However, 
the construction here is considerably more involved due to the fact that 
${\cal V} \neq \emptyset$ and we need a ${\cal C}^2$ function whose 
second-order partial derivatives are uniformly bounded. 
In particular, since the distance function is not ${\cal C}^2$, 
  we need to establish the existence of sufficiently smooth approximations to  
the distance function that satisfy both the gradient properties and the bound on the  second-order derivatives. 
While discussing approximations, for conciseness we will use the Schwarz notation for multi-indices: given $\alpha \in 
\Z_+^J$ and a function $f$ on some open set $\Omega \subset  \R^J$, recall that 
\[ D^\alpha f = \dfrac{\partial^{\alpha_1 + \ldots \alpha_J} f}{\partial x_1^{\alpha_1} \ldots \partial x_J^{\alpha_J}} 
\quad \mbox{ and } \quad 
D_w^\alpha f = \left(\dfrac{\partial^{\alpha_1 + \ldots \alpha_J} f}{\partial x_1^{\alpha_1} \ldots x_J^{\alpha_J}}\right)_w 
\] 
 denote the ordinary and weak derivative operators of the function $f$ of order $\alpha$ on $\Omega$ and 
 $|\alpha|$  denotes $\alpha_1 + \ldots + \alpha_J$  
(see \cite[Definition 2, page 19]{bur} for the definition of weak derivatives). 
With some abuse of notation, we will  say $h = D_w^\alpha f$ to mean $h$ is a 
weak derivative of the function $f$ of order $\alpha$ on $\Omega$.

\begin{lemma}
\label{lem-gc}  For every $C \in {\cal P} ({\cal U})$, 
given  any $0 < \down_C <  \up_C < \infty$, $\tilde{\ve}_C > 0$ 
and 
\[ \domain \doteq \left(\lcdc^{\up_C}\right)^\circ \sm \lcdc^{\down_C} = \{x \in \R^J: \down_C <  d(x,\lcdc ) < \up_C \},   
\]
 there exists a function $g_C: \domain \ra \R$, that satisfies the following four properties. 
\begin{enumerate} 
\item 
$g_C \in {\cal C}^\infty (\domain)$; 
\item 
there exists $\theta_C > 0$ such that  
\be
\label{gradgc} 
\lan \nabla g_C(x), p \ran \leq -\theta_C \quad \quad \mbox{ for } p \in K_C^{\delta_C/3} 
  \mbox{ and } x \in \domain ; 
\ee
\item 
for every $j,k \in \{1, \ldots, J\}$, 
\[ \sup_{x \in \domain}  
 \left| \dfrac{\partial^2 g_C (x)}{\partial x_j \partial x_k}\right| < \infty ; 
\]
\item 
$\sup_{x \in \Lambda_C} \left( |g_C(x)  - d(x,L_{C,\delta_C})| \vee \left(|\nabla g_C(x)| - 1 \right) \right)\leq \tilde{\ve}_C$. 
\end{enumerate} 
\end{lemma} 
\begin{proof}  Fix $C \in {\cal P}({\cal U})$,  $0 < \down_C < \up_C < \infty$ and $\tilde{\ve}_C > 0$. 
  For ease of notation, for the rest of this proof we will usually just write $\newlc$, $\Lambda$, $\up$ and $\down$
 for $\lcdc$, $\domain$, $\up_C$ and $\down_C$, respectively.  
Define $\newg:\R^J \ra \R_+$  to be $\newg(\cdot) \doteq d(\cdot, \newlc)$, the 
 distance function to the cone $\newlc$.  The proof is comprised of four steps. 
In the first step, we collect some known properties of $\newg$. 
In the second step, we establish gradient properties of the type (\ref{gradgc}) for $\newg$ on $\domain$, and 
 in the third step, we obtain bounds on the growth of the second derivatives of $\newg$ on a subset of  $\domain$.  
In the last step we show that there exists a  smooth approximation $g_C$ of $\newg$ that 
satisfies the conditions of the theorem. \\  
{\bf Step 1.} 
Let $P_{\newlc}:\R^J \ra \newlc$ be the  metric projection onto the cone $\newlc$ (which assigns to each point $x \in \R^J$ 
the point on $\newlc$ that is closest to $x$), and let $H_C \doteq \{x \in \R^J: P_{\newlc} (x) = 0 \}$ be the set of points that get projected 
to the vertex $0$ under the map $P_{\newlc}$. 
Since $\newlc$ is a closed convex set, it is well-known     
 (and seems to have been first proved in \cite[p.\ 286]{mor})  that on $\R^J\sm \newlc$, 
  $\newg$ is ${\cal C}^1$, $P_{\newlc}$ is Lipschitz with constant $1$, $\newg(\cdot) = |x - P_{\newlc}(\cdot)|$ and 
\begin{equation}
\label{nablagc}
 \nabla \newg(x) = \dfrac{ x - P_{\newlc}(x)}{|x - P_{\newlc}(x)|} \quad \mbox{ for } x \in \R^J \sm \newlc. 
\end{equation}  
We now argue that  $\newg$ is (at least) ${\cal C}^3$ and $P_{\newlc}$ is (at least) ${\cal C}^2$ on  $\R^J \sm [\newlc \cup \partial H_C]$.  
Indeed, first note that the fact that $\newg(x) = |x|$ and $P_{\newlc}(x) = 0$ for $x$ in the interior of $H_C$ implies that both functions 
are  $C^\infty$ on $H_C^\circ$. 
Next, observe that since the boundary of $\kcdc$ is ${\cal C}^\infty$ (by construction), the boundary of $\newlc$ is 
also ${\cal C}^\infty$ everywhere except at the vertex $0$.  
Theorem 2 of \cite{hol} asserts that for $p \geq 1$, if the boundary of $\newlc$ is ${\cal C}^p$ 
in a neighbourhood of a point $x \in \partial \newlc$ then the distance function $\newg$ is ${\cal C}^p$ and the 
projection $\P_{\newlc}$ is   ${\cal C}^{p-1}$ in a neighbourhood of the open normal ray to $\newlc$ at the point $x$.  
In particular, this guarantees that $\tilde{g}_C$ is (at least) ${\cal C}^3$ and $P_{L_\delta}$ is (at least) 
${\cal C}^2$ on $\R^J\sm [L_\delta \cup \partial H_C]$.

\noi 
{\bf Step 2.}  
 For any $x \in \R^J\sm \newlc$, because $x - P_{\newlc}(x)$ is an outward normal to 
$\newlc$ at $P_{\newlc}(x)$ and $\newlc$ is convex, we have 
\be
\label{ineq-proj}
\lan x - P_{\newlc} (x), w -  P_{\newlc} (x) \ran \leq 0 \quad \mbox{ for every } w \in \newlc. 
\ee
Since $\newlc$ is a convex cone with vertex at the origin that contains $\kcdc$, we have 
$P_{\newlc}(x) + p \in \newlc$ for $p \in \kcdc$. Setting $w = P_{\newlc}(x) + p$ 
in the last inequality we see that $\lan x - P_{\newlc}(x), p \ran \leq 0$ for every 
$p \in \kcdc$. 
Since $K_C^{\delta_C/2} \subset (K_{C,\delta_C})^\circ$ by (\ref{def-kcdc}) and  
$K_C^{\delta_C/2}$ is compact,  
  combining this with the expression for $\nabla \newg$ given in (\ref{nablagc}), 
 we conclude that there must exist $\theta_C$ that satisfies 
\be 
\label{newgradgc} 
\lan \nabla \newg (x), p \ran \leq -2 \theta_C \quad \quad \mbox{ for } p \in K_C^{\delta_C/2} 
  \mbox{ and } x \in \R^J \sm L_{C,\delta_C}^{\down_C/2}.   
\ee

\noi 
{\bf Step 3.} 
Let $\{\delta_{ij}, i,j \in \{1, \ldots, J\}\}$ be the Kronecker delta function: $\delta_{ij} = 1$ if $i = j$ and 
$0$ otherwise. 
For $x \in \R^J \sm \left[\newlc \cup \partial H_{C}\right]$, using the expression (\ref{nablagc}) we obtain  
\[ 
\ba{l}
 \dfrac{\partial^2 \newg}{\partial x_j \partial x_i} (x)  
 =   
\dfrac{\partial}{\partial x_j} \lan \nabla \newg (x), e_i \ran \\
 =  \dfrac{\partial}{\partial x_j} \left( \dfrac{x_i - \lan P_{\newlc}(x), e_i \ran}{|x - P_{\newlc}(x)|} 
 \right) \\
 =  \dfrac{1}{|x - P_{\newlc}(x)|} \dfrac{\partial }{\partial x_j} (x_i - \lan P_{\newlc}(x), e_i \ran) 
 + (x_i - \lan P_{\newlc}(x), e_i \ran ) \dfrac{\partial}{\partial x_j} \left(\dfrac{1}{|x - P_{\newlc}(x)|} \right)
\\
 = \dfrac{1}{\newg(x)}  \left(\delta_{ij} - \dfrac{\partial }{\partial x_j}\lan P_{\newlc}(x), e_i \ran \right)
- (x_i - \lan P_{\newlc}(x), e_i \ran )\dfrac{\lan \nabla \newg(x), e_j \ran}{\newg^2(x)}. 
\ea
\]
Observing that the maps $F:x \mapsto \lan x, e_i \ran$ and $P_{\newlc}$ are non-expansive (see, for example, Section 4 of \cite{hol} for the latter result), denoting the differential of $P_{\newlc}$ by $DP_{\newlc}$  and using the trivial relations 
$|x_i - \lan P_{\newlc} (x), e_i \ran| \leq |x - P_{\newlc}(x)| = \newg(x)$ and $|\nabla \newg(x)| = 1$, we obtain 
the following bound: for $x \in \R^J\sm [\newlc \cup \partial H_C]$, 
\be
\label{der2gcbound}
 \left| \dfrac{\partial^2 \newg}{\partial x_j \partial x_i} (x) \right|  \leq
 \dfrac{1 + |\nabla F(P_{\newlc}(x))| |DP_{\newlc}(x)|}{\newg(x)}  + \dfrac{1}{\newg(x)} 
\leq  \dfrac{3}{\newg(x)}.  
\ee

\noi 
{\bf Step 4.} While the constructed function $\newg$ satisfies most of the desired properties, it is not 
${\cal C}^2$ in a neighbourhood of $\partial H_C$.  We shall use an approximation argument to overcome this problem.  
Since $\newg \in {\cal C}^1 (\R^J\sm \newlc)$,   by Lemma 2 on page 19 of 
\cite{bur}, $D^\alpha \newg = D^\alpha_w \newg$ on $\R^J \sm \newlc$ 
 whenever $|\alpha| = 1$. 
When $|\alpha| = 2$, we only know that 
$D^\alpha \newg$ exists on the open set $\R^J \sm [\newlc \cup \partial H_C]$ (this was established 
in Step 1). 
We now claim that $D^\alpha \newg$, $|\alpha| = 2$, serve as {\em weak} derivatives of second order on the larger domain   
$\R^J \sm \newlc$.  
Although for $|\alpha| = 2$, $D^\alpha \tilde{g}_C$ is defined only on $\R^J \sm [ \newlc \cup \partial H_C]$, 
since $\partial H_C$ is a set of Lebesgue measure zero and 
weak derivatives are only defined upto sets of measure zero, 
the claim makes sense. 
Moreover, 
in order to establish the claim, it clearly suffices to show 
that  weak derivatives of second order exist on $\R^J \sm \newlc$. 
From the expression (\ref{nablagc}) for $\nabla \newg$ and the fact that 
$P_{\newlc}$ is Lipschitz, we know  that $D^\alpha g$ is 
absolutely continuous on $\R^J \sm \newlc$ when $|\alpha| = 1$. 
A simple integration by parts argument
 (along the lines, for example, of  Lemma 9 of page 34 of \cite{bur}),  
 combined with the fact that $D^\alpha \newg = D_w^\alpha \newg$  when $|\alpha| = 1$, 
then shows that $D_w^\alpha \newg$ exists on $\R^J \sm \newlc$ for $|\alpha| = 2$ and the claim follows.

The closure of $\Lambda$ lies in $\R^J \sm \newlc^{\down/2}$.  
The claim just proved,  along with (\ref{der2gcbound}), then implies that when $|\alpha| = 2$, 
\be
\label{sec-deriv}
 \mbox{ess}\sup_{x \in \Lambda} | D^\alpha_w \newg (x) |  = \sup_{x \in \Lambda \sm \partial H_C} |D^\alpha \newg (x) | 
\leq \dfrac{3}{\inf_{x \in \Lambda}\newg(x)} = \dfrac{6}{\down}. 
\ee
Therefore, for $|\alpha| = 2$, the essential supremum of $D^\alpha \newg = D^\alpha_w \newg$ on 
$\Lambda$ is finite. 
 Furthermore, $\newg$ is also uniformly bounded on $\Lambda$ (by ${\up}$).  
Thus $\newg$ lies in 
 $W_\infty^2(\Lambda)$, the Banach space of uniformly bounded functions on $\Lambda$ for which 
all weak derivatives of second order exist, equipped with the norm 
\[ ||f||_{W_\infty^2(\Lambda)} = \mbox{ess}\sup_{x \in \Lambda} |f(x)| + \sum_{\alpha \in \Z_+^J:|\alpha| = 2} 
\mbox{ess}\sup_{x \in \Lambda} 
|D^\alpha_w f(x)|. 
\]
Although the space $C^\infty(\Lambda)$ is not dense in $W_\infty^2(\Lambda)$, 
by Theorem 1 on page 48 of \cite{bur}   there exists a sequence $\{f_k\}$ of $C^\infty (\Lambda)$ functions such that for $k \in \N$, 
the following three properties hold: 
\[ 
\begin{array} {ll}
1. \sup_{x \in \Lambda} |f_k (x) - \newg (x)| \leq \dfrac{1}{k} ; \\
2.  \sup_{x \in \Lambda} |D^\alpha f_k (x) - D^\alpha \newg (x)| \leq \dfrac{1}{k} \quad \quad \quad  & \mbox{ if } |\alpha| = 1 ; \\ 
3. \left|\sup_{x \in \Lambda} |D^\alpha f_k(x)| - \sup_{x \in \Lambda} |D^\alpha \newg(x)|\right| \leq \dfrac{1}{k} \quad \quad \quad & \mbox{ if } |\alpha| = 2. 
\end{array}
\]
Combining this with (\ref{newgradgc}),  (\ref{sec-deriv}) 
and the fact that $|\nabla \tilde{g}_C(x)|=1$ for all $x \in \Lambda$ (see (\ref{nablagc})), it is clear that there 
exists a large enough integer $k > 1/\tilde{\ve}_C$ such that $g_C \doteq f_k$  
satisfies the properties of the lemma.   
\end{proof}

\begin{figure}
\centerline{\epsfxsize1in}\epsfbox{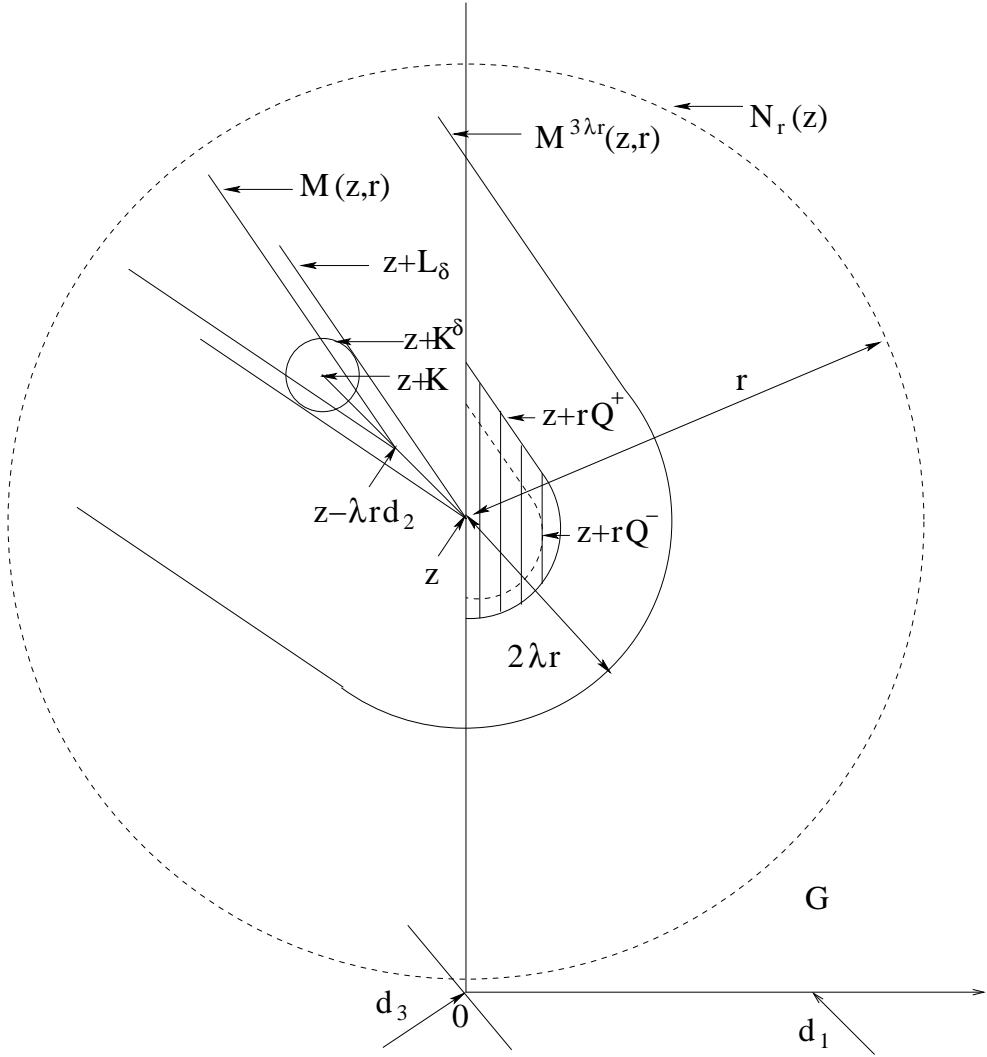} 
\caption{Construction of $g_{z,r}$ 
for the 2-dimensional GPS ESP}
\label{fig-lemma} 
\end{figure}

In the proof of Theorem \ref{lem-locgfn} given below, we show that  
a family of functions $\{g_{z,r}, z \in {\cal U}, r \in (0,r_z)\}$ 
that satisfy the necessary properties 
 can be obtained as localized, scaled versions of the  
functions $g_C, C \in {\cal P}({\cal U})$, constructed in the 
last lemma. 
The next two lemmas  will be used in the proof in order to characterize the geometry of 
the supports of these localized functions $\{g_{z,r}\}$ (see Figure \ref{fig-lemma} for 
an illustration).

\begin{lemma}
\label{lem-msets}
For $C \in {\cal P}({\cal U})$, let $\delta_C > 0$ satisfy  (\ref{ineq-min}) and 
let $q_C$ be a unit vector in $K_C$.  Then there exists $\lambda_C > 0$ such that for 
every $z \in F_C$ and $r \in (0,r_z)$, the set 
\be
\label{def-mset}
\mczr \doteq z + \lambda_C r q_C + \lcdc
\ee
satisfies the following three properties: 
\begin{enumerate} 
\item 
there exists $\eta_C \in (0, \lambda_C)$ such that $M^{\eta_C r} (z,r) \cap G = \emptyset$;  
\item 
$M_C^{3 \lambda_C r} (z,r) \cap G \subset N_r(z) \cap G$; 
\item 
for every $\ve > 0$, $\tilde{z} \in F_C$ and $\tilde{r} \in (0,r_{\tilde{z}})$, 
\[ x \in \left[ M_C^{\ve r} (z,r) - z \right] \Leftrightarrow \dfrac{\tilde{r} x}{r} \in 
\left[ M_C^{\ve \tilde{r}} (\tilde{z}, \tilde{r} ) - \tilde{z} \right]. 
\]
\end{enumerate} 
\end{lemma} 
\begin{proof} 
Fix $C \in {\cal P}({\cal U})$, $\delta_C > 0$ and $q_C \in K_C$ as in the statement of the lemma. 
Since $G$ is a convex polyhedron 
and $n_i, i \in C,$ are inward normals to $G$ at $z$, 
any $x \in G$  satisfies 
\[ 
 \min_{i \in C} \lan n_i, x - z \ran \geq 0. 
\]
Along with (\ref{eq-theta}), this shows  that 
the boundary of $G$ separates $z + L_{C,\delta_C}$ from the 
interior of $G$.  More precisely, for $z \in F_C$, we see that 
\be
\label{int-vert}
 (z + L_{C,\delta_C}) \cap G  = \{z\}  
\ee
and so the (minimal) angle $\phi_C$ between the closed convex cone $z + L_{C,\delta_C}$ and 
the closed polyhedron $G$ at $z$ satisfies  $0 < \phi_C < \pi/2$. 
Therefore for any $r \in (0,r_z)$, 
\[ \left\{x \in \R^J: d(x, z + L_{C,\delta_C}) \leq \dfrac{r}{2} \sin \phi_C  \right\} \cap \partial G \subset N_r(z) \cap \partial G. \]
Define $\lambda_C \doteq \sin \phi_C/6  \in (0,1/3)$.  Then the last display implies that 
for every $r \in (0,r_z)$,  
\be
\label{lambda}
\left\{ x \in \R^J: d (x, z + L_{C, \delta_C}) \leq  3 \lambda_C r \right\} \cap 
G \cap \partial N_r (z) = \emptyset. 
\ee
It is clear from the definition that $\lambda_C$ is independent of $z \in F_C$ 
(since the angle $\phi_C$ does not depend on $z \in F_C$) and thus 
(\ref{lambda}) holds for all $z \in F_C$ and $r \in (0,r_z)$.

Now let $M_C(z,r)$ be the cone $\lcdc$ shifted to have its vertex at $z + \lambda_C r q_C$, 
as defined in (\ref{def-mset}).  
Then, $z + \lambda_C r q_C$ lies inside the cone $z + L_{C,\delta_C}$ because  $q_C \in K_C \subset L_{C,\delta_C}$ and $\lambda_C r > 0$. 
Since $z + L_{C,\delta_C}$ is a convex cone, it follows that  
\[ M_C(z,r) \subset z + L_{C,\delta_C} \quad \mbox{ and } \quad M_C(z,r) \cap \{z\} = \emptyset. \]
Hence we infer that 
\[ M_C^{3 \lambda_C r}(z,r) \doteq \left\{ x: d(x,M_C(z,r)) \leq 3 \lambda_C r\right\} \subset \left\{x:d(x, z+L_{C,\delta_C}) \leq 3 \lambda_C r\right\}.
\]  
The last three displays together with (\ref{int-vert}) show that 
\be
\label{lambda2}
M_C(z,r) \cap G \cap N_r(z) = \emptyset  \quad \mbox{ and } \quad M_C^{3 \lambda_C r} (z,r) \cap G \cap \partial N_{r} (z) = \emptyset. 
\ee 
Together with the fact that $z \in M_C^{3\lambda_C r}(z,r)$, the second equality above 
 establishes property 2 of the lemma. 
Moreover, since $M_C(z,r)$ and $G \cap N_r(z)$ are both closed, 
from the first relation of the last display and property 2, it is clear 
that property 1 of the lemma also holds. 

It only remains to prove property 3.  For 
$z \in F_C$, 
$x \in  M_C^{\ve r} (z,r) - z$ if and only if there exists 
$v$ with $|v| \leq \ve r$ and $u \in L_{C,\delta_C}$ such that 
$x = v + \lambda_C r q_C + u$. 
This implies that 
$\tilde{r}x/r = \tilde{r}v/r + \lambda_C \tilde{r} q_C + \tilde{r}u/r 
= \tilde{v} + \lambda_C \tilde{r} q_C + \tilde{u}$, where 
$|\tilde{v}| = |\tilde{r}v/r| \leq \ve \tilde{r}$ and 
$\tilde{u} = \tilde{r} u/r$ lies in  $L_{C,\delta_C}$ since 
 $L_{C,\delta_C}$ is a cone with vertex at $0$, which concludes the proof. 
\end{proof}

We now identify the sets $Q_C^-$, $Q_C^+$, $C \in {\cal P}({\cal U})$,  
that arise in Theorem \ref{lem-locgfn}.

\begin{lemma}
\label{lem-qsets}
Fix $C \in {\cal P}({\cal U})$, let 
$\{M_C(z,r), z \in F_C, r \in (0,r_z)\}$ be as in Lemma \ref{lem-msets}. 
For any  
$\alpha \in (0,1)$,  
there exist closed, convex sets $Q_C^-$ and $Q_C^+$ such that 
$Q_C^- \subset (Q_C^+)^\circ \subset N_1(0)$ and for every $z \in F_C$ and $r \in (0,r_z)$, 
\be
\label{def-qset}
 M_C^{2 \lambda_C r} (z,r) \cap G = z + r Q_C^+  \quad \quad \quad 
M_C^{(1+\alpha)\lambda_C r} (z,r) \cap G = z + r Q_C^- 
\ee
and 
\be
\label{rint}
z \in \rint [(z + r^\prime Q_C^-) \cap \partial G]. 
\ee 
\end{lemma}
\begin{proof} 
Fix $C \in {\cal P}({\cal U})$, choose a particular $z^\prime \in F_C$ and $r^\prime \in (0,r_{z^\prime})$, 
and define $Q_C^-$ and $Q_C^+$
 using (\ref{def-qset}),   with $z$ and $r$ replaced by $z^\prime$ and $r^\prime$, respectively.  
Since the set $M^{\theta \lambda_C r^\prime}_C (z^\prime,r^\prime)$, for $\theta = \alpha, 2$, and the set 
$G$ are closed and convex, clearly $Q_C^+$ and $Q_C^-$ are also closed and convex. 
The inclusion $Q_C^- \subset (Q_C^+)^\circ$ also follows directly from the definition and the fact that 
$1 + \alpha < 2$.   Moreover, 
since property 2 of Lemma \ref{lem-msets} implies  
$M_C^{2 \lambda_C r^\prime} (z^\prime,r^\prime) \cap G \subset N_{r^\prime}(z^\prime)$, 
 it follows that $z^\prime + r^\prime Q_C^+ \subset N_{r^\prime}(z^\prime)$ or, equivalently, that 
$Q_C^+ \subset N_1(0)$. 

Property 3 of Lemma \ref{lem-msets}, the fact that $G$ is a cone 
and the property that  for $z \in F_C$ and $r < r_z$, 
$d({\cal V},M_C^{2 \lambda_C r} (z,r) \cap G) > d({\cal V},\partial N_{r}(z)) > 0$ (where the latter inequality 
holds since $r <  r_{z} < d(z,{\cal V})$), then show that (\ref{def-qset}) is applicable for every $z \in F_C$ and $r \in (0,r_z)$. 
In addition, since $d(z,M_C(z,r)) = \lambda_C r$ and $\alpha > 0$, it follows that for every $r \in (0,r_z)$, 
\[ z \in \rint [ M_C^{(1+ \alpha)\lambda_C r} (z,r) \cap \partial G] = \rint[(z + r Q_C^-) \cap \partial G], \]
which proves (\ref{rint}).  
\end{proof}

The sets $Q_C^+$, $Q_C^-$ and $M_C(z,r)$ for the case of the two-dimensional GPS ESP are illustrated 
in Figure \ref{fig-lemma}.  
We now combine the above results to establish Theorem \ref{lem-locgfn}. 
For $C \in {\cal P}({\cal U})$, let $\delta_C > 0$ satisfy (\ref{ineq-min}) 
and let $q_C \in K_C$, $\eta_C \in (0,1)$,  $\lambda_C \in (\eta_C/2,\eta_C)$,  
and $\{M_C(z,r), z\in F_C, r \in (0,r_z)\}$ be as in Lemma \ref{lem-msets}. 
Choose $\alpha \in (0, 1/3)$, 
and let $Q_C^-$, $Q_C^+$ be the corresponding sets described in Lemma \ref{lem-qsets}. 
Lastly, with the choice of $\down_C = \eta_C$, $\up_C = 3 \lambda_C$ and 
$\tilde{\ve}_C = \alpha \lambda_C/2$, 
let $g_C$, $\Lambda_C$ and $\theta_C > 0$ be as in Lemma \ref{lem-gc}. 

   For $z \in {\cal U}$ and $r \in (0,r_z)$, we now construct $g_{z,r}$ as a suitably scaled and localized 
version of the function $g_{I(z)}$, which vanishes on ${\cal V} = \{0\}$, while still maintaining 
the smoothness and derivative properties of $g_{I(z)}$.   First define $\tilde{\alpha} \doteq 1 + 3\alpha/2$ and 
$\tilde{\kappa} \doteq 2 - \alpha/2$ and note that $\tilde{\alpha} < \tilde{\kappa}$.  
 Then let  $h_C$  be a  non-decreasing ${\cal C}^\infty(\R_+)$ function that satisfies 
\rear{1.2}
\be
\label{hprops}
\ba{rcll} 
h_C^\prime(s) &  \geq & 0 & \mbox{ for } s \in [0,\infty) \\
h_C^\prime(s)& \geq & 1 & \mbox{ for } s \in (0,\tilde{\alpha} \lambda_C) \\
h_C(s) & = & 0  & \mbox{ for } s \in [ \tilde{\kappa} \lambda_C,\infty). 
\ea
\ee
 For  $C \in {\cal P}({\cal U})$, if $z \in F_C$ and $r \in (0,r_z)$, define 
\be
\label{def-omegac}
\Omega_C(z,r) \doteq \left( M_C^{3\lambda_C r} \right)^\circ \sm M_C^{\eta_C r} \quad \mbox{ for }  z \in F_C \mbox{ and } r \in (0,r_z)
\ee
and let  $g_{z,r}: \R^J \ra \R_+$ be given by  
\be
\label{def-gzr} 
g_{z,r} (x) \doteq \left\{ 
\ba{rl}
r h_C \left( g_{C} \left( \dfrac{x - \lambda_C q_C - z}{r}\right) \right)   & \mbox{ if } x \in \Omega_{C}(z,r), \\
0 & \mbox{ otherwise. } 
\ea
\right. 
\ee
We show below that the functions $\{g_{z,r}, z \in F_C, r \in (0,r_z)\}$ above 
satisfy the properties listed in Theorem \ref{lem-locgfn}. \\

\noi 
{\bf Proof of Theorem \ref{lem-locgfn}. } 
Define 
\be
\label{def-aprime}
 A^\prime \doteq \max_{C \in {\cal P}({\cal U})} \sup_{s \in [0,3 \lambda_C]} \left[\left|h_C(s)\right| \vee \left| h_C^\prime(s)\right| \vee \left| h_C^{\prime \prime}(s) \right| \right] \vee 2 
\ee
and  
$\theta \doteq \min_{C \in {\cal P}({\cal U})} \theta_C$ and note that $A^\prime < \infty$ and $\theta > 0$. 
Since by (\ref{hprops}), $h_C(s) = 0$ for $s \geq 3 \lambda_C \geq \tilde{\kappa} \lambda_C$, 
 the definition of $A^\prime$ in (\ref{def-aprime})  immediately implies property 3 of the theorem.

Fix $C \in {\cal P}({\cal U}), z \in F_C$ and $r \in (0,r_z)$, let $g_{z,r}$ be defined as in (\ref{def-gzr}) 
and for the rest of the proof, write $h$ for $h_C$, $\newlc$  for $L_{C,\delta_C}$, $q$ for $q_C$,  $\eta$ for $\eta_C$, 
 $\lambda$ for $\lambda_C$, 
$\tilde{\lambda}$ for $\tilde{\lambda}_C$,  $\tilde{\eta}$ for $\tilde{\eta}_C$, $\tilde{\ve}$ for $\tilde{\ve}_C$, 
$\Lambda$ for $\domain$ and 
$\Omega(z,r)$ for $\Omega_C(z,r)$.  
Note that since $\tilde{\lambda} = 3 \lambda$ and $\tilde{\eta} =  \eta$, 
\be
\label{rel-1}
x \in \Omega(z,r) \Rightarrow \dfrac{x - \lambda q - z}{r} \in \left(L_C^{\tilde{\lambda}}\right)^\circ 
\sm L_C^{\tilde{\eta}} = \Lambda. 
\ee
This guarantees that whenever $x \in \Omega (z,r)$,  $(x - \lambda q - z)/r$ lies in the domain $\Lambda$ 
of $g_C$ and so the function $g_{z,r}$ is well-defined.  

Now, observe that property 4 of Lemma \ref{lem-gc}, the identity $\tilde{\kappa} + \tilde{\ve}/\lambda = \tilde{\kappa} + \alpha/2 = 2$, 
the definitions (\ref{hprops}), (\ref{def-omegac}), (\ref{def-gzr}) and   (\ref{def-qset}) 
of $\Omega(z,r)$, $h$, $g_{z,r}$ and $Q^+_C$, respectively, 
 and properties 1 and 2 of Lemma \ref{lem-msets}, when combined, yield the relation   
\be
\label{rel-2}
\ba{l}
\supp \left[g_{z,r}\right] \cap G  \\
\quad \quad   \ds \subseteq 
 \left\{ x \in \R^J: g_C \left(\dfrac{x - \lambda q - z}{r} \right) \leq \tilde{\kappa} \lambda \right\} \cap G \\
\quad \quad   \ds \subset \left\{x \in \R^J: d\left( \dfrac{x - \lambda q - z}{r}, \newlc \right) 
 \leq \tilde{\kappa} \lambda + \tilde{\ve} \right\} \cap G \\
\quad \quad   \ds    =    \left\{x \in \R^J: d(x - \lambda q - z, \newlc) \leq 2 \lambda r \right\} \cap G\\
\quad \quad   \ds   =    M_C^{2 \lambda r} (z,r) \cap G \\
\quad \quad   \ds  =   z + rQ^+_C \\
\quad \quad   \ds   \subset   \Omega (z,r) \cap G. 
\ea
\ee 
The inclusion 
$\supp \left[g_{z,r}\right] \cap G  \subset \Omega(z,r) \cap G$, 
the relation (\ref{rel-1})  and the fact that 
 $h \in {\cal C}^\infty(\R_+)$ and $g_C \in {\cal C}^\infty (\Lambda)$,
imply that $g_{z,r} \in {\cal C}^\infty(\R^J)$. 
In addition, since   $\Omega (z,r) \cap G 
  \subset    N_r(z) \cap G$, 
the inclusion also implies  that for $r < r_z$, 
\be
\label{rel-4} 
x \in \Omega(z,r) \cap G \quad \Rightarrow \quad d^1(x) \subseteq d^1(z) 
\ee
 and  because $N_r(z) \cap {\cal V} = \emptyset$, that 
$g_{z,r} (x) = \nabla g_{z,r} (x) = 0$ for $x \in {\cal V}$. 
The last three assertions show that $g_{z,r}$ satisfies properties 1 and 2 of 
Theorem \ref{lem-locgfn}.

 Now for $x \in \R^J$, 
\be
\label{rel-5}
 \nabla g_{z,r} (x) = h^\prime \left( g_C \left( \dfrac{x - \lambda q - z}{r} \right) \right) 
\nabla g_C\left( \dfrac{x - \lambda q - z}{r} \right). 
\ee
Property 2 of Lemma \ref{lem-gc}, along with  the fact that $z \in F_C$ implies 
$d^1(z) \subset \{aK_C, a \leq -1\}$, shows that 
\[ \lan \nabla g_C(y), d \ran  \geq \theta  \quad \forall d \in d^1(z), y \in \Lambda. \]
When combined with (\ref{rel-1}) and (\ref{rel-4}), this in turn shows that 
\be
\label{rel-6}
\left\lan \nabla g_{C} \left( \dfrac{x - \lambda q - z}{r} \right), d \right\ran \geq \theta \quad 
\forall d \in d^1(x),  x \in \Omega(z,r) \cap \partial G,  
\ee  
while property (\ref{rel-2}) ensures that 
\[ \nabla g_C \left( \dfrac{x - \lambda q - z}{r} \right) = 0 \quad \mbox{ for } x \in G \sm \Omega (z,r). 
\]
Property 5 of the theorem then follows from the last display, (\ref{rel-5}), (\ref{rel-6}) and  the 
first inequality in (\ref{hprops}). 
On the other hand, using the identity $\tilde{\alpha} = 1 + \alpha + \alpha/2 = 1 + \alpha + \tilde{\ve}/\lambda$,
 and once again invoking property 4 of 
Lemma \ref{lem-gc} and recalling the defintion (\ref{def-qset}) of $Q_C^-$, we observe that  
\be
\label{rel-3}
\ba{l} 
\left\{ x \in \R^J: g_C \left( \dfrac{x - \lambda q  - z}{r}\right)  \leq  \tilde{\alpha} \lambda \right\} \cap G\\
\quad \quad = \left\{ x \in \R^J: g_C \left( \dfrac{x - \lambda q - z}{r}\right)  \leq (1+ \alpha) \lambda + \tilde{\ve}
 \right\} \cap G\\
\quad \quad  \supset  
 \left\{ x \in \R^J: d \left(\dfrac{x - \lambda q - z}{r}, \newlc\right) \leq (1+\alpha) \lambda \right\} \cap G\\
\quad \quad  = 
 M_C^{(1+\alpha)\lambda r} (z,r) \cap G \\
\quad \quad = z + r Q_C^-.  
\ea
\ee 
When combined with the second inequality in (\ref{hprops}), (\ref{rel-5}) and (\ref{rel-6}), 
this implies property 4 of the theorem. 


Differentiating $g_{z,r}$ twice and using the chain rule,  we see that for $x \in \Omega (z,r)$, 
\[
\ba{l}
\dfrac{\partial^2 g_{z,r}}{\partial x_i \partial x_j} (x) \\
 =  \dfrac{1}{r} h^{\prime \prime} \left(g_{C}\left(\dfrac{x - \lambda q r - z}{r} \right) \right)
\dfrac{\partial g_{C}}{\partial x_i} \left(\dfrac{x - \lambda q r - z}{r} \right)  \dfrac{\partial g_{C}}{\partial x_j} 
 \left(\dfrac{x - \lambda q r - z}{r} \right) \\
\quad + h^\prime \left(g_{C}\left(\dfrac{x - \lambda q r - z}{r} \right)  \right)
\dfrac{\partial^2 g_{C}}{\partial x_i \partial x_j} \left(\dfrac{x - \lambda q r - z}{r} \right). 
\ea
\]
Let $L^\prime < \infty$ be an upper bound (independent of $j,k$) for the second derivatives in property 
3 of Lemma \ref{lem-gc}. 
By property 4 of Lemma \ref{lem-gc}, we know that $\sup_{x \in \Lambda} |\nabla g_C(x)| \leq 1 + \tilde{\ve} \leq 2$. 
 Along with (\ref{rel-2}) and the definition of $A^\prime$, 
this implies the bound 
\[ \sup_{x \in G} \left| \dfrac{\partial^2 g_{z,r}}{\partial x_i \partial x_j} (x) \right| = 
\sup_{x \in \Omega(z,r) \cap G}  \left| \dfrac{\partial^2 g_{z,r}}{\partial x_i \partial x_j} (x)   \right| 
\leq \dfrac{4A^\prime}{r} + A^\prime L^\prime \leq \dfrac{4A^\prime + A^\prime L^\prime r_z}{r}. 
\]
Hence properties 6 and 7 of the theorem are satisfied with $A(z) \doteq A^\prime (4 + L^\prime r_z)$ 
and $\tilde{A}_R \doteq A^\prime (4 + L^\prime R)$. 
This completes the proof of the theorem. 
\ink \\

\subsection{Proof of Lemma \ref{lem-cover}}
\label{subs-ap2}

In the last section, we constructed a family of $C^\infty(\R^J)$ 
 functions $\{g_{z,r}: z \in {\cal U}, r \in (0,r_z)\}$,  
with each $g_{z,r}$ satisfying certain gradient and second derivative properties in a  
neighbourhood of $z$. 
Lemma \ref{lem-cover} below shows that a countable set $S$ and scalars $\rho_z, z \in S$,  
can be chosen such that any $x \in G$ lies in the support of at most a finite number 
(independent of $x$) of functions $g_{z,\rho_z}, z \in S$, and  the 
sets $z + \rho_z Q_{I(z)}^-, z \in S$, cover $N_r({\cal V})\sm {\cal V}$. 
This ensures that the function  $g$ defined as the countable sum  $\sum_{z\in S} g_{z,\rho_z}$ 
is well-defined on $G$ and  satisfies the necessary derivative conditions on all of $G$. 
Although the notation in the proof of Lemma \ref{lem-cover} is a bit involved, 
the basic idea behind the proof is quite simple. 
One first identifies a finite number of points $z$ on 
a hyperplane $H_s$ a  distance $s$ 
away from ${\cal V} = \{0\}$, and  associated scalars $\rho_z$ such that the union of the  corresponding  
neighbourhoods covers the intersection of $G$ 
with the fattening of a hyperplane $H_s$ (see the claim below). 
Using scaling arguments 
one then identifies a corresponding finite number of points on 
a suitable translation of that hyperplane along its 
normal. 
The covering is then obtained by taking the union of the associated (finite number) of  
neighbourhoods in each of a  
countable number of hyperplanes $H_{s_j}, s_j \downarrow 0$.  \\

\noi
{\bf Proof of Lemma \ref{lem-cover}.} 
By assumption, $G$ is a convex cone with vertex at the origin 
and  angle  less than $\pi$. 
Thus there exists $v \in G^{\circ}$ with $|v| = 1$, such that the maximum angle 
$\omega$ between $v$ and any  $x \in \partial G$ is less than $\pi/2$. 
This implies that there exists $\zeta = \cos \omega > 0$ such that 
\be
\label{zeta} 
 |x| < \dfrac{\lan x, v \ran}{ \zeta}  \quad \mbox{ for every } x \in G \sm \{0\}. 
\ee
Now for $s > 0$, let $H_s$ be the hyperplane defined by 
\[ H_s \doteq \{x \in \R^J: \lan x, v \ran = s \} \]
and for $0 \leq s < \tilde{s} < \infty$ define the ``slabs''  $H[s,\tilde{s}]$ and 
$ H(s,\tilde{s}]$ to be  
\[ H[s,\tilde{s}] \doteq \{x: s \leq \lan x,v \ran \leq \tilde{s} \}  \quad 
\mbox{ and  } \quad 
 H(s,\tilde{s}] \doteq \{x: s < \lan x, v \ran \leq \tilde{s} \}. \]
For $C \in {\cal P}({\cal U})$, let $Q_C^+$, $Q_C^-$ be the closed, bounded, convex sets with 
$0 \in (Q_C^-)^\circ \subset (Q_C^+ )^\circ$ specified in the statement of 
the lemma.  Given $0  < r < R < \infty$, choose 
$s_1 \in (0,R)$ such that 
\be
\label{def-s0}
  G \cap N_{r} (0) \subset G \cap H[0,s_1] \subset [N_{R} (0)]^\circ.
\ee
Also, define 
\be
\label{eq-kappa}
 \kappa \doteq \sup \{\rho > 0: z + \rho Q_{I(z)}^+ \subset 
N_{R}(0) \mbox{ for every } z \in H_{s_1} \cap G \}, 
\ee
where we recall that $I(z) \doteq \{i \in \bI: \lan z, n_i \ran = 0 \}$. 

We now show that  the lemma is a consequence of the following claim, 
and defer the proof of the claim to the end. \\
\noi 
{\it Claim.}  There exist $0 < \beta <  \gamma < 1$,   
 a finite set $S_1 = \{z^{(1)}_i, i = 1, \ldots, L\} \subset H_{s_1}$    
 with associated scalars $\rho_{z_i^{(1)}}$, $i = 1, \ldots, L$,  
such that the sets defined for $i = 1, \ldots, L$, by 
\be 
\label{qi} 
 Q_{z^{(1)}_i}^+ \doteq z_i^{(1)} + \rho_{z_i^{(1)}} Q_{I(z^{(1)}_i)}^+  \gap \mbox{ and } \gap 
Q_{z^{(1)}_i}^- \doteq z_i^{(1)} + \rho_{z_i^{(1)}} Q_{I(z^{(1)}_i)}^- 
\ee 
satisfy the following two properties: 
\begin{enumerate} 
\item 
$ H [s_1 (1 - \beta), s_1(1 + \beta)] \cap \partial G 
 \subset \left[ \cup_{z \in S_1} Q^-_z \right] \cap \partial G$;  
\item
$\left[ \cup_{z \in S_1} Q^+_z \right] \cap G 
\subset  H[s_1(1- \gamma), s_1(1+\gamma)]\cap G.$
\end{enumerate} 
Suppose the claim is true. 
Then for $k = 2, 3, \ldots$, and $i = 1, \ldots, L$, define
\be
\label{zidefs}
 z_i^{(k)} \doteq (1- \beta)^{k-1} z_i^{(1)}  \quad \mbox{ and } \quad  \rho_{z_i^{(k)}} \doteq (1-\beta)^{k-1} 
\rho_{z_i^{(1)}},   
\ee
 let $S_k \doteq \{z_i^{(k)}, i = 1, \ldots, L\}$ and define  
$Q_{z_i^{(k)}}^-$ and $Q_{z_i^{(k)}}^+$ as in (\ref{qi}), with $(1)$ replaced everywhere by $(k)$. 
   Then, since $G$ is a conical polyhedron  it is clear that for $i = 1, \ldots, L$, 
$I(z_i^{(k)}) = I(z_i^{(1)})$ and therefore we have 
\[ Q_{z_i^{(k)}}^- = (1-\beta)^{k-1} Q_{z^{(1)}_i}^-  \quad \mbox{ and } 
\quad Q_{z_i^{(k)}}^+ = (1-\beta)^{k-1} Q_{z^{(1)}_i}^+. 
\]
Combining this with the claim it follows that for 
$k \in \N$, 
\be
\label{incl5}
 H [s_1 (1 - \beta)^{k}, s_1(1 - \beta)^{k-1}(1+\beta)] \cap \partial G 
 \subset \left[ \cup_{z \in S_k} Q^-_z \right] \cap \partial G 
\ee
and 
\be
\label{incl6}
 \left[ \cup_{z \in S_k} Q^+_z \right] \cap G 
\subset  H[s_1(1-\beta)^{k-1}(1- \gamma), s_1(1-\beta)^{k-1}(1+\gamma)]\cap G. 
\ee  
Define $S \doteq \cup_{k \in \N} S_k$.  Since each $S_k$ has $L$ elements, $S$ is countable.

Let $n$ be the smallest integer such that 
\be
\label{n-ineq}
 n > \log [ (1 - \gamma)/(1 + \gamma)]/\log(1 - \beta). 
\ee 
Fix $k \in \N$, $z \in S_k$ and $x \in Q_z^+$ for this paragraph.  
We first show that   
\be
\label{hs-ineq}
 \tilde{z} \not \in \bigcup_{j=(k-n+1)\vee 1}^{k+n-1} S_j \quad \Rightarrow \quad x \not \in Q_{\tilde{z}}^+.  
\ee
Indeed, note that (\ref{incl6}) implies that  
\[ x \in  H[s_1(1-\beta)^{k-1}(1- \gamma), s_1(1-\beta)^{k-1}(1+\gamma)]\cap G. \]
If $i \geq k + n$ then $(i - 1) - (k-1)   \geq n$ and so 
(\ref{n-ineq}) yields the inequality 
\[ (1-\beta)^{i-1}(1+\gamma) < (1-\beta)^{k-1} (1-\gamma). \]
The inclusion  (\ref{incl6}) then implies that $x \not \in \cup_{y \in S_i} Q_y^+$. 
Likewise, if $1 \leq i \leq (k-n)\vee 1$ then for $k \geq n+1$, 
$(k - 1) - (i - 1) \geq n$, and so  (\ref{n-ineq}) implies that 
\[ (1-\beta)^{k-1} (1 + \gamma) < (1-\beta)^{i-1} (1 - \gamma). \]
Once again, we have  $x \not \in \cup_{y \in S_i} Q_y^+$ by  (\ref{incl6}) and  
so (\ref{hs-ineq}) follows.  Since each $S_i$ has $L$ elements, this establishes 
property 1 with $N = 2(n-1)L$. 
Next, observe that since $z \in S_k$, from (\ref{incl6}) and (\ref{zidefs}) we obtain  
\[ \dfrac{\lan z, \nu \ran}{\rho_z} = \dfrac{(1-\beta)^{k-1} s_1}{(1-\beta)^{k-1} \rho_{I(z)}^{(1)}}  \leq  \dfrac{s_1}{\rho_*} \]
where  
 $\rho_* \doteq \min_{z \in S_1} \rho^{z} > 0$. 
Moreover, since $x \in Q_z^+$ by (\ref{zeta}) and (\ref{incl6}), we have 
\[ d(x,{\cal V}) = |x| \leq \dfrac{\lan x, \nu \ran}{\zeta} \leq \dfrac{(1+\gamma)\lan z, \nu \ran}{\zeta}  
\leq \dfrac{(1+\gamma)s_1 \rho_z}{\rho_* \zeta}, \]
which shows that property 2 is satisfied with $\nu \doteq  (1+\gamma) s_1/\rho_* \zeta  < \infty$. 
Furthermore, property 3 is a simple consequence of the fact that since $Q_C^+ \subset N_1(0)$ for every 
$C \in {\cal P}({\cal U})$ and $\rho_z < r_z$, we have 
$x \in Q_z^+ \subset N_{r_z} (z)$.

Lastly, note that since $\beta \in (0,1)$, we have 
\[ 
\cup_{k \in \N} H [s_1 (1 - \beta)^{k}, s_1(1 - \beta)^{k-1}(1+\beta)] \cap  \partial G 
 =   H(0,s_1(1+\beta)] \cap \partial G,
\]
which is contained in  $[ \cup_{z \in S} Q_z^- ] \cap \partial G$ by (\ref{incl5}). 
Recalling that   ${\cal V} = \{0\}$ and using 
(\ref{def-s0}) 
we then obtain 
\[ {\cal U} \cap N_r({\cal V}) =  [\partial G \sm \{0\}] \cap N_{r}(0)  \subset 
 \partial G \cap H(0,s_1]  \subset \partial G \cap [ \cup_{z \in S} Q_z^- ]. 
\]
Also, the fact that $\rho_C < \kappa$  and the definition  (\ref{eq-kappa}) of $\kappa$ 
yield 
\[ [\cup_{z \in S} Q_z^+ ] \cap G \subset N_{R} (0) = N_R ({\cal V}). \]
This proves properties 4 and 5. 
To complete the proof of the lemma, 
it only remains to establish the claim. \\

\noi
{\em Proof of Claim.}   
For $j = 1, \ldots, J-1$, let ${\cal P}_j ({\cal U})$ be the collection of 
subsets corresponding to $j$-dimensional faces in $\partial G$: 
\[ 
{\cal P}_j ({\cal U}) \doteq \{C \in {\cal P}({\cal U}): 
\dim [F_C] = j\},  
\] 
where $\dim[A]$ represents the dimension of the 
affine hull of $A$. 
Note that ${\cal U} = \cup_{j=1}^{J-1} \cup_{C \in {\cal P}_j({\cal U})} F_C$. 
In order to establish the first property of the claim it therefore suffices to show that  
 for every $j = 1, \ldots, J-1$, there exist $\beta_j \in (0,1)$ 
and a finite set of points $S_1^j \subset H_{s_1} \cap [\cup_{C \in {\cal P}_j({\cal U})} F_C]$ 
and associated scalars $\rho_z \in (0,\kappa)$, $z \in S_1^j$,  such that 
\be
\label{induction}
 H[s_1(1 - \beta_j), s_1 (1 + \beta_j)]  \cap 
 \left[ \cup_{C\in {\cal P}_j(\cu)} F_{C} \right] 
 \subset \left[ \cup_{z \in \cup_{i=1}^j S_1^i} Q^-_z \right] \cap 
\left[ \cup_{C \in {\cal P}_j(\cu)} F_{C} \right]  
\ee
where, as usual, $Q_z^- \doteq z + \rho_z Q_{I(z)}^-$. 
Indeed then setting  $\beta \doteq \min_j \beta_j$ and $S_1 = \cup_{j=1}^J S_1^j$,  
the union of (\ref{induction}) over $j = 1, \ldots, J-1$ yields 
property 1 of the claim.

We shall prove (\ref{induction}) by induction on $j$. 
As explained below, it is easy to see that (\ref{induction}) holds when $j = 1$. 
For each $C \in {\cal P}_1(\cu)$, note that  
$H_{s_1} \cap F_C$ is equal to a point. 
Let $S^1_1 \doteq \{H_{s_1} \cap F_C: C \in {\cal P}(\cu)\}$
 be  the finite collection of these points as $C$ varies over 
${\cal P}({\cal U})$.  
Now define 
\be
\label{choice}
\rho^{(1)} \doteq \inf_{z \in S_1^1 } r_z > 0,  
\ee 
and for each $z \in S_1^1$, set $\rho_{z} \doteq \rho^{(1)}$ and  
 let $Q_z^-$ be defined as above. 
By assumption, $0 \in (Q_C)^\circ$ and so  $Q_{z}^- \cap F_C$ is a line segment containing 
$z$ in its relative interior and  for each $z \in S_1^1$, 
there exists a neighbourhood of $z$ that is contained in 
$Q_{z}^-$. Hence there exists $\beta_1 \in (0,1)$ and $\ve_1 > 0$ such that for every 
$C \in {\cal P}_1({\cal U})$, 
\[ H[s_1(1 - \beta_1), s_1 (1 + \beta_1)]  \cap 
 F_{C}  \subset Q_{z}^- \cap F_{C} 
\]
(where $z = H_{s_1} \cap F_C$) and 
\be 
\label{def-ej}
 H[s_1(1 - \beta_1), s_1 (1 + \beta_1)]  \cap 
 N_{\ve_1}(F_{C})  \subset Q_{z}^-.   
\ee
Taking the union of the penultimate relation over $C \in {\cal P}(\cu)$ and $z \in S_1^1$ yields 
(\ref{induction}) for $j = 1$.

Now suppose there exists $k \leq J-1$ such that for all $j \leq k -1$,  
there exist $\beta_j \in (0,1), \ve_j > 0$ and a finite set of points 
$S_1^j \subset H_{s_1} \cap [\cup_{C \in {\cal P}_j(\cu)} F_C]$ such that  
(\ref{induction})  holds and  for every $j \leq k-1$,  the relation 
\be
\label{ej}
  H[s_1(1 - \beta_j), s_1 (1 + \beta_j)]  \cap 
 \left[\cup_{C \in {\cal P}_j(\cu)} N_{\ve_j}(F_{C})\right]  \subset \ds \cup_{z \in \cup_{1 \leq j \leq k-1} S_1^j} Q_{z}^-   
\ee
is satisfied. 
We will show that then (\ref{induction}) also holds for $j = k$ and (\ref{ej}) holds with $k-1$ replaced by $k$.  
Indeed, the argument in this case is analogous to the case $j=1$.  
The  difference now is that it is no longer true that the intersection of 
the $k$-dimensional faces with $H_{s_1}$ are mutually disjoint and therefore
the direct analog of (\ref{choice}) does not hold. 
Indeed, the intersections of multiple $k$-dimensional faces yield lower-dimensional faces. 
But by the second induction assumption (\ref{ej}), we have already covered a neighbourhood of the intersection of
$H_{s_1}$ with  these lower-dimensional faces.  
 Thus once these neighbourhoods are excised from 
the $k$-dimensional faces, we once again obtain disjoint sets and an inequality analogous
to (\ref{choice}) will be true..  To make this reasoning precise, fix $C \in {\cal P}_k(\cu)$. 
By the definition of $r_z$, it follows that if 
\[ z \in \tilde{F}_C \doteq [F_C \cap H_{s_1}] 
\sm [ \cup_{1 \leq j \leq k-1} \cup_{\tilde{C} \in {\cal P}_j(\cu)} N_{\ve_j} (F_{\tilde{C}}) ], 
\]
then $r_z > \min_{1 \leq j < k} \ve_j$. 
Thus choose $0 < \rho^{(k)} < \min_{1 \leq j < k} \ve_j \wedge \kappa$ (note that this is the 
analogue of (\ref{choice}) that we wanted). 
Since $\cup_{C \in {\cal P}_k(\cu)} \tilde{F}_C$ is bounded and 
$0 \in (Q_C^-)^\circ$ for every $C$, 
it follows that there exist $\beta_k \in (0,1)$, $\ve_k > 0$ and
 a finite number of points $S_1^k = \cup_{C \in{\cal P}_k(\cu)} S_1^{k,C}$ with 
 $S_1^{k,C} \subset \tilde{F}_C$ such that if 
$\rho_z \doteq \rho^{(k)}$  and (as before)  
$Q_{z}^- = z + \rho_z Q_{I(z)}^-$ for $z \in S_1^k$, then 
(\ref{induction}) is satisfied with $j = k$ and (\ref{ej}) holds with $k-1$ replaced by $k$. 
By induction, this shows that (\ref{induction}) holds for all $j=1, \ldots, J-1$,  
and the proof of the first property of the claim is complete.  
The proof of the second property is analogous and therefore omitted.  
\ink \\

\noi
{\bf Acknowledgments.} The author would like to thank an anonymous referee for a careful reading 
of the paper and for making many constructive comments. 
The author is also grateful to Giovanni Leoni and Steve Shreve for a useful discussion that led to a simplification 
of the proof of Lemma \ref{lem-esp1},  and  to Marty Reiman for his interest in this work.

\bibliography{ref} 
\bibliographystyle{plain}

\end{document}